\documentclass[11pt,a4paper]{article}

\usepackage[T1]{fontenc}
\usepackage[utf8]{inputenc}
\usepackage{amsmath,amssymb,amsthm}
\usepackage[margin=2.6cm]{geometry}
\usepackage[numbers,sort&compress]{natbib}
\usepackage{booktabs,array}
\usepackage{tikz}
\usetikzlibrary{patterns}
\usepackage{hyperref}
\theoremstyle{plain}
\newtheorem{theorem}{Theorem}[section]
\newtheorem{proposition}[theorem]{Proposition}
\newtheorem{lemma}[theorem]{Lemma}
\newtheorem{corollary}[theorem]{Corollary}
\theoremstyle{definition}
\newtheorem{definition}[theorem]{Definition}
\newtheorem{example}[theorem]{Example}
\theoremstyle{remark}
\newtheorem{remark}[theorem]{Remark}

\newcommand{\yout}{y^{\mathrm{out}}}
\newcommand{\yin}{y^{\mathrm{in}}}
\newcounter{proofstep}
\newcommand{\LL}{\mathcal{L}}
\newcommand{\N}{\mathbb{N}}
\newcommand{\R}{\mathbb{R}}
\newcommand{\C}{\mathbb{C}}
\NewDocumentCommand{\abs}{s m}{\IfBooleanTF{#1}{\lvert #2\rvert}{\left\lvert #2\right\rvert}}
\NewDocumentCommand{\norm}{s m}{\IfBooleanTF{#1}{\lVert #2\rVert}{\left\lVert #2\right\rVert}}
\newcommand{\Id}{\mathrm{Id}}
\newcommand{\dist}{\operatorname{dist}}
\newcommand{\Ima}{\operatorname{Im}}
\newcommand\Phihat{\widehat{\Phi}}
\newcommand\clos[1]{\overline{#1}}
\newcommand\conj[1]{\overline{#1}}
\NewDocumentCommand{\ens}{s m}{\IfBooleanTF{#1}{\mathopen{\{}\begingroup\def\stdelim{}#2\endgroup\mathclose{\}}}{\left\{\def\stdelim{\middle}#2\right\}}}
\newcommand{\stdelim}{}
\newcommand{\st}{\;\stdelim\vert\;}
\newcommand{\Span}{\operatorname{span}}
\newcommand{\diag}{\operatorname{diag}}
\NewDocumentCommand{\pair}{s m m}{\IfBooleanTF{#1}{\langle #2,#3\rangle}{\left\langle #2,#3\right\rangle}}
\newcommand{\conca}[1]{\mathbin{\diamond_{\scriptscriptstyle #1}}}
\newcommand\les{\lesssim}
\newcommand\ges{\gtrsim}
\newcommand{\dive}{\mathrm{div} \,}
\newcommand\tr{{\mathsf{T}}}
\newcommand\indic{\mathsf{1}}
\DeclareMathOperator{\codim}{codim}

\renewcommand{\Re}{\operatorname{Re}}

\AtBeginDocument{\def\MR#1{}}

\providecommand{\keywords}[1]{%
  \par\addvspace{1.2ex}%
  {\small\noindent\textit{Key words.}\ #1\par}}
\providecommand{\subjclass}[2][2020]{%
  \par\addvspace{0.4ex}%
  {\small\noindent\textit{#1 Mathematics Subject Classification.}\ #2\par}%
  \addvspace{1.2ex}}

\title{Compactness--uniqueness method for null controllability}

\author{Guillaume Olive\thanks{Faculty of Mathematics and Computer Science, Jagiellonian
University, ul. {\L}ojasiewicza 6, 30-348 Krak\'{o}w, Poland.
E-mail: math.golive@gmail.com or guillaume.olive@uj.edu.pl}}

\date{November 10, 2026}

\begin{document}
\maketitle

\begin{abstract}
We develop the compactness--uniqueness method for the null controllability of abstract linear control systems on Banach spaces.
Our first main result characterizes the null controllable states for systems that are null controllable up to finite dimensional spaces, and shows that null controllability is then equivalent to the Fattorini--Hautus test.
Our second main result gives a dual counterpart of this assumption, an observability inequality with a suitable compact remainder.
These results unify and extend several known results, such as those on heat equations with nonlocal terms and on the simultaneous controllability of systems.
Finally, for a transport equation with nonlocal terms from population dynamics, we characterize the null controllable states after a sharp time, including in the non-reflexive space $L^1$.
\end{abstract}

\keywords{Null controllability; compactness--uniqueness; Fattorini--Hautus test; nonlocal equations; simultaneous controllability; age-structured population dynamics}

\subjclass[2020]{93B05; 93B07, 93C25, 93C20, 35R09, 92D25.}

\tableofcontents


\section{Introduction}\label{sec:intro}

The compactness--uniqueness method is one of the most effective tools to establish the controllability of linear partial differential equations, usually through observability inequalities.
One first proves the desired inequality up to a \emph{compact} remainder, by whatever means the equation allows (multipliers, or the known controllability of a simpler system), and one then removes this remainder by a \emph{unique continuation} property.
Early uses include \cite{RT74} for a decay problem and \cite{Zua87,Lio88} for controllability, and the method has since been applied to a large variety of systems, see e.g.\ \cite{Zua91,Kom89,BLR92,Ros97,Meh04,CT10,LR17,CDM23}.
For exact controllability, the method has a general abstract form \cite{DO18}, whose analytic core is a lemma of Peetre \cite{Pee61}: the uniqueness step reduces to the Fattorini--Hautus test \cite{Fat66}, a condition on the stationary adjoint problem which is necessary for any controllability property, and the reachable states are described exactly when this test fails.

For null controllability, the situation is different, and the reason is structural.
The null observability inequality only estimates the final state of the adjoint system, that is, a seminorm, to which the Peetre lemma does not directly apply.
Indeed, whereas a compact perturbation of the generator preserves exact controllability as soon as the perturbed system remains approximately controllable, an arbitrarily small one can destroy null controllability in a given time \cite[Section 4]{DR77}.
The method has thus been adapted case by case, each time through a device specific to the problem: a weak space with exponentially decaying weights, adapted to the heat equation \cite{FCLZ16}, a restriction to the reachable space of a system null controllable in arbitrary time \cite{ELT22}, an additional estimate of the adjoint state between two intermediate times \cite{CN21,HO22-JDE}, or, closest to the present work, a reduction to null controllability up to finite dimensional spaces followed by a unique continuation argument, a principle isolated in \cite[Proposition~A.1]{KL26} for parabolic--transport systems.

\paragraph{Main results.}
In this article, we give the general form of the compactness--uniqueness method for null controllability, for abstract systems $y'=Ay+Bu$ on Banach spaces, with no reflexivity, possibly unbounded control operators and general classes of controls.
Our first main result, Theorem~\ref{thm:main}, assumes that the system is \emph{null controllable up to finite dimensional spaces}: in some time $T_0$, every state of a subspace of finite codimension can be steered into a fixed finite dimensional subspace.
It then \emph{completely characterizes} the null controllable states, in every time $T>T_0$: they form the pre-annihilator of the space spanned by the unobservable generalized eigenvectors of $(A^*,B^*)$, which is finite dimensional.
In particular, the system is null controllable in some time $T>T_0$ if, and only if, it is in every such time, if, and only if, the Fattorini--Hautus test holds.
Since its assumption is also necessary for null controllability, Theorem~\ref{thm:main} reduces the problem exactly to two \emph{independent} questions: null controllability up to finite dimensional spaces, and the Fattorini--Hautus test.
This recovers \cite[Proposition~A.1]{KL26}, stated there in Hilbert spaces, and goes beyond it by describing the null controllable states when the test fails.
The proof uses neither compactness nor observability inequalities, and rests on one key idea: to \emph{decouple} the duration of the free evolution from the length of the control window.
This is what allows us to dispense with the additional estimate of \cite[Theorem 7.7]{HO22-JDE} (Remark~\ref{rem:HO22}).

Our second main result, Theorem~\ref{thm:obs}, expresses this assumption on the adjoint system, which is the form often used in practice: an observability inequality with a compact remainder $KLz$, in which both the final state of the adjoint system and $Lz$ are estimated.
When the control space is reflexive, this is equivalent to the assumption of Theorem~\ref{thm:main}, so that it identifies exactly the compact remainders that can be tolerated when only a seminorm is estimated.
The Peetre lemma enters here, applied to $L$ augmented by the observation, and not to the seminorm.
The two halves of the method thus become two separate statements, compactness in Theorem~\ref{thm:obs} and uniqueness in Theorem~\ref{thm:main}.

\paragraph{Applications.}
These two results cover in one stroke situations treated so far by different routes, and go beyond them (Section~\ref{sec:comparison}).
For the heat equation with a nonlocal term \cite{FCLZ16}, they give the null controllability under the Fattorini--Hautus test instead of the analyticity of the kernel, as also obtained in \cite{ELT22} by a different route, and they describe the null controllable states when the test fails, which happens for simple kernels.
For two systems $(A_1,B_1)$ and $(A_2,B_2)$ driven by the same control, the first one steered to an arbitrary target and the second one to rest \cite{Tuc26}, they replace the simultaneous approximate observability, obtained there by analyticity, by the test, which holds as soon as $A_1^*$ and $A_2^*$ share no eigenvalue.
This settles one of the questions left open in \cite[Section 7]{Tuc26} and partly answers another one.
Finally, we treat a linearized model of age-structured population dynamics \cite{GM74}, a transport equation with nonlocal terms in the equation and at the boundary, controlled from an arbitrary band (Section~\ref{sec:nonloc t}).
We completely describe the null controllable states after an explicit time, which is the \emph{minimal null control time} when the Fattorini--Hautus test holds, in every space $L^r$, $1\le r<\infty$.
This includes the natural space $L^1$, which is not reflexive, so that the observability route of Theorem~\ref{thm:obs} is not available, and Theorem~\ref{thm:main}, a statement about controls only, is applied directly.

\paragraph{Organization of the paper.}
The rest of this section fixes the notation and the framework, and recalls the known abstract results on which the present work builds (Section~\ref{sec:CU}).
Section~\ref{sec:mainthm} states Theorems~\ref{thm:main} and~\ref{thm:obs}, which are proved in Sections~\ref{sec:proofmain} and~\ref{sec:proofobs}, respectively.
Section~\ref{sec:comparison} is devoted to the applications.

\subsection{Notation}
\label{sec:nota}

\paragraph{Spaces and operators.}
Unless otherwise stated, all the spaces occurring below are complex.
For a normed space $E$, we denote by $E'$ its dual, that is, the space of continuous conjugate linear forms on $E$, equipped with the norm $\norm*{f}_{E'}:=\sup\ens*{|\pair*{f}{e}_{E',E}| \st \norm*{e}_{E}\le1}$, where $\pair*{f}{e}_{E',E}:=f(e)$ is linear in $f$ and conjugate linear in $e$.
Recall that $E'$ is a Banach space, even if $E$ is only a normed space, the scalar field being complete.
A Hilbert space is identified with its dual through the Riesz theorem, this identification being linear.
We denote by $\LL(E,F)$ the space of bounded linear operators from $E$ to $F$, with $\LL(E):=\LL(E,E)$, and by $L^*\in\LL(F',E')$ the adjoint of $L\in\LL(E,F)$, defined by $\pair*{L^*f}{e}_{E',E}:=\pair*{f}{Le}_{F',F}$, so that $(\lambda L)^*=\conj{\lambda}L^*$ for $\lambda\in\C$.
Similarly, if $G:D(G)\subset E\to E$ is densely defined, then $f\in E'$ belongs to $D(G^*)$ if, and only if, $y\mapsto\pair*{f}{Gy}_{E',E}$ is continuous on $D(G)$ for the norm of $E$, and then $\pair*{G^*f}{y}_{E',E}=\pair*{f}{Gy}_{E',E}$ for every $y\in D(G)$.
The domain of a closed operator $G$ always carries the graph norm $\norm*{y}_{D(G)}:=\norm*{y}_{E}+\norm*{Gy}_{E}$, for which it is a Banach space.
We write $\Ima$ and $\ker$ for the range and the kernel of an operator, $R^{-1}M:=\ens*{x\st Rx\in M}$ for the preimage of a set $M$ under a linear map $R$, no injectivity being implied, $\Span M$ for the linear span of $M$, and $\dist_{E}(x,M):=\inf_{a\in M}\norm*{x-a}_{E}$ for $x\in E$ and $M\subset E$.
We denote by $\Id$ the identity operator and, for a linear operator $G$, by $\rho(G)$, $\sigma(G)$ and $\sigma_p(G)$ its resolvent set, its spectrum and the set of its eigenvalues.

\paragraph{Finite codimension.}
For a subspace $M\subset E$, we write $\codim M:=\dim(E/M)$.
Thus, $M$ is of finite codimension if, and only if, $E=M\oplus M_0$ for some finite dimensional subspace $M_0$.

\paragraph{Annihilators.}
For subsets $M\subset E$ and $N\subset E'$, we denote by
\[
M^{\perp}:=\ens{f\in E'\st \pair{f}{e}_{E',E}=0,\quad \forall e\in M},
\qquad
N_{\perp}:=\ens{e\in E\st \pair{f}{e}_{E',E}=0,\quad \forall f\in N},
\]
the annihilator of $M$ and the pre-annihilator of $N$, which are closed subspaces of $E'$ and $E$, respectively.
We have (see e.g.\ \cite[Proposition 1.9]{Bre11})
\begin{equation}\label{prop perp}
(M^\perp)_\perp=\clos{\Span M}.
\end{equation}
For subspaces $V,W\subset E$ and $L\in\LL(E)$, with $W$ closed, we have
\begin{equation}\label{eq:dualincl}
LV\subset W\quad\Longleftrightarrow\quad L^*W^{\perp}\subset V^{\perp}.
\end{equation}
Indeed, for $f\in W^{\perp}$ and $v\in V$, we have $\pair*{L^*f}{v}_{E',E}=\pair*{f}{Lv}_{E',E}$.
If $LV\subset W$, this vanishes.
Conversely, if $L^*W^{\perp}\subset V^{\perp}$, this shows that $Lv\in(W^{\perp})_{\perp}=W$, by \eqref{prop perp}, $W$ being closed.
Finally, if $M\subset E$ is a closed subspace of finite codimension, then $M^{\perp}$ is finite dimensional, with (see e.g.\ \cite[Proposition~11.13(b)]{Bre11})
\begin{equation}\label{eq:dimMperpcodim}
\dim M^{\perp}=\codim M.
\end{equation}

\paragraph{Other notation.}
The set $\N$ of natural numbers contains $0$.
We denote by $\indic_{\omega}$ the indicator function of a set $\omega$, by $|I|$ the length of an interval $I$, by $s_+:=\max\ens*{s,0}$ the positive part of $s\in\R$, by $\C_-$ and $\C_+$ the open half-planes $\ens*{\Re\lambda<0}$ and $\ens*{\Re\lambda>0}$, by $v^\tr$ the transpose of a vector or a matrix $v$, and by $\partial_\nu$ the outward normal derivative.
Finally, $f\les g$ means that $f\le cg$ for some constant $c>0$ which does not depend on the variables involved (typically, the adjoint state $z$), $f\ges g$ means $g\les f$, and $f\simeq g$ means $f\les g$ and $g\les f$.
The numbers $c,c',\ldots>0$ have the same meaning.

\subsection{Framework}\label{sec:setting}

\subsubsection{Solutions by transposition}\label{sec:control}

Let $Y$ and $U$ be complex Banach spaces, let $A:D(A)\subset Y\to Y$ be the generator of a $C_0$-semigroup $(S(t))_{t\ge0}$ on $Y$, and let $B\in\LL(U,D(A^*)')$ be the \emph{control operator}.
We consider the control system
\begin{equation}\label{syst}
\left\{
\begin{aligned}
\frac{d}{dt}y(t) &=Ay(t)+Bu(t), &&t\in(0,T),\\
y(0) &=y^0,
\end{aligned}
\right.
\end{equation}
where $T>0$, $y(t)$ is the state at time $t$, $y^0$ the initial data and $u(t)$ the \emph{control} at time $t$.
A reader interested only in the Hilbert setting, for which we refer to \cite{Cor07,TW09}, may read this section with $Y$ and $U$ Hilbert spaces identified with their duals and controls in $U_T=L^2(0,T;U)$.
The framework below is detailed with some care because the adjoint semigroup $(S(t)^*)_{t\ge0}$ is, in general, not strongly continuous on $Y'$ \cite[Paragraph~I.5.14]{EN00}.

\paragraph{The space $D(A^*)'$.}
The space $Y$ can always be identified with a subspace of $D(A^*)'$.
Indeed, every $y\in Y$ defines a continuous conjugate linear form on $D(A^*)$ via the map
\begin{equation}\label{eq:iota}
j:Y\longrightarrow D(A^*)' ,\qquad
\pair{j y}{z}_{D(A^*)',D(A^*)}:=\conj{\pair{z}{y}}_{Y',Y} ,
\end{equation}
which is linear, continuous and injective, since $(\lambda-A^*)^{-1}$ is an isomorphism from $Y'$ onto $D(A^*)$ for $\lambda\in\rho(A^*)$.
We identify $Y$ with $j(Y)$ accordingly, writing $Y\subset D(A^*)'$.

This allows us to extend $S(t)$ to $D(A^*)'$, still denoted by $S(t)$, by setting
\[
\pair{S(t)\phi}{z}_{D(A^*)',D(A^*)}:=\pair{\phi}{S(t)^*z}_{D(A^*)',D(A^*)}.
\]
Note, however, that the extended semigroup is not strongly continuous on $D(A^*)'$ in general.

Finally, we use the same notation for the adjoint $B^*\in\LL(D(A^*)'',U')$ and for its restriction to $D(A^*)$, which is naturally identified with a subspace of its bidual $D(A^*)''$.

\paragraph{Notion of solution.}
Let $t\in(0,T]$ and $z\in D(A^*)$, and let $\varphi(s):=S(t-s)^*z$ for $s\in[0,t]$, which solves, in the weak$^*$ sense, the adjoint system $-\varphi'=A^*\varphi$ in $(0,t)$, $\varphi(t)=z$.
Pairing \eqref{syst} with $\varphi(s)$ and integrating by parts on $(0,t)$, we formally obtain
\begin{equation}\label{eq:weaksol}
\pair{y(t)}{z}_{D(A^*)',D(A^*)}-\pair{y^0}{S(t)^*z}_{D(A^*)',D(A^*)}=\int_0^t\conj{\pair{B^*S(t-s)^*z}{u(s)}}_{U',U}\,ds .
\end{equation}
Its right-hand side requires the integrand to be measurable, which is not automatic.
We shall thus assume throughout that, for every $z\in D(A^*)$, the function
\begin{equation}\label{hyp measurable}
t\in(0,+\infty)\longmapsto B^*S(t)^*z\in U'\quad\text{is weakly$^*$ measurable,}
\end{equation}
that is, $t\mapsto\pair*{B^*S(t)^*z}{w}_{U',U}$ is measurable for every $w\in U$.
This holds for instance if $B\in\LL(U,Y)$ or if $Y$ is reflexive, but neither property will be part of our standing assumptions.
Under \eqref{hyp measurable}, the right-hand side of \eqref{eq:weaksol} is indeed integrable:
\begin{equation}\label{RGEbasic}
\int_0^t \left|\pair{B^*S(t-s)^*z}{u(s)}_{U',U}\right|\,ds
\le c\norm{z}_{D(A^*)}\norm{u}_{L^1(0,t;U)},
\end{equation}
with $c=\sup_{s \in [0,T]} \norm{S(s)}_{\LL(Y)} \norm{B}_{\LL(U,D(A^*)')}$.
This estimate also shows that the right-hand side of \eqref{eq:weaksol} is a continuous conjugate linear form of $z\in D(A^*)$.

\begin{definition}[solution by transposition]\label{def:solution}
Let $T>0$, $y^0\in Y$ and $u\in L^1(0,T;U)$.
The \emph{solution} to \eqref{syst} is the function $y:[0,T]\to D(A^*)'$ defined by \eqref{eq:weaksol} for every $t\in[0,T]$ and $z\in D(A^*)$.
\end{definition}

The key controllability properties are encoded in the map defined by the right-hand side of \eqref{eq:weaksol}.

\begin{definition}[control map]
The \emph{control map} of the pair $(A,B)$ is the map defined, for $T>0$ and $u\in L^1(0,T;U)$, by
\begin{equation}\label{eq:inputmapdef}
\pair{\Phi_T u}{z}_{D(A^*)',D(A^*)}
:=\int_0^T\conj{\pair{B^*S(T-t)^*z}{u(t)}}_{U',U}\,dt, \qquad z\in D(A^*),
\end{equation}
and we set $\Phi_0:=0$.
\end{definition}

We have $\Phi_T\in\LL(L^1(0,T;U),D(A^*)')$ by \eqref{RGEbasic}.
Comparing with \eqref{eq:weaksol}, and by definition of the extension of $S(t)$ to $D(A^*)'$, the solution to \eqref{syst} is given by the \emph{Duhamel formula}
\begin{equation}\label{duhamel}
y(t)=S(t)y^0+\Phi_t u_{|(0,t)},\qquad t\in[0,T].
\end{equation}

\paragraph{The composition property.}
For $T>0$ and $t\in(0,T)$, the \emph{concatenation} at time $t$ of $u\in L^1(0,t;U)$ and $v\in L^1(0,T-t;U)$ is $u\conca{t}v\in L^1(0,T;U)$, defined by
\[
\left(u\conca{t}v\right)(s):=
\begin{cases}
u(s) & \text{if }s\in(0,t),\\
v(s-t) & \text{if }s\in(t,T).
\end{cases}
\]
Then, we have the \emph{composition property}
\begin{equation}\label{eq:composition}
\Phi_{T}\left(u\conca{t}v\right)=S(T-t)\,\Phi_{t}u+\Phi_{T-t}v.
\end{equation}
This follows by splitting the integral \eqref{eq:inputmapdef} at $t$ and a change of variables on $(t,T)$.

\subsubsection{Control classes and the regularity-gain condition}\label{sec:scale}

By \eqref{eq:inputmapdef}, the control map is defined on the whole of $L^1(0,T;U)$ and takes its values in $D(A^*)'$.
Applications usually afford better: taking the controls in a smaller space improves the space in which the state is found.

\begin{definition}[control class]\label{def:controlclass}
A \emph{control class} is a family $(U_T)_{T>0}$ of Banach spaces, each $U_T$ being a dense subspace of $L^1(0,T;U)$ with continuous inclusion, such that
\begin{equation}\label{eq:ALCSconcat}
U_t\conca{t}U_{T-t}\subset U_T,
\end{equation}
for every $T>0$ and $t\in(0,T)$.
\end{definition}

\begin{example}[Lebesgue and Sobolev controls]\label{ex:UT}
\emph{Lebesgue.}
For any $p \in [1,\infty]$, the family $U_T:=L^p(0,T;U)$ is a control class, and the one most often used.

\emph{Sobolev.}
Let $\omega\subset\R^N$ be open, $U:=L^2(\omega)$ identified with its dual, $\omega_T:=(0,T)\times\omega$, $k\in\N$, and $U_T:=H^k_0(\omega_T)$, the closure of $C_c^\infty(\omega_T)$ in $H^k(\omega_T)$.
It is continuously included in $L^2(\omega_T)=L^2(0,T;U)$, hence in $L^1(0,T;U)$, and it is dense in $L^1(0,T;U)$, since it contains $C_c^\infty(\omega_T)$.
For \eqref{eq:ALCSconcat}, the extension by zero from $\omega_t$ to $\omega_T$ and the shift $(s,x)\mapsto(s+t,x)$ from $\omega_{T-t}$ to $\omega_T$ map $C_c^\infty$ into $C_c^\infty$ and preserve the norm of $H^k$, hence extend by density to isometries into $H^k_0(\omega_T)$, which are still the extension by zero and the shift, as limits in $L^2$.
Products of spaces of this type are used in \cite[proof of Proposition 3.1]{KL26}.
\end{example}

\begin{definition}[regularity-gain condition]\label{def:wellposed}
Let $(U_T)_{T>0}$ be a control class.
We say that the \emph{regularity-gain condition} holds if
\begin{equation*}\label{eq:RGC}\tag{RGC}
\Phi_T\left(U_T\right)\subset Y, \quad \forall T>0.
\end{equation*}
\end{definition}

The inclusion means that $\Phi_T(U_T)\subset j(Y)$, with $j$ defined in \eqref{eq:iota}.
By the Duhamel formula \eqref{duhamel}, condition \eqref{eq:RGC} says that a state starting in $Y$ is still in $Y$ at time $T$ when driven by a control in $U_T$.
This pure range condition is turned into an estimate by the closed graph theorem:
\[
{\Phi_T}_{|U_T}\in\LL\left(U_T,Y\right).
\]
Indeed, the graph of ${\Phi_T}_{|U_T}$ is closed in $U_T\times Y$, since $\Phi_T\in\LL(L^1(0,T;U),D(A^*)')$, the inclusion of $U_T$ into $L^1(0,T;U)$ is continuous, and $j$ is continuous and injective.
From now on, $\Phi_T$ denotes this restriction, and we write $\Ima\Phi_T:=\Phi_T(U_T)$, a subspace of $Y$, and $\Phi_T^*\in\LL(Y',U_T')$ for its adjoint.
By \eqref{eq:iota}, the definition \eqref{eq:inputmapdef} then reads, in the duality between $Y'$ and $Y$,
\begin{equation}\label{eq:link}
\pair{z}{\Phi_Tu}_{Y',Y}=\int_0^T\pair{B^*S(T-t)^*z}{u(t)}_{U',U}\,dt ,\qquad\forall z\in D(A^*),\ \forall u\in U_T .
\end{equation}
Finally, for $u\in U_t$ and $v\in U_{T-t}$, we have $u\conca{t}v\in U_T$ by \eqref{eq:ALCSconcat}, and all the terms of the composition property \eqref{eq:composition} belong to $Y$.

\begin{remark}[bounded control operator]
\label{rem:bounded B}
If $B$ is bounded, then \eqref{hyp measurable} holds, as mentioned in Section~\ref{sec:control}, and, for every $u\in U_T$, the Bochner integral $\int_0^TS(T-t)Bu(t)\,dt$ belongs to $Y$ and satisfies $\Phi_Tu=j(\int_0^TS(T-t)Bu(t)\,dt)$, so that \eqref{eq:RGC} holds.
In that case, the solution $y$ to \eqref{syst} belongs to $C^0([0,T];Y)$, by the dominated convergence theorem.
\end{remark}

\subsubsection{Controllability}\label{sec:controllability}

Let $(U_T)_{T>0}$ be a control class, and assume that the regularity-gain condition \eqref{eq:RGC} holds.

\begin{definition}[controllability]\label{def:control}
For $T>0$, we say that the system \eqref{syst} is:
\begin{itemize}
\item \emph{exactly controllable in time $T$} if, for every $y^0,y^T\in Y$, there exists $u\in U_T$ such that the corresponding solution to \eqref{syst} satisfies $y(T)=y^T$; equivalently, $\Ima\Phi_T=Y$.

\item \emph{null controllable in time $T$} if the above property holds at least for $y^T=0$; equivalently, $\Ima S(T) \subset \Ima\Phi_T$.
\end{itemize}
The \emph{null controllable states} in time $T$ are the elements of the preimage
\[
S(T)^{-1}\,\Ima\Phi_T=\ens{y^0\in Y\st \exists u\in U_T,\quad y(T)=0}.
\]
\end{definition}

\begin{remark}[concatenation and controllability]
\label{rem:concat for NC}
For every $t,T>0$, we have
\begin{align}
\label{eq:monoP1}
S(t)\,\Ima\Phi_T&\subset\Ima\Phi_{T+t} ,\\
\label{eq:monoP2}
\Ima\Phi_T&\subset\Ima\Phi_{T+t} .
\end{align}
Indeed, let $u\in U_T$.
By \eqref{eq:ALCSconcat}, the extension by zero $u\conca{T}0$ and the forward shift $0\conca{t}u$ belong to $U_{T+t}$, and the composition property \eqref{eq:composition} gives $\Phi_{T+t}\left(u\conca{T}0\right)=S(t)\,\Phi_T u$ and $\Phi_{T+t}\left(0\conca{t}u\right)=\Phi_T u$.
In particular, if $S(T)y^0\in\Ima\Phi_T$, then $S(T+t)y^0\in\Ima\Phi_{T+t}$ by \eqref{eq:monoP1}, so that null controllability in some time implies null controllability at larger times.
\end{remark}

\subsubsection{Dual conditions}\label{sec:duality}

\paragraph{The Fattorini--Hautus test.}
The following condition, which goes back to \cite{Fat66}, plays a central role in this article.

\begin{definition}[Fattorini--Hautus test]\label{def:FH}
Let $Z$ and $W$ be complex vector spaces, and let $G:D(G)\subset Z\to Z$ and $C:D(G)\to W$ be linear operators.
We say that the pair $(G,C)$ satisfies the \emph{Fattorini--Hautus test} if
\begin{equation*}\label{eq:FH}\tag{FH}
\ker(\lambda-G)\cap\ker C=\ens{0},\qquad\forall\lambda\in\C .
\end{equation*}
\end{definition}

\begin{proposition}[necessity of the test]\label{prop:FHnec}
If the system \eqref{syst} is null controllable in some time $T>0$, then $(A^*,B^*)$ satisfies the Fattorini--Hautus test.
\end{proposition}

\begin{proof}
Let $\lambda\in\C$ and $z\in\ker(\lambda-A^*)\cap\ker B^*$ with $z\neq0$.
We first show that $S(t)^*z=e^{\lambda t}z$ for every $t\ge0$.
Let $t>0$ and $y\in D(A)$.
The function $s\mapsto e^{\lambda s}\pair*{z}{S(t-s)y}_{Y',Y}$ is differentiable on $[0,t]$, with derivative $e^{\lambda s}\pair*{(\lambda-A^*)z}{S(t-s)y}_{Y',Y}=0$, so that its values at $s=0$ and $s=t$ coincide, that is, $\pair*{S(t)^*z}{y}_{Y',Y}=\pair*{e^{\lambda t}z}{y}_{Y',Y}$, and we conclude by density of $D(A)$.
Now, let $y^0\in Y$ be such that $\pair*{z}{y^0}_{Y',Y}\neq0$, which exists since $z\neq0$.
For every $u\in U_T$, since $B^*S(T-t)^*z=e^{\lambda(T-t)}B^*z=0$, \eqref{eq:link} gives $\pair*{z}{\Phi_Tu}_{Y',Y}=0$, so that $\pair*{z}{S(T)y^0+\Phi_Tu}_{Y',Y}=\pair*{S(T)^*z}{y^0}_{Y',Y}=e^{\lambda T}\pair*{z}{y^0}_{Y',Y}\neq0$.
Thus, $y(T)\neq0$ for every $u\in U_T$.
\end{proof}

\paragraph{Observability inequalities.}
Controllability properties are properties of ranges, and they are characterized by estimates on the adjoint operators, which is what we call \emph{duality}.

We start with exact controllability, which is the surjectivity of $\Phi_T$.
A bounded operator between Banach spaces is onto if, and only if, its adjoint is bounded below \cite[Theorem~2.20]{Bre11}.
Thus, the system \eqref{syst} is exactly controllable in time $T$ if, and only if, the following \emph{exact observability inequality} holds:
\begin{equation}\label{eq:exactobs}
\norm{z}_{Y'}\les \norm{\Phi_T^* z}_{U_T'} ,\qquad\forall z \in Y'.
\end{equation}

Null controllability is an inclusion between ranges, and such inclusions are characterized by the following classical result, which goes back to \cite{Dou66,DR77}.

\begin{theorem}[duality]
\label{thm:duality}
Let $X_1,X_2,X_3$ be Banach spaces, and let $\Lambda_1\in\LL(X_1,X_3)$ and $\Lambda_2\in\LL(X_2,X_3)$.
Consider the following two properties:
\begin{enumerate}
\item\label{classical i1}
$\Ima\Lambda_1\subset\Ima\Lambda_2$.

\item\label{classical i2}
There exists $c>0$ such that $\norm{\Lambda_1^*z}_{X_1'}\le c\,\norm{\Lambda_2^*z}_{X_2'}$ for every $z\in X_3'$.
\end{enumerate}
Then \ref{classical i1} implies \ref{classical i2}.
If moreover $X_2$ is reflexive, then \ref{classical i2} implies \ref{classical i1}, so that the two properties are equivalent.
\end{theorem}

See \cite[Part~IV, \S2.1]{Zab08} for a proof.
We apply this with $\Lambda_1=S(T)$, $\Lambda_2=\Phi_T$, $X_1=X_3=Y$ and $X_2=U_T$.
The null controllability of the system \eqref{syst} in time $T$ therefore always implies the \emph{null observability inequality}
\begin{equation}\label{eq:nullobs}
\norm{S(T)^*z}_{Y'}\les \norm{\Phi_T^*z}_{U_T'} ,\qquad\forall z\in Y',
\end{equation}
and it follows from it as soon as $U_T$ is reflexive.

When $Y$ and $U$ are Hilbert spaces, $U$ being identified with its dual, and $U_T=L^2(0,T;U)$, the function $t\mapsto B^*S(T-t)^*z$ is continuous for $z\in D(A^*)$ since the adjoint semigroup is then strongly continuous, and \eqref{eq:link} gives the familiar norm
\[
\norm{\Phi_T^*z}_{U_T'}=\left(\int_0^T\norm{B^*S(T-t)^*z}_{U}^2\,dt\right)^{1/2},\qquad\forall z\in D(A^*),
\]
by the Cauchy--Schwarz inequality and the choice $u(t)=B^*S(T-t)^*z$.

\subsection{Known abstract compactness--uniqueness results}\label{sec:CU}

\paragraph{Exact controllability.}
We end this section with the form of the compactness--uniqueness method obtained in \cite{DO18}, whose counterpart for null controllability is the subject of the present article.
In this paragraph, we place ourselves in the setting of \cite{DO18}: $Y$ and $U$ are Hilbert spaces, identified with their duals, $U_T=L^2(0,T;U)$, and \eqref{eq:RGC} holds.

The analytic heart of the method is the following lemma of Peetre \cite[Lemma 3]{Pee61}.

\begin{lemma}[Peetre]\label{lem:peetre}
Let $X_1,X_2$ be Banach spaces and let $\Lambda\in\LL(X_1,X_2)$.
The following two properties are equivalent:
\begin{enumerate}
\renewcommand{\theenumi}{(\alph{enumi})}\renewcommand{\labelenumi}{\theenumi}
\item\label{peetre i1}
there exist a compact operator $K\in\LL(X_1,X_3)$ ($X_3$ normed space) and $c>0$ such that
\begin{equation}\label{peetre estim}
\norm*{z}_{X_1}\le c\left(\norm*{\Lambda z}_{X_2}+\norm*{Kz}_{X_3}\right),\qquad\forall z\in X_1;
\end{equation}

\item\label{peetre i2}
$\ker\Lambda$ is finite dimensional and $\Ima\Lambda$ is closed.
\end{enumerate}
\end{lemma}

In particular, the compact remainder \emph{can be removed} in \eqref{peetre estim} if (and only if) $\ker\Lambda=\ens*{0}$, by the bounded inverse theorem.
To be precise, Peetre assumes that $K$ is the injection of $X_1$ into $X_3$, and that all the spaces are reflexive, but his proof applies word for word and uses no reflexivity.

This lemma is the core of \cite{DO18}, where it is applied as follows.
The method compares two systems $(A_1,B_1)$ and $(A_2,B_2)$, with control maps $\Phi_{1,T}$ and $\Phi_{2,T}$, the first one being exactly controllable in some time $T_0$.
The exact observability inequality \eqref{eq:exactobs} for the first system gives
\[
\norm{z}_{Y} \les \norm{\Phi_{2,T_0}^*z}_{L^2(0,T_0;U)}+\norm{\Phi_{1,T_0}^*z -\Phi_{2,T_0}^*z}_{L^2(0,T_0;U)},
\qquad\forall z\in Y,
\]
an observability inequality for the second system \emph{with a remainder} $K:=\Phi_{1,T_0}^*-\Phi_{2,T_0}^*$.
If $K$ is compact, the Peetre lemma (Lemma~\ref{lem:peetre}) shows that the remainder can be removed if, and only if, $\ker\Phi_{2,T_0}^*=\ens*{0}$.
By \eqref{prop perp}, this is the \emph{approximate controllability} of the second system in time $T_0$, that is, $\clos{\Ima\Phi_{2,T_0}}=Y$, a unique continuation property which may be as hard to establish as the controllability one started from.
These are the two halves of the method: \emph{compactness} of the remainder, and \emph{uniqueness} for the system of interest.
It was shown in \cite{DO18} that the uniqueness half can be weakened to the Fattorini--Hautus test, at the price of an arbitrarily small loss in the control time.
We recall this result in the form of \cite[Theorem 4.1]{DO18}, for a single system, after introducing the objects involved in its conclusion.

\begin{definition}[unobservable eigenvalues and generalized eigenvectors]\label{def:unobs}
Let $G$ and $C$ be as in Definition~\ref{def:FH}.
We define
\[
\sigma_F:=\ens{\lambda\in\C\st \ker(\lambda-G)\cap\ker C\ne\ens{0}},
\]
the set of \emph{unobservable eigenvalues} of the pair $(G,C)$ and, for $\lambda\in\C$,
\[
E_\lambda:=\ens{z\in\ker(\lambda-G)^m\st m\in\N,\quad C(\lambda-G)^jz=0,\quad\forall j<m},
\]
the subspace of \emph{unobservable generalized eigenvectors} of $(G,C)$ associated with $\lambda$.
\end{definition}

By definition, \eqref{eq:FH} holds if, and only if, $\sigma_F=\emptyset$.
The set $E_\lambda$ is a subspace, invariant under $G$.
Moreover, $E_\lambda\ne\{0\}$ if, and only if, $\lambda\in\sigma_F$.
Finally, the subspaces $E_\lambda$, $\lambda\in\C$, are in direct sum, since generalized eigenvectors associated with distinct eigenvalues are linearly independent.
Consequently, the space spanned by the unobservable generalized eigenvectors of $(G,C)$ is $\sum_{\lambda\in\C}E_\lambda=\bigoplus_{\lambda\in\sigma_F}E_\lambda$, since $E_\lambda=\ens{0}$ for $\lambda\notin\sigma_F$.
We can now state the result of \cite{DO18}.

\begin{theorem}[{\cite[Theorem 4.1]{DO18}}]
\label{thm:DO18control}
Let $(A,B)$ be a system as above and let $T_0>0$.
Assume that there exist a compact operator $K:Y\to V$, $V$ a normed space, and $c>0$ such that
\[
\norm*{z}_{Y}\le c\left(\norm*{\Phi_{T_0}^*z}_{L^2(0,T_0;U)}+\norm*{Kz}_{V}\right),\qquad\forall z\in Y.
\]
Then, we have:
\begin{enumerate}
\item\label{DO18:i0}
For every $T\ge T_0$, the set $\Ima\Phi_T$ is a closed subspace of $Y$ of finite codimension.

\item\label{DO18:i1}
The set $\sigma_F$ of unobservable eigenvalues of $(A^*,B^*)$ is finite and, for every $\lambda\in\sigma_F$, the subspace $E_\lambda$ of unobservable generalized eigenvectors of $(A^*,B^*)$ associated with $\lambda$ is finite dimensional.

\item\label{DO18:i2}
The sets of reachable states are given by
\[
\Ima\Phi_{T}=E_\perp,\qquad\forall T>T_0,\qquad\text{where }E:=\bigoplus_{\lambda\in\sigma_F}E_\lambda.
\]
In particular, if $(A^*,B^*)$ satisfies the Fattorini--Hautus test, then the system $(A,B)$ is exactly controllable in every time $T>T_0$.
\end{enumerate}
\end{theorem}

In finite dimension, the assumption of Theorem~\ref{thm:DO18control} always holds, with $K=\Id$ and any $T_0>0$, and the theorem reduces to the classical description of the reachable states $\Ima\Phi_T=\Ima\begin{pmatrix}B&AB&\cdots&A^{n-1}B\end{pmatrix}$ \cite[Corollary~1.4.9]{TW09}, where $n:=\dim Y$.

\paragraph{Null controllability.}
For this property, an abstract principle of the same kind was isolated in \cite{KL26}.
It is also one of the starting points of the present work.

\begin{proposition}[{\cite[Proposition~A.1]{KL26}}]\label{prop:KL26}
Let $T_0>0$.
Let $Y$ be a complex Hilbert space and let $A$ be the generator of a $C_0$-semigroup $(S(t))_{t\ge0}$ on $Y$.
Let $U$ be a complex Hilbert space and let $B\in\LL(U,Y)$.
For every $T>0$, let $U_T$ be a Hilbert space, continuously and densely embedded in $L^2(0,T;U)$, that satisfies the following ``extension by zero property'': if $u\in U_T$ and $a,b\ge0$, then the function equal to $0$ on $(0,a)$, to $u(\cdot-a)$ on $(a,a+T)$ and to $0$ on $(a+T,a+T+b)$ belongs to $U_{T+a+b}$.
Assume that there exist a finite dimensional subspace $F\subset Y$ invariant under the semigroup and a closed subspace $G\subset Y$ of finite codimension such that:
\begin{itemize}
\item
(control up to finite dimension) for every $y^0\in G$, there exists $u\in U_{T_0}$ such that the solution $y$ of $y'=Ay+Bu$, $y(0)=y^0$, satisfies $y(T_0)\in F$,

\item
(unique continuation) for every $\varepsilon>0$ and every finite linear combination $z$ of generalized eigenvectors of $A^*$, we have $B^*S(t)^*z=0$ for every $t\in(0,\varepsilon)$ $\Longrightarrow$ $z=0$.
\end{itemize}
Then, for every $T>T_0$ and every $y^0\in Y$, there exists $u\in U_T$ such that the solution $y$ of $y'=Ay+Bu$, $y(0)=y^0$, satisfies $y(T)=0$.
\end{proposition}


\section{Main results}\label{sec:mainthm}

This section states our main results.
The first one, Theorem~\ref{thm:main}, completely characterizes the set of null controllable states, assuming null controllability up to finite dimensional spaces.
The second one, Theorem~\ref{thm:obs}, gives this assumption in terms of observability, which is the form in which it is often verified in applications.

\subsection{Null controllability up to finite dimensional spaces}\label{sec:controlside}

The first main result of this article is the following.

\begin{theorem}[null controllable states]\label{thm:main}
Let $(U_T)_{T>0}$ be a control class, and assume that the regularity-gain condition \eqref{eq:RGC} holds.
Assume moreover that there exist $T_0>0$, a subspace $G\subset Y$ of finite codimension and a finite dimensional subspace $F\subset Y$ such that
\begin{equation}\label{eq:primalH}
\forall y^0\in G,\ \exists u\in U_{T_0},\qquad S(T_0)y^0+\Phi_{T_0}u\in F .
\end{equation}
Then, we have:
\begin{enumerate}
\item\label{main:i1}
For every $T\ge T_0$, the set $S(T)^{-1}\,\Ima\Phi_T$ is a closed subspace of $Y$ of finite codimension.

\item\label{main:i2}
The set $\sigma_F$ of unobservable eigenvalues of $(A^*,B^*)$ is finite and, for every $\lambda\in\sigma_F$, the subspace $E_\lambda$ of unobservable generalized eigenvectors of $(A^*,B^*)$ associated with $\lambda$ is finite dimensional.

\item\label{main:i3}
The sets of null controllable states are given by
\begin{equation}\label{eq:reach}
S(T)^{-1}\,\Ima\Phi_{T}=E_\perp,\qquad\forall T>T_0,\qquad\text{where }E:=\bigoplus_{\lambda\in\sigma_F}E_\lambda.
\end{equation}
In particular, if $(A^*,B^*)$ satisfies the Fattorini--Hautus test, then the system $(A,B)$ is null controllable in every time $T>T_0$.
\end{enumerate}
\end{theorem}

\begin{remark}[comparison with \cite{KL26}]\label{rem:scope}
Theorem~\ref{thm:main} extends the result of \cite{KL26}, recalled in Proposition~\ref{prop:KL26}, in several directions.
\begin{itemize}
\item
\emph{Null controllable states.}
The unique continuation property assumed in \cite{KL26}, although stronger in appearance than \eqref{eq:FH}, is in fact equivalent to it, so that the conclusion of \cite{KL26} is the last statement of Theorem~\ref{thm:main}.
But Theorem~\ref{thm:main} says \emph{more}: without any unique continuation assumption, it describes exactly the null controllable states in every time $T>T_0$.
\item
\emph{Setting.}
Unlike in \cite{KL26}, the spaces $Y$, $U$ and $U_T$ are Banach spaces, and no reflexivity is assumed.
This covers the state space $L^1$, the natural setting for transport equations, used in Section~\ref{sec:nonloc t} below.
\item
\emph{Dual formulation.}
Assumption \eqref{eq:primalH} is an existence statement about controls, which is not always easy to establish directly.
Beyond Theorem~\ref{thm:main}, we provide in Theorem~\ref{thm:obs} below a condition equivalent to \eqref{eq:primalH} in terms of observability estimates.
Together, Theorems~\ref{thm:main} and~\ref{thm:obs} cover a wide range of applications, as displayed in Section~\ref{sec:comparison} below.
\end{itemize}
\end{remark}

\begin{remark}[a necessary condition]
Item~\ref{main:i2} gives a necessary condition for \eqref{eq:primalH} which is often cheap to check: if $\sigma_F$ is infinite, or if a single $E_\lambda$ is infinite dimensional, then \eqref{eq:primalH} fails at every time.
More quantitatively, the proof of item~\ref{main:i2} gives $\sum_{\lambda\in\sigma_F}\dim E_\lambda\le\codim G+\dim F$ (see \eqref{eq:dimE} below).
\end{remark}

\subsection{An equivalent observability condition}\label{sec:nullobsfd}

To apply Theorem~\ref{thm:main}, one has to produce, for every initial state in $G$, a control steering it into $F$.
It is often more convenient to work on the adjoint system, where controllability turns into observability.
Our second main result gives a condition equivalent to \eqref{eq:primalH} in this form.

\begin{theorem}[equivalent observability condition]\label{thm:obs}
Let $(U_T)_{T>0}$ be a control class such that the regularity-gain condition \eqref{eq:RGC} holds and let $T_0>0$.
Consider the following two properties:
\begin{enumerate}
\item\label{obs:primal}
assumption \eqref{eq:primalH} of Theorem~\ref{thm:main} holds for some subspace $G\subset Y$ of finite codimension and some finite dimensional subspace $F\subset Y$;

\item\label{obs:CU}
there exist a compact operator $K:V_1 \to V_2$, $V_1$ a Banach space, $V_2$ a normed space, a linear map $L:D\to V_1$, $D$ a dense subspace of $Y'$, and $c>0$ such that, with $R:=KL$,
\begin{equation}\label{eq:single estim}
\norm*{S(T_0)^*z}_{Y'}+\norm*{Lz}_{V_1}\le c\left(\norm*{\Phi_{T_0}^*z}_{U_{T_0}'}+\norm*{Rz}_{V_2}\right),\qquad\forall z\in D.
\end{equation}
\end{enumerate}
Then \ref{obs:primal} implies \ref{obs:CU}.
If moreover $U_{T_0}$ is reflexive, then \ref{obs:CU} implies \ref{obs:primal}, so that the two properties are equivalent.
\end{theorem}

The two properties express the same fact from two points of view: controllability up to finite dimensional spaces for property~\ref{obs:primal}, and observability up to a compact remainder for property~\ref{obs:CU}.
The latter is the compactness part of the \emph{compactness--uniqueness method}, whose uniqueness part, the Fattorini--Hautus test \eqref{eq:FH}, enters through Theorem~\ref{thm:main}.

Property~\ref{obs:CU} is often the form used in practice: one usually first estimates $S(T_0)^*z$ by $\Phi_{T_0}^*z$ up to a remainder term $Rz$, then finds a splitting $R=KL$ with $K$ compact, and finally estimates $Lz$ by the observation and the remainder itself.
This is the case in many works of the literature, for instance \cite{CN21,HO22-JDE,FCLZ16,Tuc26} (see also Section~\ref{sec:comparison}), so that several results obtained separately, in different settings and by different arguments, all follow from Theorems~\ref{thm:obs} and~\ref{thm:main}.

Contrary to our first result, Theorem~\ref{thm:obs} is a purely functional analytic result: its proof, given in Section~\ref{sec:proofobs}, uses no property of control systems, neither the semigroup law, nor the composition property, nor the generator.

\begin{remark}[comparison with \cite{CN21,HO22-JDE}]\label{rem:HO22}
In \cite{CN21}, a compactness--uniqueness argument is used to prove the null controllability of one-dimensional first-order linear hyperbolic systems with a triangular internal coupling.
This argument was put in an abstract Hilbert framework, which fits ours, in \cite[Theorem 7.7]{HO22-JDE}, for a different purpose \cite[Theorem 7.1]{HO22-JDE}, with the Fattorini--Hautus test in place of the triangularity.
The other hypotheses are three estimates \cite[(52), (53) and (54)]{HO22-JDE}.
In our notation, with $Y$ and $U$ Hilbert spaces identified with their duals and $U_T=L^2(0,T;U)$, the first two read: there exist a compact operator $K\in\LL(V_1,V_2)$ ($V_1$ a Banach space, $V_2$ a normed space), a linear map $L:D(A^*)\to V_1$ and $c>0$ such that, with $R:=KL$,
\begin{align}
\label{est1 HO22}
\norm*{S(T_0)^*z}_{Y}&\leq c\left(\norm*{\Phi_{T_0}^*z}_{L^2(0,T_0;U)}+\norm*{Rz}_{V_2}\right),
\\
\label{est2 HO22}
\norm*{Lz}_{V_1}&\leq c\left(\norm*{\Phi_{T_0}^*z}_{L^2(0,T_0;U)}+\norm*{S(T_0)^*z}_{Y}\right),
\end{align}
for every $z\in D(A^*)$.
The third one requires, in addition, that
\begin{equation}\label{eq:h2}
\forall t_2>t_1>T_0,\ \exists c'>0,\quad
\norm*{S(t_1)^*z}_{Y}^2\le c'\left(\norm*{S(t_2)^*z}_{Y}^2+\int_{t_1}^{t_2}\norm*{B^*S(t)^*z}_{U}^2\,dt\right),
\end{equation}
for every $z\in D(A^*)$.
Inserting \eqref{est1 HO22} into \eqref{est2 HO22}, and then adding \eqref{est1 HO22}, gives \eqref{eq:single estim} with $D=D(A^*)$.
Therefore, Theorems~\ref{thm:obs} and~\ref{thm:main} draw the same conclusion \emph{without} \eqref{eq:h2}, and this can be essential in applications.
For instance, for the transport equation with nonlocal terms of Section~\ref{sec:nonloc t}, \eqref{eq:h2} fails for every $T_0<1$ (see Remark~\ref{rem:nbh2} below), so that \cite[Theorem 7.7]{HO22-JDE} cannot reach the \emph{minimal control time} of this equation, since this time is smaller than $1$.
\end{remark}


\section{Proof of Theorem~\ref{thm:main}}\label{sec:proofmain}

This section proves our first main result.

\subsection{Finite codimension and closedness}\label{sec:proof1}

We start with the simple observation that, if assumption \eqref{eq:primalH} holds at the time $T_0$, then it holds at every time $T\ge T_0$, with the same $G$ and with $F$ replaced by $S(T-T_0)F$.
Indeed, let $T>T_0$, and let $y^0\in G$ and $u\in U_{T_0}$ be as in \eqref{eq:primalH}.
Then $u\conca{T_0}0\in U_T$ by \eqref{eq:ALCSconcat}, and the composition property \eqref{eq:composition} gives
\[
S(T)y^0+\Phi_T\left(u\conca{T_0}0\right)=S(T-T_0)\left(S(T_0)y^0+\Phi_{T_0}u\right)\in S(T-T_0)F .
\]
In particular, item~\ref{main:i1} needs to be proved only for $T=T_0$.
To this end, we introduce the subspace
\begin{equation}\label{eq:MT0}
M_{T_0}:=S(T_0)^{-1}\,\Ima\Phi_{T_0}=\ens{y^0\in Y\st \exists u\in U_{T_0},\quad S(T_0)y^0+\Phi_{T_0}u=0}.
\end{equation}

\subsubsection{Finite codimension}

The following result reformulates assumption \eqref{eq:primalH} in terms of the set $M_{T_0}$.

\begin{proposition}\label{prop:reduction}
Assumption \eqref{eq:primalH} holds for some subspace $G\subset Y$ of finite codimension and some finite dimensional subspace $F\subset Y$ if, and only if, $\codim M_{T_0}<+\infty$.
In that case,
\begin{equation}\label{eq:codimbound}
\codim M_{T_0}\le\codim G+\dim F ,
\end{equation}
and \eqref{eq:primalH} holds with $G=M_{T_0}$ and $F=\ens*{0}$.
\end{proposition}

\begin{proof}
Assume \eqref{eq:primalH}.
It reads $S(T_0)G\subset\Ima\Phi_{T_0}+F$, that is, $G$ is contained in the set $G_F:=S(T_0)^{-1}(\Ima\Phi_{T_0}+F)$ of states that can be steered into $F$ in time $T_0$.
Since $M_{T_0}\subset G_F\subset Y$, the map $y+M_{T_0}\mapsto y+G_F$, from $Y/M_{T_0}$ to $Y/G_F$, is well defined, linear and surjective, with kernel $G_F/M_{T_0}$.
Thus, $M_{T_0}$ has finite codimension as soon as $Y/G_F$ and $G_F/M_{T_0}$ are finite dimensional, and then
\begin{equation}\label{eq:codimsum}
\codim M_{T_0}=\dim Y/G_F+\dim G_F/M_{T_0} .
\end{equation}
For the first space, the map $y+G\mapsto y+G_F$, from $Y/G$ to $Y/G_F$, is well defined, since $G\subset G_F$, and surjective, so that $\dim Y/G_F\le\codim G$.
For the second space, consider the linear maps
\[
G_F/M_{T_0}\to\left(\Ima\Phi_{T_0}+F\right)/\Ima\Phi_{T_0},
\qquad
F\to\left(\Ima\Phi_{T_0}+F\right)/\Ima\Phi_{T_0},
\]
defined by $y+M_{T_0}\mapsto S(T_0)y+\Ima\Phi_{T_0}$ and $f\mapsto f+\Ima\Phi_{T_0}$ respectively.
The first one is well defined and injective, since $M_{T_0}=S(T_0)^{-1}\,\Ima\Phi_{T_0}$, and the second one is surjective, so that $\dim G_F/M_{T_0}\le\dim F$.
Collecting these bounds in \eqref{eq:codimsum} gives \eqref{eq:codimbound}.

Conversely, assume $\codim M_{T_0}<+\infty$.
For $y^0\in M_{T_0}$, there is $u\in U_{T_0}$ such that $S(T_0)y^0+\Phi_{T_0}u=0$.
Hence \eqref{eq:primalH} holds with $G=M_{T_0}$, which has finite codimension, and $F=\ens*{0}$.
\end{proof}

\subsubsection{Closedness}

The finite codimension of $M_{T_0}$ is thus equivalent to assumption \eqref{eq:primalH}.
The closedness relies on the following lemma.

\begin{lemma}\label{lem:finitecodim}
Let $X_1,X_2$ be Banach spaces and let $\Lambda\in\LL(X_1,X_2)$.
If $\Ima\Lambda$ has finite codimension in $X_2$, then it is closed.
\end{lemma}

\begin{proof}
Set $M:=\Ima\Lambda$ and let $M_0\subset X_2$ be a finite dimensional subspace such that $X_2=M\oplus M_0$ algebraically.
The map
\[
\widehat\Lambda:\left(X_1/\ker \Lambda\right)\times M_0\to X_2,\qquad \left(x+\ker \Lambda,f\right)\mapsto \Lambda x+f,
\]
is well defined, linear, continuous and bijective, between Banach spaces, $\ker \Lambda$ being closed and $M_0$ finite dimensional.
By the bounded inverse theorem, it is an isomorphism, so that $M=\widehat\Lambda((X_1/\ker \Lambda)\times\ens*{0})$ is closed.
\end{proof}

With this lemma, we can prove that $M_{T_0}$ is closed in $Y$ and thus finish the proof of item~\ref{main:i1} of Theorem~\ref{thm:main}.

\begin{proposition}\label{prop:closed}
Assume \eqref{eq:primalH}.
Then $M_{T_0}$ is closed in $Y$.
\end{proposition}

\begin{proof}
Consider on $X:=\ens*{(y^0,u)\in Y\times U_{T_0}\st S(T_0)y^0+\Phi_{T_0}u=0}$ the projection $P:X\to Y$ onto the first component, defined by $P(y^0,u):=y^0$.
The space $X$ is closed in $Y\times U_{T_0}$, and thus a Banach space.
By definition \eqref{eq:MT0}, we have $\Ima P=M_{T_0}$.
By Proposition~\ref{prop:reduction}, $M_{T_0}$ has finite codimension.
Lemma~\ref{lem:finitecodim} then shows that it is closed.
\end{proof}

\subsection{Identification of the obstruction}\label{sec:proof2}

This subsection proves items~\ref{main:i2} and~\ref{main:i3} of Theorem~\ref{thm:main}.
We recall that $E:=\bigoplus_{\lambda\in\sigma_F}E_\lambda\subset Y'$ is the space spanned by the unobservable generalized eigenvectors of the pair $(A^*,B^*)$, with the notation of Definition~\ref{def:unobs}.

\subsubsection{Finite dimensional invariant subspaces}\label{sec:knownE}

Some of the following facts are classical when the state space $Y$ is finite dimensional.
However, in our setting, the adjoint semigroup $(S(t)^*)_{t\ge0}$ is in general not strongly continuous on $Y'$, so that we give complete proofs.

\begin{proposition}\label{prop:knownE}
\begin{enumerate}
\item\label{knownE:inv}
A finite dimensional subspace $N\subset Y'$ satisfies $S(t)^*N\subset N$ for every $t\ge0$ if, and only if, $N\subset D(A^*)$ and $A^*N\subset N$.
In that case, denoting by $L$ the restriction of $A^*$ to $N$, we have
\begin{equation}\label{eq:expN}
S(t)^*f=e^{tL}f,\qquad\forall t\ge0,\ \forall f\in N,
\end{equation}
and in particular $S(t)^*N=N$ for every $t\ge0$.

\item\label{knownE:max}
For every $\tau>0$, the space $E$ is the union of all finite dimensional subspaces $N\subset(\Ima\Phi_\tau)^{\perp}$ such that $S(t)^*N\subset N$ for every $t\ge0$.
\end{enumerate}
\end{proposition}

\begin{proof}
\begin{enumerate}
\item
Assume first that $S(t)^*N\subset N$ for every $t\ge0$.
Then $((S(t)^*)_{|N})_{t\ge0}$ is a semigroup on $N$, and we check that it is uniformly continuous.
Let $f_1,\dots,f_n$ be a basis of $N$, and let $y_1,\dots,y_n\in Y$ with $\pair*{f_i}{y_j}_{Y',Y}=\delta_{ij}$ (see e.g.\ \cite[Section~2.4, Example~2]{Bre11}).
Then every $f\in N$ can be decomposed as $f=\sum_j\pair*{f}{y_j}_{Y',Y}f_j$.
Applying this to $f=S(t)^*z-z$ with $z\in N$, and using $\abs*{\pair*{S(t)^*z-z}{y_j}_{Y',Y}}=\abs*{\pair*{z}{S(t)y_j-y_j}_{Y',Y}}\le\norm*{z}_{Y'}\norm*{S(t)y_j-y_j}_{Y}$, we obtain
\[
\norm*{S(t)^*z-z}_{Y'}\le \norm*{z}_{Y'}\sum_{j=1}^n \norm*{S(t)y_j-y_j}_{Y} \norm*{f_j}_{Y'} \longrightarrow0\quad\text{as }t\to0^+ .
\]
Hence $((S(t)^*)_{|N})_{t\ge0}$ is a uniformly continuous semigroup on $N$, so that $(S(t)^*)_{|N}=e^{tL}$ for some linear operator $L:N\to N$ \cite[Theorem~I.3.7]{EN00}, and $(S(t)^*f-f)/t\to Lf$ in $Y'$ for every $f\in N$.
For $f\in N$ and $y\in D(A)$, this gives
\[
\pair{Lf}{y}_{Y',Y}=\lim_{t\to0^+}\pair{f}{\frac{S(t)y-y}{t}}_{Y',Y}=\pair{f}{Ay}_{Y',Y} ,
\]
so that $y\mapsto\pair*{f}{Ay}_{Y',Y}$ is continuous on $D(A)$ for the norm of $Y$, that is, $f\in D(A^*)$, with $A^*f=Lf\in N$.
Thus $N\subset D(A^*)$ and $A^*N\subset N$.

Assume now that $N\subset D(A^*)$ and $A^*N\subset N$, and let us prove \eqref{eq:expN}, which gives the converse implication.
Let $f\in N$, and put $\psi(t):=e^{tL}f$ for $t\ge0$.
The space $N$ being equipped with the norm of $Y'$, we have $\psi \in C^1([0,+\infty);N)$ with $\psi'=A^*\psi$ and $\psi(0)=f$.
Let $y\in D(A)$ and $t>0$.
The function $s\mapsto\pair*{\psi(s)}{S(t-s)y}_{Y',Y}$ is differentiable on $[0,t]$, with derivative $\pair*{A^*\psi(s)}{S(t-s)y}_{Y',Y}-\pair*{\psi(s)}{AS(t-s)y}_{Y',Y}=0$.
Hence it is constant, and $\pair*{\psi(t)}{y}_{Y',Y}=\pair*{f}{S(t)y}_{Y',Y}=\pair*{S(t)^*f}{y}_{Y',Y}$.
As $D(A)$ is dense in $Y$, this gives \eqref{eq:expN}.
Finally, $e^{tL}$ is invertible on $N$, so that $S(t)^*N=N$.

\item
Let $\tau>0$ and $f\in E$, and write $f=\sum_{\lambda\in\Lambda}f_\lambda$, with $\Lambda\subset\C$ finite and $f_\lambda\in E_\lambda$.
By definition of $E_\lambda$, there exists $m_\lambda\geq 1$ such that $(\lambda-A^*)^{m_\lambda}f_\lambda=0$.
Thus, $f$ belongs to the span $N$ of the vectors $(\lambda-A^*)^jf_\lambda$, $\lambda\in\Lambda$, $j<m_\lambda$.
Clearly, $N$ is finite dimensional, is contained in $D(A^*)$, is invariant under $A^*$, since $(\lambda-A^*)^{m_\lambda}f_\lambda=0$, and is contained in $\ker B^*$.
By item~\ref{knownE:inv}, $S(t)^*N=N\subset\ker B^*$ for every $t\ge0$, so that $N\subset(\Ima\Phi_\tau)^{\perp}$, from \eqref{eq:link}.

Conversely, let $N\subset(\Ima\Phi_\tau)^{\perp}$ be a finite dimensional subspace such that $S(t)^*N\subset N$ for every $t\ge0$.
Let $f\in N$, and let $\varphi(t):=S(\tau-t)^*f$.
Since $f\in(\Ima\Phi_\tau)^{\perp}$, \eqref{eq:link} gives
\begin{equation}\label{eq:uzeroobs}
\int_0^{\tau}\pair{B^*\varphi(t)}{u(t)}_{U',U}\,dt=\pair{f}{\Phi_\tau u}_{Y',Y}=0 ,\qquad\forall u\in U_{\tau} .
\end{equation}
Now we observe that, by item~\ref{knownE:inv}, $t\mapsto\varphi(t)=e^{(\tau-t)L}f$ is continuous from $[0,\tau]$ to $D(A^*)$, the graph norm and the norm of $Y'$ being equivalent on $N$.
Consequently, $B^*\varphi$ is continuous and the left-hand side of \eqref{eq:uzeroobs} is a bounded linear form of $u\in L^1(0,\tau;U)$, which vanishes on $L^1(0,\tau;U)$, since it vanishes on the dense subspace $U_\tau$.
Thus, $B^*\varphi=0$ and $B^*\varphi(\tau)=B^*f=0$.
Hence, $N\subset\ker B^*$.
Then, $N$ is the direct sum of the generalized eigenspaces of $L$ (recall that we work in complex spaces).
If $f$ lies in the one associated with $\lambda\in\sigma(L)$, then $(\lambda-A^*)^mf=0$ for some $m\in\N$, and $(\lambda-A^*)^jf\in N\subset\ker B^*$ for every $j\in\N$, that is, $f\in E_\lambda$.
Hence $N\subset E$.
\end{enumerate}
\end{proof}

\paragraph{Unobservable invariant subspaces are obstructions.}
Let $T>0$.
Every finite dimensional subspace $N\subset Y'$ such that $S(t)^*N\subset N$ for every $t\ge0$ and $N\subset(\Ima\Phi_T)^{\perp}$ is an obstruction to null controllability in time $T$:
$$
N\subset\left(S(T)^{-1}\,\Ima\Phi_T\right)^{\perp}.
$$
Indeed, if $f\in N$ and $y^0\in S(T)^{-1}\,\Ima\Phi_T$, then $f=S(T)^*z$ for some $z\in N$, by item~\ref{knownE:inv} of Proposition~\ref{prop:knownE}, and $S(T)y^0=\Phi_Tu$ for some $u\in U_T$, so that $\pair*{f}{y^0}_{Y',Y}=\pair*{z}{\Phi_Tu}_{Y',Y}=0$.

Using now item~\ref{knownE:max} of Proposition~\ref{prop:knownE} with $\tau=T$, we get the inclusion $\subset$ in \eqref{eq:reach} for any $T>0$ and without any assumption:
\begin{equation}\label{eq:easyE}
E \subset \left(S(T)^{-1}\,\Ima\Phi_T\right)^{\perp},\qquad\forall T>0.
\end{equation}

\subsubsection{Finite dimension of $E$}\label{sec:Efinite}

\begin{proof}[Proof of Theorem~\ref{thm:main}, item~\ref{main:i2}]
Assume \eqref{eq:primalH}.
By \eqref{eq:easyE}, we have $E\subset M_{T_0}^{\perp}$.
By Proposition~\ref{prop:reduction}, $M_{T_0}$ has finite codimension, so that $\dim M_{T_0}^{\perp}\le\codim M_{T_0}$, by \eqref{eq:dimMperpcodim}.
Thus, $E$ is finite dimensional.
Then, every subspace $E_\lambda$ is finite dimensional and, since they are nonzero and in direct sum for $\lambda\in\sigma_F$, the set $\sigma_F$ is finite.
\end{proof}

From the proof above and \eqref{eq:codimbound} we have the announced bound
\begin{equation}\label{eq:dimE}
\dim E\le\codim G+\dim F .
\end{equation}

\subsubsection{Decoupling the free evolution time from the control time}\label{sec:strategy}

Let $T>T_0$.
Our aim is to describe the set $S(T)^{-1}\,\Ima\Phi_T$, or equivalently its annihilator, since it is closed by item~\ref{main:i1}.
We already know that this annihilator contains $E$, by \eqref{eq:easyE}, and it remains to show the converse inclusion.
To this end, we use the characterization in item~\ref{knownE:max} of Proposition~\ref{prop:knownE}.
We already know that the annihilator is finite dimensional, by item~\ref{main:i1} and \eqref{eq:dimMperpcodim}.
It remains to show two properties: stability under every $S(t)^*$ and inclusion in some $(\Ima\Phi_\tau)^{\perp}$.
Since $S(T)^{-1}\,\Ima\Phi_T$ is closed, this amounts, by \eqref{eq:dualincl}, to showing that $S(T)^{-1}\,\Ima\Phi_T$ is stable under every $S(t)$ and contains some $\Ima\Phi_\tau$.
For \emph{exact} controllability, the same two properties are needed for $\Ima\Phi_T$, and they follow from a plateau of the annihilators, since we always have $S(t)\,\Ima\Phi_T\subset\Ima\Phi_{T+t}$ and $\Ima\Phi_\tau\subset\Ima\Phi_T$ for $\tau\le T$, by Remark~\ref{rem:concat for NC}.
For the sets $S(T)^{-1}\,\Ima\Phi_T$, the analogous inclusions cannot be read off directly.
Indeed, for $t,\tau>0$, Remark~\ref{rem:concat for NC} only gives
\begin{equation}\label{eq:strategyincl}
S(t)\left(S(T)^{-1}\,\Ima\Phi_T\right) \subset S(T)^{-1}\,\Ima\Phi_{T+t},
\qquad
\Ima \Phi_\tau \subset S(T)^{-1}\,\Ima\Phi_{T+\tau}.
\end{equation}
We see that the free evolution time is unchanged, whereas the control window has become longer.
Thus, starting from $S(T)^{-1}\,\Ima\Phi_T$, we end up in sets of a different kind, in which the two durations no longer coincide.
This motivates the \emph{key idea} of the proof: we decouple the duration of the free evolution and the length of the control window, and we let only the latter grow.
We thus fix the former at the time $T_0$ and consider
\begin{equation}\label{eq:PT}
M_T:=S(T_0)^{-1}\,\Ima\Phi_T=\ens{y^0\in Y\st \exists u\in U_T,\quad S(T_0)y^0+\Phi_Tu=0}, \qquad T \geq T_0,
\end{equation}
which agrees with \eqref{eq:MT0} for $T=T_0$.
We will see that nothing is lost by doing so.

\paragraph{Basic properties of the sets $M_T$.}
With the free evolution time fixed at $T_0$, the inclusions \eqref{eq:strategyincl} now stay within the family $(M_T)_{T\ge T_0}$: for $T\ge T_0$ and $t,\tau>0$,
\begin{equation}\label{eq:PTprops}
S(t)\,M_{T}\subset M_{T+t},\qquad
\Ima\Phi_\tau\subset M_{T_0+\tau} .
\end{equation}
Indeed, let $y^0\in M_{T}$, that is, $S(T_0)y^0\in\Ima\Phi_T$.
Then $S(T_0)S(t)y^0=S(t)S(T_0)y^0$ belongs to $S(t)\,\Ima\Phi_T\subset\Ima\Phi_{T+t}$, by \eqref{eq:monoP1}.
This property also gives $S(T_0)\,\Ima\Phi_\tau\subset\Ima\Phi_{T_0+\tau}$, that is, the second inclusion.
The sets $M_T$ are also nondecreasing by \eqref{eq:monoP2}:
\begin{equation}\label{eq:MTnested}
M_{T_1}\subset M_{T_2},\qquad T_1<T_2.
\end{equation}
Finally, under assumption \eqref{eq:primalH}, the sets $M_T$ are closed subspaces of finite codimension of $Y$, like $M_{T_0}$.
Indeed, $M_T$ contains $M_{T_0}$, by \eqref{eq:MTnested}, so that its codimension is at most that of $M_{T_0}$, and it can then be deduced that $M_T$ is closed as in the proof of Proposition~\ref{prop:closed}.

The obstruction will be described by the annihilators
$$
N_T:=M_T^{\perp}\subset Y',\qquad T\ge T_0 .
$$
They inherit the properties dual to those of the sets $M_T$.
Under assumption \eqref{eq:primalH}, they are finite dimensional, by \eqref{eq:dimMperpcodim}, and, by \eqref{eq:dualincl}, \eqref{eq:PTprops} and \eqref{eq:MTnested} give
\begin{equation}\label{eq:NTprops}
S(t)^*\,N_{T+t}\subset N_T,\qquad
N_{T_0+\tau}\subset\left(\Ima\Phi_\tau\right)^{\perp},\qquad
N_{T_2}\subset N_{T_1},
\end{equation}
for every $T\ge T_0$, $t,\tau>0$ and $T_1\le T_2$.

\subsubsection{Identification of the annihilators}\label{sec:identif}

\begin{proposition}\label{prop:identif}
We have $N_T \supset E$ for every $T\ge T_0$.
If \eqref{eq:primalH} holds, then
$$
N_{T}=E ,\qquad\forall T>T_0 .
$$
\end{proposition}

\begin{proof}
\emph{The inclusion $\supset$.}
For any $T\ge T_0$, the argument giving \eqref{eq:easyE}, with $S(T_0)$ in place of $S(T)$, gives $E\subset(S(T_0)^{-1}\,\Ima\Phi_T)^{\perp}=M_T^{\perp}=N_T$.

\emph{The inclusion $\subset$.}
Assume now that \eqref{eq:primalH} holds, and fix $T>T_0$.
\begin{enumerate}
\item \emph{A plateau of the annihilators.}
By nesting (the third property in \eqref{eq:NTprops}) and the finite dimensionality of $N_{T_0}$, the function $s\mapsto\dim N_{s}$ is nonincreasing on $[T_0,+\infty)$, with values in $\N$.
It thus takes finitely many values on $[T_0,T]$, and it is constant on some open interval $I\subset(T_0,T)$.
By nesting, $N_{s}=:N$ does not depend on $s\in I$, and $N_{T}\subset N$.

\item \emph{Invariance under the semigroup.}
Let $s\in I$ and $\varepsilon>0$ be such that $[s,s+\varepsilon]\subset I$.
For every $t\in[0,\varepsilon]$, the first property in \eqref{eq:NTprops} gives $S(t)^*N=S(t)^*N_{s+t}\subset N_{s}=N$.
By the semigroup property, $S(t)^*N\subset N$ for every $t\ge0$.

\item \emph{Orthogonality to the reachable states.}
Let $s\in I$ and set $\tau:=s-T_0>0$.
As $N=N_{s}=N_{T_0+\tau}$, the second property in \eqref{eq:NTprops} gives $N\subset(\Ima\Phi_\tau)^{\perp}$.

\item \emph{Conclusion.}
By item~\ref{knownE:max} of Proposition~\ref{prop:knownE}, with the $\tau$ of the previous step, $N$ is contained in $E$, and so is $N_T$, since $N_T\subset N$.
\end{enumerate}
\end{proof}

\subsubsection{End of the proof}\label{sec:proofmainend}

\begin{proof}[Proof of Theorem~\ref{thm:main}, item~\ref{main:i3}]
Let $T>T_0$.
We have
\[
S(T)^{-1}\,\Ima\Phi_T
=S(T-T_0)^{-1}\,M_T
=S(T-T_0)^{-1}\,E_{\perp}
=\left(S(T-T_0)^*E\right)_{\perp}
=E_{\perp} .
\]
The first equality follows from the definition \eqref{eq:PT} of $M_T$ and $S(T)=S(T_0)S(T-T_0)$.
The second one holds because $M_T=(M_T^{\perp})_{\perp}=(N_T)_{\perp}$, by \eqref{prop perp} and the closedness of $M_T$, and $N_T=E$, by Proposition~\ref{prop:identif}.
The third one holds for any subset of $Y'$: for $y\in Y$ and $t\ge0$, we have $S(t)y\in E_{\perp}$ if, and only if, $\pair*{S(t)^*f}{y}_{Y',Y}=\pair*{f}{S(t)y}_{Y',Y}=0$ for every $f\in E$.
For the last equality, the space $E$ is finite dimensional (Section~\ref{sec:Efinite}) and invariant under $A^*$, so that $S(t)^*E=E$ for every $t\ge0$ by item~\ref{knownE:inv} of Proposition~\ref{prop:knownE}.
\end{proof}


\section{Proof of Theorem~\ref{thm:obs}}\label{sec:proofobs}

This section proves our second main result.
As already mentioned, the proof only uses that $S(T_0)$ and $\Phi_{T_0}$ are bounded operators with the same target space $Y$, and we abbreviate $S:=S(T_0)$, $\Phi:=\Phi_{T_0}$ and $W:=U_{T_0}$, so that
$$S\in\LL(Y), \qquad \Phi\in\LL(W,Y).$$

\subsection{Reformulation by duality}\label{sec:hubdual}

We first translate property~\ref{obs:primal}, an existence statement, into an estimate on the adjoint operators, by duality.

\begin{proposition}\label{prop:primaltoaug}
If property~\ref{obs:primal} of Theorem~\ref{thm:obs} holds, then there exist a finite dimensional subspace $G_0\subset Y$ and $c>0$ such that
\begin{equation}\label{eq:hub}
\norm*{S^*z}_{Y'}\le c\left(\norm*{\Phi^*z}_{W'}+\norm*{\iota^* S^*z}_{G_0'}\right), \qquad \forall z\in Y',
\end{equation}
where $\iota\in\LL(G_0,Y)$ is the inclusion map.
Conversely, if $W$ is reflexive and \eqref{eq:hub} holds for some finite dimensional subspace $G_0\subset Y$ and $c>0$, then property~\ref{obs:primal} holds.
\end{proposition}

\begin{proof}
\begin{enumerate}
\item
If property \ref{obs:primal} holds, then it holds with $F=\ens*{0}$, up to changing $G$ into another subspace of finite codimension, by Proposition~\ref{prop:reduction}.
Writing then $Y=G\oplus G_0$, with $\dim G_0<+\infty$, we obtain
\begin{equation}\label{eq:augrange}
\Ima S\subset\Ima\Phi+SG_0.
\end{equation}
Conversely, if \eqref{eq:augrange} holds for some finite dimensional subspace $G_0\subset Y$, then this gives \eqref{eq:primalH} with $G=Y$ and $F=SG_0$.
Property \ref{obs:primal} is thus equivalent to the existence of a finite dimensional subspace $G_0\subset Y$ such that \eqref{eq:augrange} holds.

\item
Let $G_0\subset Y$ be a finite dimensional subspace, equipped with the norm of $Y$.
The right-hand side of \eqref{eq:augrange} is the range of the augmented map
$$
\Phihat: G_0\times W\to Y,
\qquad
\Phihat\left(g,u\right):=Sg+\Phi u,
$$
so that \eqref{eq:augrange} is the range inclusion $\Ima S\subset\Ima\Phihat$.
We equip $G_0\times W$ with the max norm $\norm*{(g,u)}:=\max\ens*{\norm*{g}_{Y}, \, \norm*{u}_{W}}$, so that its dual is $G_0'\times W'$, with the sum norm.
A computation shows that $\Phihat^*z=(\iota^*S^*z, \, \Phi^*z)$.
The implications then follow from duality (Theorem~\ref{thm:duality}).
\end{enumerate}
\end{proof}

One implication of Theorem~\ref{thm:obs} is then immediate.

\begin{proof}[Proof of Theorem~\ref{thm:obs}, first implication]
Assume property~\ref{obs:primal}.
By Proposition~\ref{prop:primaltoaug}, we have the estimate \eqref{eq:hub} for some finite dimensional subspace $G_0\subset Y$ and $c>0$.
Then, property~\ref{obs:CU} holds, up to the value of $c$, with $L:=S^*:Y' \to Y'$ and $K:=\iota^*:Y' \to G_0'$, which is compact, as a finite rank operator.
\end{proof}

\subsection{Compactness and the Peetre lemma}\label{sec:obstohub}

It remains to establish the converse implication, and we do so by proving the estimate \eqref{eq:hub}.
The next step is the \emph{core of the proof}, where the compactness of $K$ and the Peetre lemma are used.

\begin{proposition}\label{prop:obstopeetre}
If property~\ref{obs:CU} of Theorem~\ref{thm:obs} holds, then there exist a finite dimensional subspace $N\subset Y'$ and $c>0$ such that
\begin{equation}\label{eq:distT0}
\dist_{Y'}\left(S^*z,N\right)\le c\,\norm*{\Phi^*z}_{W'},\qquad\forall z\in Y'.
\end{equation}
\end{proposition}

\begin{proof}
We first decouple \eqref{eq:single estim} into the two estimates
\begin{align}
\label{eq:decL}
\norm*{Lz}_{V_1}&\les\norm*{\Phi^*z}_{W'}+\norm*{KLz}_{V_2},
\\
\label{eq:decS}
\norm*{S^*z}_{Y'}&\les\norm*{\Phi^*z}_{W'}+\norm*{Lz}_{V_1},
\end{align}
for every $z\in D$.
These two estimates have different purposes.
Estimate \eqref{eq:decL} is a compactness--uniqueness estimate, to which the Peetre lemma applies (steps 1 and 2).
Estimate \eqref{eq:decS} involves no compactness, and it transfers the conclusion from $L$ to $S^*$ (step 3).
\begin{enumerate}
\item
\emph{Augmentation of $L$.}
We would like to see \eqref{eq:decL} as an estimate on the range of $L$.
This is possible if the map $Lz\mapsto\Phi^*z$ is well defined, but this is not the case in general.
Fortunately, we can augment $L$ by $\Phi^*$, which costs nothing in \eqref{eq:decL}.
We equip $V_1\times W'$ with the sum norm, which makes it a Banach space, $V_1$ and the dual $W'$ being complete.
We then define $\widehat L:D\to V_1\times W'$ and $\widehat K:V_1\times W'\to V_2$ by
\[
\widehat Lz:=\left(Lz,\Phi^*z\right),
\qquad
\widehat K(v,w):=Kv.
\]
The operator $\widehat K$ is compact, and $\widehat K\widehat Lz=KLz$ for every $z\in D$.
Adding $\norm*{\Phi^*z}_{W'}$ to both sides of \eqref{eq:decL}, we obtain
\begin{equation}\label{eq:HobsK}
\norm*{\widehat Lz}_{V_1\times W'}\les\norm*{\Phi^*z}_{W'}+\norm*{\widehat K\widehat Lz}_{V_2},\qquad\forall z\in D.
\end{equation}
Let $X$ be the closure of $\Ima\widehat L$ in $V_1\times W'$, which is a Banach space.
The map $\Lambda:\widehat Lz\mapsto\Phi^*z$ is well defined and it extends to $X$, still denoted by $\Lambda\in\LL(X,W')$.

\item
\emph{Application of the Peetre lemma.}
By density, estimate \eqref{eq:HobsK} extends to
\[
\norm*{\xi}_{V_1\times W'}\les\norm*{\Lambda\xi}_{W'}+\norm*{\widehat K\xi}_{V_2},\qquad\forall\xi\in X.
\]
The restriction of $\widehat K$ to $X$ is compact and $X$ is a Banach space.
By the Peetre lemma (Lemma~\ref{lem:peetre}), the kernel of $\Lambda$ is finite dimensional and its range is closed in $W'$.
Clearly, the kernel is of the form $\ker\Lambda=N_1\times\ens*{0}$, for some subspace $N_1\subset V_1$, which is thus finite dimensional.
On the other hand, it is well known that the range of $\Lambda$ is closed if, and only if, there exists $c'>0$ such that
\begin{equation}\label{eq:distker}
\dist_{V_1\times W'}(\xi,\ker\Lambda)\le c'\norm*{\Lambda\xi}_{W'},\qquad\forall\xi\in X
\end{equation}
(see e.g.\ \cite[Exercise~2.14 and its solution p.~382]{Bre11}).
For every $z\in D$ and $v\in N_1$, we have $\norm*{Lz-v}_{V_1}\le\norm*{\widehat Lz-(v,0)}_{V_1\times W'}$, so that, taking the infimum over $v\in N_1$, we deduce from \eqref{eq:distker} with $\xi=\widehat Lz$ that
\[
\dist_{V_1}(Lz,N_1)\le c'\norm*{\Phi^*z}_{W'},\qquad\forall z\in D.
\]

\item
\emph{Transfer to $S^*$.}
By \eqref{eq:decS}, we have $\norm*{S^*z}_{Y'}\les\norm*{\widehat Lz}_{V_1\times W'}$ for every $z\in D$.
Thus, the linear map $\Sigma:\widehat Lz\mapsto S^*z$ is well defined on $\Ima\widehat L$, and it extends to $X$, still denoted by $\Sigma\in\LL(X,Y')$.
Let $N:=\Sigma(\ker\Lambda)$, which is finite dimensional, and let $z\in D$.
For every $\xi=(v,0)\in\ker\Lambda$, we have
\[
\norm*{S^*z-\Sigma\xi}_{Y'}=\norm*{\Sigma(\widehat Lz-\xi)}_{Y'}\les\norm*{\widehat Lz-\xi}_{V_1\times W'}=\norm*{Lz-v}_{V_1}+\norm*{\Phi^*z}_{W'}.
\]
Taking the infimum over $v\in N_1$ gives $\dist_{Y'}(S^*z,N)\les\dist_{V_1}(Lz,N_1)+\norm*{\Phi^*z}_{W'}$, hence \eqref{eq:distT0}, by step~2 and density of $D$ in $Y'$.
\end{enumerate}
\end{proof}

The converse implication now follows easily.

\begin{proof}[Proof of Theorem~\ref{thm:obs}, converse implication]
Assume property~\ref{obs:CU}.
In view of Proposition~\ref{prop:primaltoaug}, it suffices to show the estimate \eqref{eq:hub} for some finite dimensional subspace $G_0\subset Y$.
Proposition~\ref{prop:obstopeetre} gives a finite dimensional subspace $N\subset Y'$ and $c>0$ such that \eqref{eq:distT0} holds.
Since $N$ is finite dimensional, there exists a finite dimensional subspace $G_0\subset Y$ such that $N\cap G_0^{\perp}=\ens*{0}$.
Indeed, if $f_1,\dots,f_n$ is a basis of $N$, then there exist $y_1,\dots,y_n\in Y$ such that $\pair*{f_i}{y_j}_{Y',Y}=\delta_{ij}$ (see e.g.\ \cite[Section~2.4, Example~2]{Bre11}), and we can take $G_0=\Span\ens*{y_1,\ldots,y_n}$.
Let $\iota\in\LL(G_0,Y)$ be the inclusion map.
Then, the map $f\mapsto\norm*{\iota^*f}_{G_0'}$ is a norm on $N$, which is equivalent to $\norm*{\cdot}_{Y'}$, $N$ being finite dimensional.
As a result, for every $z\in Y'$ and $f\in N$, using twice the decomposition $S^*z=(S^*z-f)+f$, we obtain
\[
\norm*{S^*z}_{Y'}\les\norm*{S^*z-f}_{Y'}+\norm*{\iota^*f}_{G_0'}\les\norm*{S^*z-f}_{Y'}+\norm*{\iota^*S^*z}_{G_0'}.
\]
Taking the infimum over $f\in N$ and using \eqref{eq:distT0}, we obtain \eqref{eq:hub}.
\end{proof}


\section{Applications}\label{sec:comparison}

This section gives applications of Theorems~\ref{thm:main} and~\ref{thm:obs}.
The first two revisit results of the literature.
For the heat equation with a nonlocal term of \cite{FCLZ16}, the analyticity assumption is replaced by the Fattorini--Hautus test, and the null controllable states are described when the test fails.
For the simultaneous controllability problem of \cite{Tuc26}, the simultaneous approximate controllability is replaced by the test, without analyticity.
The last one, in Section~\ref{sec:nonloc t}, is completely new: for a transport equation with nonlocal terms in the domain and at the boundary, controlled from an interior band, we obtain the null controllability in the minimal time, in every space $L^r$ with $r<\infty$, including $L^1$, together with the description of the null controllable states, a problem which is out of reach of the methods available so far.
These examples are not meant as a catalogue, but to show that our results cover a \emph{wide range} of partial differential equations, hyperbolic, parabolic, or coupling both types, and that they provide a \emph{unified point of view} and often \emph{simpler proofs}.

The partial differential equations below have real-valued coefficients, and we apply our results, as usual, to their complexifications, which have the same controllability properties by invariance under conjugation.


\subsection{Heat equation with a nonlocal term}\label{sec:compFCLZ}

We begin with parabolic equations with a nonlocal term in space.
Their controllability has attracted some attention \cite{FCLZ16,LZ18,MT18,ELT22}, the main difficulty being that the nonlocal term, although compact, cannot be absorbed in the Carleman estimates used for local terms (see \cite[Remark 4]{FCLZ16}).
Here, we show that the null controllability result of \cite{FCLZ16} falls naturally within our main results, and we describe in addition the null controllable states when the Fattorini--Hautus test fails, as can happen even for simple kernels.

\subsubsection{Setting}

Following \cite{FCLZ16}, we consider the heat equation with a nonlocal term
\begin{equation}\label{eq:FCLZsyst}
\left\{
\begin{aligned}
&y_t-\Delta y+\int_\Omega k(x,\xi)\,y(t,\xi)\,d\xi=u(t,x)\indic_\omega(x)&&\text{in }(0,T)\times\Omega,\\
&y=0&&\text{on }(0,T)\times\partial\Omega,\\
&y(0,\cdot)=y^0&&\text{in }\Omega.
\end{aligned}
\right.
\end{equation}
Here $\Omega\subset\R^N$ is a bounded connected open set with boundary of class $C^2$ and $\omega\subset\Omega$ is a nonempty open subset.
The kernel $k$ belongs to $L^2(\Omega\times\Omega)$ and is real valued.
The state space and the control space are $Y=U:=L^2(\Omega)$, and the control class is $U_T:=L^2(0,T;U)=L^2((0,T)\times\Omega)$.
All these Hilbert spaces are identified with their duals.
The system \eqref{eq:FCLZsyst} reads $y'=Ay+Bu$, where
\[
(Ay)(x):=\Delta y(x)-\int_\Omega k(x,\xi)\,y(\xi)\,d\xi,\qquad (Bu)(x):=\indic_\omega(x)\,u(x),
\]
with $D(A):=H^2(\Omega)\cap H^1_0(\Omega)$.
It is well-known that the operator $A$ generates a $C_0$-semigroup on $Y$ and, since $B\in\LL(U,Y)$, the system \eqref{eq:FCLZsyst} fits into the framework of Section~\ref{sec:setting} (Remark~\ref{rem:bounded B} and Example~\ref{ex:UT}).

A computation shows that $D(A^*)=D(A)$ and
\begin{equation}\label{eq:def AB nonloca heat}
(A^*z)(x)=\Delta z(x)-\int_\Omega k(\xi,x)\,z(\xi)\,d\xi,\qquad B^*=B.
\end{equation}

For $T>0$ and $z\in L^2(\Omega)$, we denote by $\varphi$ the solution of the backward adjoint system
\[
\left\{
\begin{aligned}
&-\varphi_t(t,x)-\Delta\varphi(t,x)+\int_\Omega k(\xi,x)\,\varphi(t,\xi)\,d\xi=0&&\text{in }(0,T)\times\Omega,\\
&\varphi=0&&\text{on }(0,T)\times\partial\Omega,\\
&\varphi(T,\cdot)=z&&\text{in }\Omega.
\end{aligned}
\right.
\]
We denote by $p$ the solution of the same system with $k=0$.
We have
\[
\norm*{S(T)^*z}_{L^2(\Omega)}=\norm*{\varphi(0,\cdot)}_{L^2(\Omega)},
\qquad
\norm*{\Phi_T^*z}_{U_T}^2=\int_0^T\norm*{\indic_\omega\varphi(t,\cdot)}_{L^2(\Omega)}^2\,dt.
\]

Let $(\phi_j)_{j\ge1}$ be an orthonormal basis of $L^2(\Omega)$ made of eigenfunctions of the Dirichlet Laplacian, with eigenvalues $0<\lambda_1\le\lambda_2\le\cdots$.
We denote by $\pair*{f}{g}_{L^2(\Omega)}:=\int_\Omega f\,\conj{g}\,dx$ the inner product of $L^2(\Omega)$.
We equip $H^1_0(\Omega)$ with the norm $\norm{\nabla\cdot}_{L^2(\Omega)}$, so that $\norm{f}_{H^{-1}(\Omega)}^2=\sum_{j\ge1}\lambda_j^{-1}|\pair*{f}{\phi_j}_{H^{-1}(\Omega),H^1_0(\Omega)}|^2$.

The analysis of \cite{FCLZ16} rests on weighted observability inequalities for the heat equation ($k=0$) of the following form:
\begin{equation}\label{eq:FCLZspec}
\sum_{j\ge1} w_j^2 \abs{\pair{z}{\phi_j}_{L^2(\Omega)}}^2\les \int_0^T \norm*{\indic_\omega p(t,\cdot)}_{L^2(\Omega)}^2\,dt,
\qquad \forall z\in L^2(\Omega),
\end{equation}
where $w_j>0$, $j \geq 1$, are weights.

It is well known that the heat equation is null observable in any time $T>0$ \cite{LR95,FI96}.
In other words, \eqref{eq:FCLZspec} is true with $w_j=e^{-\lambda_jT}$.
Actually more holds: as reminded in \cite[Lemma 2]{FCLZ16}, there exists $r>0$ such that \eqref{eq:FCLZspec} holds with $w_j=e^{-r\sqrt{\lambda_j}} \ges e^{-\lambda_jT}$, and it is this weight that is used therein.

As in that paper, we consider the following class of kernels:
\begin{equation}\label{eq:FCLZsum}
\norm*{k}_w^2:=\sum_{m\ge1} \frac{1}{\lambda_m w_m^2} \left(\sum_{j\ge1}\frac{\abs{k_{mj}}^2}{\lambda_j}\right)<+\infty,
\qquad k_{mj}:=\int_{\Omega\times\Omega}k(x,\xi)\,\phi_m(x)\,\phi_j(\xi)\,dx\,d\xi.
\end{equation}
The purpose of this class is to bound the nonlocal term
\begin{equation}\label{obsheat:def f}
f(t,x):=-\int_\Omega k(\xi,x)\,p(t,\xi)\,d\xi,
\end{equation}
by the left-hand side of \eqref{eq:FCLZspec} (see, e.g., \cite[(21)]{FCLZ16} for the computations with the weights $w_j=e^{-r\sqrt{\lambda_j}}$):
\begin{equation}\label{obsheat:estim f}
\norm*{f}_{L^2(0,T;H^{-1}(\Omega))}^2\le\frac{1}{2}\norm*{k}_w^2\sum_{j\ge1}w_j^2\abs{\pair{z}{\phi_j}_{L^2(\Omega)}}^2.
\end{equation}

Let us point out that, a priori, whether \eqref{eq:FCLZspec} holds for given weights depends on $T$, but we will see that this is not the case (Lemma~\ref{lem:FCLZR} below).

The authors then proved in \cite[Theorem 3]{FCLZ16} the null controllability of the nonlocal heat equation \eqref{eq:FCLZsyst} in time $T$ for kernels $k$ satisfying \eqref{eq:FCLZsum}, with weights for which \eqref{eq:FCLZspec} holds at the time $T$, and under the following analyticity assumption on the integral term:
\begin{equation}\label{hyp FCLZ}
x\mapsto\int_\Omega k(\xi,x)\,z(\xi)\,d\xi\text{ is analytic in }\Omega,\qquad\forall z\in L^2(\Omega).
\end{equation}
They did it for the particular weights $w_j=e^{-r\sqrt{\lambda_j}}$, but the proof is the same for general weights.
We use general weights since this makes clear that the weights $w_j=e^{-\lambda_jT}$ of the null observability are allowed as well.

\subsubsection{Main result and comments}

The main result of Section~\ref{sec:compFCLZ} is the following extension of \cite[Theorem 3]{FCLZ16}.

\begin{proposition}[control of the nonlocal heat equation]
\label{prop:FCLZ}
Assume that \eqref{eq:FCLZspec} holds at some time, with some weights $w_j>0$, $j\ge1$, and let $k \in L^2(\Omega\times\Omega)$ satisfy \eqref{eq:FCLZsum}.
Then, we have:
\begin{enumerate}
\item\label{FCLZ:i1}
The set $\sigma_F$ of unobservable eigenvalues of the pair $(A^*,B^*)$ given by \eqref{eq:def AB nonloca heat} is finite and, for every $\lambda\in\sigma_F$, the subspace $E_\lambda$ of unobservable generalized eigenvectors of $(A^*,B^*)$ associated with $\lambda$ is finite dimensional.

\item\label{FCLZ:i2}
For every $T>0$ and $y^0\in L^2(\Omega)$, there exists $u\in L^2((0,T)\times\Omega)$ such that the solution $y$ of \eqref{eq:FCLZsyst} satisfies $y(T,\cdot)=0$ if, and only if,
\[
y^0\in E_\perp,\qquad\text{where }E:=\bigoplus_{\lambda\in\sigma_F}E_\lambda.
\]
In particular, if $(A^*,B^*)$ satisfies the Fattorini--Hautus test, then \eqref{eq:FCLZsyst} is null controllable in every time $T>0$.
\end{enumerate}
\end{proposition}

The proof is given in Section~\ref{sec proof obsheat} below.
We first give some comments.

\begin{remark}[on the analyticity assumption]
The assumption \eqref{hyp FCLZ} implies the Fattorini--Hautus test, and it is used in \cite{FCLZ16} only to prove the \emph{unique continuation property} for the adjoint system.
This property is in fact equivalent to the test, by \cite[Corollary 3.3]{Fat66}, whose assumptions hold here.
Indeed, $A$ is a bounded perturbation of the Dirichlet Laplacian, so that it has a compact resolvent and, by the Keldysh completeness theorem, a complete system of generalized eigenvectors.
Thus, the proof of \cite{FCLZ16} combined with \cite{Fat66} gives the implication from the test to the null controllability, as already observed in \cite{LZ18}.
By contrast, the approach of Theorems~\ref{thm:obs} and~\ref{thm:main} requires neither the completeness of the generalized eigenvectors nor the analyticity of the semigroup, and it describes in addition the null controllable states.
The null controllability under the test was also obtained in \cite[Proposition 4.2]{ELT22}, together with a stronger result: the reachable space of \eqref{eq:FCLZsyst} coincides with that of the heat equation ($k=0$).
\end{remark}

\begin{remark}[the one-dimensional case]
In one space dimension, the null controllability of \eqref{eq:FCLZsyst} is settled in \cite{MT18}, by the moment method.
The result is established in \cite[Theorem 1.2]{MT18}, for kernels of the form $k(x,\xi)=g(x)h(\xi)$ with $h$ not identically zero on $\omega$.
However, the rank one structure is used there only to ensure the Fattorini--Hautus test, so that their arguments apply to every kernel $k\in L^2$ for which the test holds, without the summability condition \eqref{eq:FCLZsum}.
\end{remark}

\begin{remark}[kernels producing unobservable eigenvectors]
When the Fattorini--Hautus test fails, there is no result in \cite{FCLZ16,MT18}, and \cite{ELT22} only shows that the reachable space then differs from that of the heat equation.
Yet, it is not difficult to construct examples of kernels satisfying \eqref{eq:FCLZsum} for which the test fails.
Indeed, take a kernel of the form $k(x,\xi)=g(x)h(\xi)$ and let $\lambda \in \R$ be arbitrary.
We look for a function $\psi\in H^2(\Omega)\cap H^1_0(\Omega)$, not identically zero, such that
\[
\left\{
\begin{aligned}
&\Delta\psi(x)-\int_\Omega k(\xi,x)\,\psi(\xi)\,d\xi=\lambda\psi(x)&&\text{in }\Omega,\\
&\psi=0&&\text{in }\omega.
\end{aligned}
\right.
\]
We have $\int_\Omega k(\xi,x)\,\psi(\xi)\,d\xi=I\,h(x)$ with $I:=\int_\Omega g(\xi)\,\psi(\xi)\,d\xi$.
If $I=0$, this means that $\Delta\psi=\lambda\psi$ in $\Omega$, and $\psi=0$ in $\omega$ gives $\psi=0$, by a classical unique continuation property.
Hence, we must look for an example with $I\neq0$, and then we have no choice for $h$ but $h:=(\Delta\psi-\lambda\psi)/I$.
Next, we have $k_{mj}=\pair*{g}{\phi_m}_{L^2(\Omega)}\pair*{h}{\phi_j}_{L^2(\Omega)}$, so that $\norm{k}_w^2=(\sum_{m\ge1} \lambda_m^{-1} w_m^{-2} |\pair*{g}{\phi_m}_{L^2(\Omega)}|^2)\norm{h}_{H^{-1}(\Omega)}^2$, and \eqref{eq:FCLZsum} holds if, and only if, the series in parentheses converges.
If $g=\phi_j$ for some $j\ge1$, the series reduces to one term, so that this holds for every choice of weights.
The condition $I\neq0$ then reads $\pair*{\psi}{\phi_j}_{L^2(\Omega)}\neq0$, and such a $j$ exists as soon as $\psi\neq 0$.
In summary, we take any $\psi\in C^\infty_c(\Omega\setminus\clos{\omega})$ which is not identically zero, assuming that $\Omega\setminus\clos{\omega}$ is not empty.
Finally, for this example, one can even compute $\sigma_F=\ens{\lambda}$ and $E_\lambda=\Span\ens{\psi}$, so that, by Proposition~\ref{prop:FCLZ}, the null controllable states are exactly the $y^0\in L^2(\Omega)$ orthogonal to this $\psi$.
\end{remark}

\subsubsection{Proofs}
\label{sec proof obsheat}

We start with the following lemma, which shows in particular that, if \eqref{eq:FCLZspec} holds at some time, then it holds at every time with the same weights.

\begin{lemma}\label{lem:FCLZR}
Let $(S_0(t))_{t \geq 0}$ be the heat semigroup, that is, the one generated by $A$ with $k=0$.
For every $T,T'>0$, there exists $c>0$ such that, for every $z\in L^2(\Omega)$,
$$
\int_0^T\norm*{\indic_\omega S_0(t)^*z}_{L^2(\Omega)}^2\,dt
\leq c\int_0^{T'}\norm*{\indic_\omega S_0(t)^*z}_{L^2(\Omega)}^2\,dt.
$$
\end{lemma}

Note that this is nothing but the time invariance of the reachable space for the heat equation, by duality.

\begin{proof}
For $T\le T'$, this is obvious.
If $T>T'$, we split the integral over $(0,T)$ at $t=T'$:
\begin{align*}
\int_{T'}^{T}\norm*{\indic_\omega S_0(t)^*z}_{L^2(\Omega)}^2\,dt
&=\int_0^{T-T'}\norm*{\indic_\omega S_0(s)^*S_0(T')^*z}_{L^2(\Omega)}^2\,ds\\
&\les\norm*{S_0(T')^*z}_{L^2(\Omega)}^2
\les\int_0^{T'}\norm*{\indic_\omega S_0(t)^*z}_{L^2(\Omega)}^2\,dt.
\end{align*}
The first identity follows from the semigroup property, the first estimate from the boundedness of $S_0$ on $[0,T]$, and the last estimate is the null observability of the heat equation in time $T'$.
\end{proof}

\begin{proof}[Proof of Proposition~\ref{prop:FCLZ}]
Fix $T>0$.
We prove property~\ref{obs:CU} of Theorem~\ref{thm:obs}, with $T$ in place of $T_0$ and $D$ as the full state space.
More precisely, we show that there exist $R=KL$, with $K$ compact, such that
\begin{align}
\label{eq:Ha+}
\norm*{S(T)^*z}_{Y}&\les \norm*{\Phi_{T}^*z}_{U_{T}}+\norm*{Rz}_{V_2}+\norm*{Lz}_{V_1},
\\
\label{eq:Hbp+}
\norm*{Lz}_{V_1}&\les \norm*{\Phi_{T}^*z}_{U_{T}}+\norm*{Rz}_{V_2},
\end{align}
for every $z \in Y$.
Inserting \eqref{eq:Hbp+} into \eqref{eq:Ha+}, and then adding \eqref{eq:Hbp+}, gives \eqref{eq:single estim} (it is in fact equivalent to \eqref{eq:Ha+} and \eqref{eq:Hbp+}).
Since $T$ is arbitrary, the conclusion then follows from Theorems~\ref{thm:obs} and~\ref{thm:main}.
The estimates below are essentially those of the proof of \cite[Theorem 3]{FCLZ16}, organized so as to identify the splitting $R=KL$ required by Theorem~\ref{thm:obs}.

Let $z\in L^2(\Omega)$, and let $\zeta:=\varphi-p$.
It solves, with $f$ given by \eqref{obsheat:def f},
\begin{equation}\label{eq:FCLZzeta}
\left\{
\begin{aligned}
&-\zeta_t(t,x)-\Delta\zeta(t,x)+\int_\Omega k(\xi,x)\,\zeta(t,\xi)\,d\xi=f(t,x)&&\text{in }(0,T)\times\Omega,\\
&\zeta=0&&\text{on }(0,T)\times\partial\Omega,\\
&\zeta(T,\cdot)=0&&\text{in }\Omega.
\end{aligned}
\right.
\end{equation}

\emph{Estimate of $\varphi(0,\cdot)$.}
Since $\varphi=p+\zeta$, we have
\[
\norm*{\varphi(0,\cdot)}_{L^2(\Omega)}^2\le2\norm*{p(0,\cdot)}_{L^2(\Omega)}^2+2\norm*{\zeta(0,\cdot)}_{L^2(\Omega)}^2.
\]
Estimating the first term by the classical null observability of the heat equation, then writing $p=\varphi-\zeta$, we obtain
\begin{equation}\label{heat:first estim}
\norm*{\varphi(0,\cdot)}_{L^2(\Omega)}^2\les\int_0^T\norm*{\indic_\omega\varphi(t,\cdot)}_{L^2(\Omega)}^2\,dt+\int_0^T\norm*{\zeta(t,\cdot)}_{L^2(\Omega)}^2\,dt+\norm*{\zeta(0,\cdot)}_{L^2(\Omega)}^2.
\end{equation}

\emph{The remainder and its splitting.}
In view of the previous estimate, we introduce
\[
R:L^2(\Omega)\to L^2(0,T;L^2(\Omega)),\qquad Rz:=\zeta.
\]
It has the natural splitting $R=KL$, with
\[
L:L^2(\Omega)\to L^2(0,T;H^{-1}(\Omega)),\qquad Lz:=f,
\]
where $f$ is given by the integral term \eqref{obsheat:def f}, and $K:L^2(0,T;H^{-1}(\Omega))\to L^2(0,T;L^2(\Omega))$ is the map which associates with \emph{any} $f\in L^2(0,T;H^{-1}(\Omega))$ the solution of \eqref{eq:FCLZzeta}.
Let us recall the classical energy estimate
$$
\sup_{t\in[0,T]}\norm*{\zeta(t,\cdot)}_{L^2(\Omega)}^2+\norm*{\zeta}_{L^2(0,T;H^1_0(\Omega))}^2\les\norm*{f}_{L^2(0,T;H^{-1}(\Omega))}^2,
$$
obtained by multiplying the first equation of \eqref{eq:FCLZzeta} by $\zeta$, integrating over $(t,T)\times\Omega$, and using the Young inequality and the Gronwall lemma.
In particular, the term $\zeta(0,\cdot)$ in \eqref{heat:first estim} is estimated by $f=Lz$, so that \eqref{heat:first estim} gives the first estimate \eqref{eq:Ha+}.
The energy estimate and the first equation of \eqref{eq:FCLZzeta} also give the same bound for $\zeta_t$ in $L^2(0,T;H^{-1}(\Omega))$.
Thus $K$ is bounded from $L^2(0,T;H^{-1}(\Omega))$ into $\ens*{\zeta\in L^2(0,T;H^1_0(\Omega))\st\zeta_t\in L^2(0,T;H^{-1}(\Omega))}$, which is compactly embedded in $L^2(0,T;L^2(\Omega))$ by the Aubin--Lions lemma, so that $K$ is compact.

\emph{Estimate \eqref{eq:Hbp+}.}
As mentioned before the proof, the condition \eqref{eq:FCLZsum} gives the estimate \eqref{obsheat:estim f}.
Therefore, using \eqref{eq:FCLZspec}, which holds at the time $T$ with the same weights $w_j$ by Lemma~\ref{lem:FCLZR}, and then writing $p=\varphi-\zeta$, we obtain
\[
\norm*{f}_{L^2(0,T;H^{-1}(\Omega))}^2\les\int_0^T\norm*{\indic_\omega\varphi(t,\cdot)}_{L^2(\Omega)}^2\,dt+\int_0^T\norm*{\zeta(t,\cdot)}_{L^2(\Omega)}^2\,dt.
\]
This gives the second estimate \eqref{eq:Hbp+}.
\end{proof}


\subsection{Simultaneous exact--null controllability}\label{sec:Tuc}

We now turn to a problem of a different nature: two uncoupled systems driven by one and the \emph{same control}.
Problems of this kind have a long history, going back to \cite{Rus86} and \cite[Chapitre V]{Lio88}, see also \cite{TW00}.
Here, we revisit and extend the very recent work \cite{Tuc26}, where the first system is required to reach an arbitrary state, while the second one is driven to zero.
Other mixed requirements of this kind are also considered in \cite{Lis26}.

\subsubsection{A motivating example}\label{sec:Tucmotiv}

We begin with the motivating example of \cite{Tuc26}: a wave equation and an advection--diffusion equation driven by the same boundary control,
\begin{equation}\label{syst tuc26}
\left\{
\begin{aligned}
&w_{tt}-\Delta w=0 , &&\theta_t-\Delta\theta-b(x)\cdot\nabla\theta-c(x)\,\theta=0 , &&\qquad\text{in }(0,T)\times\Omega,\\
&w=u\indic_\Gamma, &&\theta=u\indic_\Gamma , &&\qquad\text{on }(0,T)\times\partial\Omega,\\
&\left(w,w_t\right)(0,\cdot)=\left(w^0,w^1\right) , &&\theta(0,\cdot)=\theta^0 , &&\qquad\text{in }\Omega.
\end{aligned}
\right.
\end{equation}
Here $\Omega\subset\R^N$ is a bounded connected open set with boundary of class $C^2$ and $\Gamma\subset\partial\Omega$ is a nonempty relatively open subset; the coefficients satisfy $b\in W^{1,\infty}(\Omega)^N$ and $c \in L^\infty(\Omega)$, and the control is $u\in L^2((0,T)\times\partial\Omega)$.
The state spaces are $L^2(\Omega)\times H^{-1}(\Omega)$ for the first equation and $H^{-1}(\Omega)$ for the second, these being the spaces making the system with Dirichlet boundary control well-posed \cite[Proposition 10.9.1 and Theorem 10.8.3]{TW09}.
In this setting, it is shown in \cite[Proposition 2.3]{Tuc26} that, if $\Gamma$ satisfies the geometric control condition and if the advection--diffusion semigroup is exponentially stable, which holds for instance if $\dive b=0$ and $c\le0$ \cite[Remark 2.4]{Tuc26}, then for every $T$ larger than the geometric control time there is one control $u$ steering $(w,w_t)$ to an arbitrary target and $\theta$ to rest at time $T$.

This statement is obtained in \cite{Tuc26} from an abstract theorem on two systems driven by a common control, and it is that theorem which we extend here.
The main difference lies in the uniqueness part: the simultaneous approximate controllability, established there by an analytic continuation argument, is replaced by the Fattorini--Hautus test, which can often be decided by a direct computation on the spectra.
Two of the questions left open in \cite[Section 7]{Tuc26} are answered along the way.
The structure of the problem, and its relevance for systems of partial differential equations, were uncovered in \cite{Tuc26}.

\subsubsection{The abstract result}\label{sec:Tucabstract}

\paragraph{Setting.}
We now place ourselves in the abstract setting.
Let $Y_1,Y_2,U$ be complex reflexive Banach spaces and, for $k=1,2$, let $A_k$ generate a $C_0$-semigroup $(S_k(t))_{t\ge0}$ on $Y_k$ and let $B_k\in\LL(U,D(A_k^*)')$ be a control operator.
Since $Y_k$ is reflexive, \eqref{hyp measurable} holds for $(A_k,B_k)$ (see Section~\ref{sec:control}), and we let $\Phi_{k,T}$ be the control map of $(A_k,B_k)$, $k=1,2$.
We take $U_T=L^p(0,T;U)$ with $1<p<\infty$.
We assume that $\Phi_{k,T}(U_T) \subset Y_k$ for $k=1,2$.
Finally, we define
$$
\sigma_F:=\ens{\lambda\in\C\st A_1^*z_1=\lambda z_1,\quad A_2^*z_2=\lambda z_2,\quad B_1^*z_1+B_2^*z_2=0, \quad\text{for some }(z_1,z_2)\neq0},
$$
and
$$
E_\lambda:=\ens{(z_1,z_2)\st
\begin{aligned}
&(\lambda-A_1^*)^mz_1=0,\quad(\lambda-A_2^*)^mz_2=0, \quad\text{for some }m\in\N,\\
&B_1^*(\lambda-A_1^*)^jz_1+B_2^*(\lambda-A_2^*)^jz_2=0, \quad\forall j<m.
\end{aligned}}.
$$
The use of the same notation as in Definition~\ref{def:unobs} will be justified below.
Above, the vectors involved are implicitly assumed to belong to the appropriate domains.

The notion of controllability we are interested in here is the following one \cite{Tuc26}.

\begin{definition}[simultaneous exact--null controllability]\label{def:simul}
We say that the two systems $(A_1,B_1)$ and $(A_2,B_2)$ are \emph{simultaneously exact--null controllable} in time $T$ if, for every $y^0_1,y_1^T \in Y_1$ and $y^0_2\in Y_2$, there exists $u\in U_T$ such that the corresponding solutions $y_1$ and $y_2$ of the two systems satisfy $y_1(T)=y_1^T$ and $y_2(T)=0$.
\end{definition}

\paragraph{Main result and comments.}
The main result of Section~\ref{sec:Tuc} is the following.

\begin{theorem}[simultaneous exact--null controllability]
\label{thm:Tuc}
Assume that:
\begin{enumerate}
\item\label{Tuc:i1} $(A_1,B_1)$ is exactly controllable in some time $T_1$, and $\Ima S_1(t)=Y_1$ for every $t\ge0$;

\item\label{Tuc:i2} $(A_2,B_2)$ is null controllable in some time $T_2$, and
$S_2(\varepsilon)^*\in\LL(Y_2')$ is compact for some $\varepsilon>0$.
\end{enumerate}
Then, we have:
\begin{enumerate}
\item\label{Tuc:c1}
The set $\sigma_F$ defined above is finite and, for every $\lambda\in\sigma_F$, the subspace $E_\lambda$ is finite dimensional.

\item\label{Tuc:c2}
Let $T>T_1+T_2+\varepsilon$.
For every $y_1^0,y_1^T\in Y_1$ and $y_2^0\in Y_2$, there exists $u\in U_T$ such that $y_1(T)=y_1^T$ and $y_2(T)=0$ if, and only if,
\[
\left(y_1^T-S_1(T)y_1^0,\,-S_2(T)y_2^0\right)\in E_\perp,\qquad\text{where }E:=\bigoplus_{\lambda\in\sigma_F}E_\lambda.
\]
In particular, $(A_1,B_1)$ and $(A_2,B_2)$ are simultaneously exact--null controllable in every time $T>T_1+T_2+\varepsilon$ if the Fattorini--Hautus test holds: for every $\lambda\in\C$, $z_1\in D(A_1^*)$ and $z_2\in D(A_2^*)$,
\begin{equation}\label{FH for tuc26}
A_1^*z_1=\lambda z_1,\quad A_2^*z_2=\lambda z_2,\quad B_1^*z_1+B_2^*z_2=0
\ \Longrightarrow\ z_1=z_2=0 .
\end{equation}
\end{enumerate}
\end{theorem}

The proof is given after some comments on this result.

\begin{remark}[joint eigenvalues]
\label{rem:FH at joint}
The Fattorini--Hautus test can only fail at joint eigenvalues, because each system is assumed controllable.
Indeed, if for instance $\lambda \in \sigma_p(A_1^*)$ but $\lambda \not\in \sigma_p(A_2^*)$, then $z_2=0$ and $z_1 \in \ker(\lambda-A_1^*) \cap \ker B_1^*$ which is reduced to $\ens*{0}$ by controllability of $(A_1,B_1)$.
\end{remark}

\begin{remark}[surjectivity of $S_1$]
\label{rem:EC implies surj}
We recall that, for a bounded control operator, Hilbert spaces and $L^2$ controls, exact controllability implies the surjectivity of the semigroup at every time, see e.g.\ \cite[Theorem~1]{Zwa13}.
In that case, asking for the surjectivity of $S_1(t)$ in \ref{Tuc:i1} is redundant.
\end{remark}

\begin{remark}[comparison with \cite{Tuc26}]
\label{rem:compa tuc26}
Besides having a stronger conclusion and being set in Banach spaces, Theorem~\ref{thm:Tuc} differs from the results of \cite{Tuc26} in one crucial respect: \emph{the Fattorini--Hautus test}.
In particular, our result implies at once \cite[Corollary 4.5]{Tuc26} and \cite[Theorem 2.1]{Tuc26}.
Indeed, in \cite[Corollary 4.5]{Tuc26} it is assumed that the systems $(A_1,B_1)$ and $(A_2,B_2)$ are simultaneously \emph{approximately controllable}, whereas here we only assume the Fattorini--Hautus test, which is necessary and usually easier to check: it holds as soon as $A_1^*$ and $A_2^*$ have \emph{no common eigenvalue} (Remark~\ref{rem:FH at joint}).
This gives \cite[Theorem 2.1]{Tuc26}: it is assumed therein that $A_1$ generates a \emph{bounded group} and $A_2$ an \emph{exponentially stable} semigroup, so that $\sigma(A_1)\subset i\R$ and $\sigma(A_2)\subset\C_-$ are disjoint.
The simultaneous approximate controllability is proved in \cite[Section 5]{Tuc26} under the additional assumption that $S_2$ is analytic, which is not needed here.
Finally, this also answers one of the questions left open in \cite[Section 7]{Tuc26}: the exponential stability of $A_2$ may be weakened to $\sigma(A_2)\cap i\R=\emptyset$, since the two spectra are then again disjoint.
\end{remark}

\begin{example}[the motivating example]\label{ex:Tuc drifts}
As in \cite{Tuc26}, the assumptions of Theorem~\ref{thm:Tuc} hold for \eqref{syst tuc26} in the setting of Section~\ref{sec:Tucmotiv}.
Indeed, the wave operator is skew-adjoint on $Y_1=L^2(\Omega)\times H^{-1}(\Omega)$, hence generates a group \cite[Section 10.9]{TW09}, and the wave equation is exactly controllable in time $T_1$ if $\partial\Omega$ is smooth and $(\Gamma,T_1)$ satisfies the geometric control condition \cite[Corollary 4.10]{BLR92}.
The advection--diffusion operator has compact resolvent and generates an analytic semigroup, hence compact at positive times \cite[Diagram (II.4.26) and Theorem II.4.29]{EN00}, and the equation is null controllable in every time \cite[Theorem 2.3]{FI96}.
Thus, the two equations are simultaneously exact--null controllable in every time $T>T_1$ if (and only if, Proposition~\ref{prop:FHnec}) \eqref{FH for tuc26} holds, where
\[
A_1^*=\begin{pmatrix} 0 & \Delta \\ \Id & 0 \end{pmatrix},\quad B_1^*=\begin{pmatrix} 0 & -\partial_\nu \end{pmatrix},\qquad
A_2^*=\Delta-b\cdot\nabla+c-\dive b,\quad B_2^*=-\partial_\nu .
\]
Since $A_1^*$ is skew-adjoint and invertible, this test can only fail at a common eigenvalue in $i\R\setminus\ens{0}$ (Remark~\ref{rem:FH at joint}).

\emph{Gradient fields: the test always holds.}
If $b=\nabla\beta$, as always when $N=1$, then $A_2^*\psi=e^{\beta}\dive(e^{-\beta}\nabla\psi)+(c-\Delta\beta)\psi$ is symmetric in $L^2(\Omega)$ with weight $e^{-\beta}$, so that its eigenvalues are real and the test holds whatever $c$, even if $0\in\sigma(A_2)$, a case not covered by \cite{Tuc26}.

\emph{Rotations: the test may fail.}
Let $\Omega$ be the unit disk, $\Gamma=\partial\Omega$, $\mu$ the first positive zero of the Bessel function $J_1(s):=\sum_{p\ge0}\frac{(-1)^p}{p!\,(p+1)!}(\frac s2)^{2p+1}$, and $\phi(x):=e^{i\alpha}J_1(\mu r)$, where $(r,\alpha)$ are the polar coordinates of $x$, so that $-\Delta\phi=\mu^2\phi$ in $\Omega$ and $\phi=0$ on $\partial\Omega$.
For $b(x):=\mu(-x_2,x_1)$, so that $b\cdot\nabla=\mu\,\partial_\alpha$ and $\dive b=0$, and $c:=\mu^2$, we get $A_2^*\phi=-i\mu\phi$, so that $z_1:=(-i\mu\phi,\phi)^\tr$ and $z_2:=-\phi$ are nonzero and satisfy $A_1^*z_1=-i\mu z_1$, $A_2^*z_2=-i\mu z_2$ and $B_1^*z_1+B_2^*z_2=0$: the test \eqref{FH for tuc26} fails, and the two equations are never simultaneously exact--null controllable.
\end{example}

\paragraph{Proof of Theorem~\ref{thm:Tuc}.}
First of all, we recast the two systems as a single one, on the state space $Y:=Y_1\times Y_2$, with the generator and the control operator
$$
A:=\diag(A_1,A_2)=
\begin{pmatrix}
A_1 & 0 \\
0 & A_2
\end{pmatrix},
\qquad
B:=
\begin{pmatrix}
B_1 \\
B_2
\end{pmatrix}.
$$
Then $S(t)=\diag(S_1(t),S_2(t))$ and $\Phi_Tu=(\Phi_{1,T}u,\Phi_{2,T}u)$, so that the regularity gain \eqref{eq:RGC} holds, since $\Phi_{k,T}(U_T)\subset Y_k$ for $k=1,2$, and $(U_T)_{T>0}$ is a control class, by Example~\ref{ex:UT}.
On the dual side, $Y'=Y_1'\times Y_2'$, $A^*=\diag(A_1^*,A_2^*)$, and $B^*(z_1,z_2)=B_1^*z_1+B_2^*z_2$, so that $\sigma_F$ and $E_\lambda$ are indeed those of Definition~\ref{def:unobs} for the pair $(A^*,B^*)$, as announced above, and
\begin{equation}\label{eq:Tucsplit}
\Phi_T^*z=\Phi_{1,T}^*z_1+\Phi_{2,T}^*z_2,\qquad\forall z=(z_1,z_2)\in Y'.
\end{equation}

Thanks to the surjectivity of $S_1(T)$, simultaneous exact--null controllability in time $T$ is the null controllability of $(A,B)$ in time $T$.
Indeed, it asks that $(y_1^T-S_1(T)y_1^0,-S_2(T)y_2^0)\in\Ima\Phi_T$ for all $y_1^0,y_1^T\in Y_1$ and $y_2^0\in Y_2$, that is, $Y_1\times\Ima S_2(T)\subset\Ima\Phi_T$, and we have $Y_1\times\Ima S_2(T)=\Ima S(T)$, since $\Ima S_1(T)=Y_1$.
The surjectivity of $S_1(T)$ also gives that $S_1(T)^*$ is bounded below \cite[Theorem~2.20]{Bre11}, so that
\begin{equation}\label{eq:Tucsem}
\norm*{S(T)^*z}_{Y'}\simeq\norm*{z_1}_{Y_1'}+\norm*{S_2(T)^*z_2}_{Y_2'},\qquad\forall z=(z_1,z_2)\in Y'.
\end{equation}

Finally, since $U_T=L^p(0,T;U)$, the extension by zero $u\mapsto u\conca{T}0$ and the shift $v\mapsto0\conca{t}v$ are isometries from $U_T$ to $U_{T+t}$, for every $T,t>0$, and the composition property \eqref{eq:composition} gives $\Phi_{T+t}(u\conca{T}0)=S(t)\Phi_Tu$ and $\Phi_{T+t}(0\conca{t}v)=\Phi_Tv$.
Taking adjoints, we obtain, for every $T,t>0$ and $z\in Y'$,
\begin{equation}\label{eq:Tucmono}
\norm*{\Phi_T^*S(t)^*z}_{U_T'}\le\norm*{\Phi_{T+t}^*z}_{U_{T+t}'},
\qquad
\norm*{\Phi_T^*z}_{U_T'}\le\norm*{\Phi_{T+t}^*z}_{U_{T+t}'},
\end{equation}
and similarly for each system $(A_k,B_k)$.

We are now ready to prove the main result of this section.

\begin{proof}
We prove property~\ref{obs:CU} of Theorem~\ref{thm:obs}, with $D=Y'$, and then deduce the conclusion from Theorems~\ref{thm:obs} and~\ref{thm:main}.
For now, $T_0>0$ is a parameter, which the proof will constrain and which we fix at the end.

\emph{Estimate of $z_2$.}
We have
\begin{equation}\label{eq:proof Tuc est2}
\norm*{S_2(T_0)^*z_2}_{Y_2'}\les \norm*{\Phi_{2,T_0}^*z_2}_{U_{T_0}'}
\les \norm*{\Phi_{T_0}^*z}_{U_{T_0}'}+\norm*{\Phi_{1,T_0}^*z_1}_{U_{T_0}'}
\les \norm*{\Phi_{T_0}^*z}_{U_{T_0}'}+\norm*{z_1}_{Y_1'}.
\end{equation}
The first estimate is the null observability of $(A_2,B_2)$ in time $T_2$, and thus in any time $T_0 \geq T_2$ by \eqref{eq:Tucmono}; the second one is the splitting \eqref{eq:Tucsplit}; and the third one holds since $\Phi_{1,T_0}^*$ is bounded.

\emph{Estimate of $z_1$.}
Let $t \geq 0$ be another parameter, which the proof will also constrain and which we fix at the end.
We have
\begin{align}
\norm*{z_1}_{Y_1'}
&\simeq \norm*{S_1(t)^*z_1}_{Y_1'}
\les \norm*{\Phi_{1,T_1}^*S_1(t)^*z_1}_{U_{T_1}'}\notag\\
&\les \norm*{\Phi_{T_1}^*S(t)^*z}_{U_{T_1}'}+\norm*{\Phi_{2,T_1}^*S_2(t)^*z_2}_{U_{T_1}'}
\les \norm*{\Phi_{T_0}^*z}_{U_{T_0}'}+\norm*{\Phi_{2,T_1}^*S_2(t)^*z_2}_{U_{T_1}'},
\label{eq:proof Tuc est1}
\end{align}
the first equivalence is because $S_1(t)^*$ is bounded below, by assumption \ref{Tuc:i1};
the second one is the exact observability \eqref{eq:exactobs} of $(A_1,B_1)$ in time $T_1$, applied to $S_1(t)^*z_1$; the third one is the splitting \eqref{eq:Tucsplit} at the state $S(t)^*z$; the last one comes from \eqref{eq:Tucmono}, provided that $T_1+t\le T_0$.

\emph{Splitting.}
Combining \eqref{eq:Tucsem} with the estimates \eqref{eq:proof Tuc est2} and \eqref{eq:proof Tuc est1}, we have
\begin{equation}\label{eq:tuc estim S}
\norm*{S(T_0)^*z}_{Y'} \les \norm*{\Phi_{T_0}^*z}_{U_{T_0}'}+\norm*{Rz}_{U_{T_1}'},
\end{equation}
for every $z \in Y'$, where we introduce $R(z_1,z_2):=\Phi_{2,T_1}^* S_2(t)^*z_2$.
The definition of $R$ suggests a splitting between $\Phi_{2,T_1}^*$ and $S_2(t)^*$, and this is where the compactness of $S_2(\varepsilon)^*$ enters.
Assuming $t\ge\varepsilon$, we take $K:Y_2' \to U_{T_1}'$ and $L:Y' \to Y_2'$ defined by
\[
K:=\Phi_{2,T_1}^*S_2(\varepsilon)^*, \qquad
Lz:=S_2(t-\varepsilon)^*z_2,
\]
so that $R=KL$, the operator $K$ being compact.

\emph{Estimate of $L$.}
Assuming $t-\varepsilon\ge T_2$, we have
$$
\norm*{Lz}_{Y_2'}=\norm*{S_2(T_2)^*S_2(t-\varepsilon-T_2)^*z_2}_{Y_2'}
\les \norm*{\Phi_{2,T_2}^*S_2(t-\varepsilon-T_2)^*z_2}_{U_{T_2}'}
\les \norm*{\Phi_{2,t-\varepsilon}^*z_2}_{U_{t-\varepsilon}'},
$$
the first estimate is the null observability \eqref{eq:nullobs} of $(A_2,B_2)$ in time $T_2$, applied to $S_2(t-\varepsilon-T_2)^*z_2$; the second one is by \eqref{eq:Tucmono}.
We deduce
\begin{equation}\label{eq:tuc estim L}
\norm*{Lz}_{Y_2'}\les \norm*{\Phi_{t-\varepsilon}^*z}_{U_{t-\varepsilon}'}+\norm*{z_1}_{Y_1'}
\les \norm*{\Phi_{T_0}^*z}_{U_{T_0}'}+\norm*{S(T_0)^*z}_{Y'},
\end{equation}
the first estimate follows from the splitting \eqref{eq:Tucsplit}, since $\Phi_{1,t-\varepsilon}^*$ is bounded; the second one from \eqref{eq:Tucmono}, which applies since $t-\varepsilon\le T_0$ by $T_1+t\le T_0$, and \eqref{eq:Tucsem}.
Combining the estimates \eqref{eq:tuc estim L} and \eqref{eq:tuc estim S}, we obtain the estimate \eqref{eq:single estim} of property~\ref{obs:CU} of Theorem~\ref{thm:obs}.

\emph{Choice of the parameters.}
The three conditions met along the way, namely $T_0\ge T_2$, $T_1+t\le T_0$ and $t-\varepsilon\ge T_2$, are all satisfied, with the smallest possible $T_0:=T_1+T_2+\varepsilon$, by taking $t:=T_2+\varepsilon$.

\emph{Conclusion.}
By Theorem~\ref{thm:obs}, assumption \eqref{eq:primalH} holds, so that Theorem~\ref{thm:main} applies.
Its item~\ref{main:i2} gives item~\ref{Tuc:c1}, and its item~\ref{main:i3} gives $S(T)^{-1}\Ima\Phi_T=E_\perp$ for every $T>T_0$.
Let $T>T_0$, $y_1^0,y_1^T\in Y_1$ and $y_2^0\in Y_2$.
The required control exists if, and only if, $(y_1^T-S_1(T)y_1^0,-S_2(T)y_2^0)\in\Ima\Phi_T$.
Since $Y_1=\Ima S_1(T)$, we can write $y_1^T-S_1(T)y_1^0=-S_1(T)w$ for some $w\in Y_1$, so that the condition reads $-(w,y_2^0)\in S(T)^{-1}\Ima\Phi_T=E_\perp$.
Since $E_\perp=(S(T)^*E)_\perp=S(T)^{-1}(E_\perp)$, as at the end of the proof of Theorem~\ref{thm:main}, this is equivalent to $-S(T)(w,y_2^0)\in E_\perp$, that is, to $(y_1^T -S_1(T)y_1^0,-S_2(T)y_2^0)\in E_\perp$.
Finally, \eqref{FH for tuc26} is the Fattorini--Hautus test \eqref{eq:FH} for the pair $(A^*,B^*)$, so that, if it holds, then $E=\ens*{0}$ and this condition is always satisfied.
\end{proof}

\subsubsection{Several components}\label{sec:Tucseveral}

A second question left open in \cite[Section 7]{Tuc26} concerns the case of several components, where the simultaneous approximate controllability is identified as the main new issue.
Since Theorem~\ref{thm:Tuc} assumes no such property, we can answer it in part with the following result.

Let $n\ge2$, and let $(A_k,B_k)$, $k=1,\dots,n$, be systems as in Theorem~\ref{thm:Tuc}: $A_k$ generates a $C_0$-semigroup $(S_k(t))_{t\ge0}$ on a complex reflexive Banach space $Y_k$ and $B_k\in\LL(U,D(A_k^*)')$, with $U$ reflexive and $U_T=L^p(0,T;U)$, $1<p<\infty$.
We assume that $\Phi_{k,T}(U_T)\subset Y_k$ for every $k$.
All the components are driven by one and the same control.
A family of systems is said to be \emph{simultaneously exactly controllable} in time $T$ if, for every family of initial states and every family of targets, there is a single control steering each component to its target at time $T$; and \emph{simultaneously null controllable} in time $T$ if this holds with all the targets equal to $0$.

\begin{corollary}[several components]\label{cor:groups}
Let $1\le m\le n-1$, and assume that:
\begin{enumerate}
\item\label{cor:i1}
the systems $(A_k,B_k)$, $k\le m$, are simultaneously exactly controllable in some time $T_1>0$, and $\Ima S_k(t)=Y_k$ for every $t\ge0$ and $k\le m$;

\item\label{cor:i2}
the systems $(A_k,B_k)$, $k\ge m+1$, are simultaneously null controllable in some time $T_2>0$, and $S_k(\varepsilon)^*$ is compact for every $k\ge m+1$, for some $\varepsilon>0$;

\item
the Fattorini--Hautus test holds: for every $\lambda\in\C$ and $z_1\in D(A_1^*)$, \dots, $z_n\in D(A_n^*)$,
$$
\left\{
\begin{array}{r@{\;}c@{\;}l}
A_1^*z_1 & = & \lambda z_1,\\
 & \vdots & \\
A_n^*z_n & = & \lambda z_n,
\end{array}
\right.
\quad
B_1^*z_1+\dots+B_n^*z_n=0
\ \Longrightarrow\ z_1=\dots=z_n=0.
$$
\end{enumerate}
Then, for every $T>T_1+T_2+\varepsilon$, every $y^0_1\in Y_1$, \dots, $y^0_n\in Y_n$ and every $y^T_1\in Y_1$, \dots, $y^T_m\in Y_m$, there exists $u\in U_T$ such that the corresponding solutions satisfy
\[
y_1(T)=y^T_1,\ \dots,\ y_m(T)=y^T_m,\qquad y_{m+1}(T)=\dots=y_n(T)=0.
\]
\end{corollary}

\begin{proof}
Apply Theorem~\ref{thm:Tuc} to the two systems obtained by grouping the first $m$ and the last $n-m$ ones.
\end{proof}

By Remark~\ref{rem:FH at joint} applied to the grouped systems, the Fattorini--Hautus test holds as soon as the two groups share no eigenvalue:
\[
\sigma_p\left(A_k^*\right)\cap\sigma_p\left(A_l^*\right)=\emptyset,\qquad\text{for every } k\le m<l.
\]
The two simultaneous controllability assumptions, on the other hand, are genuinely stronger than the controllability of the individual components, and they are the delicate ones; see e.g.\ \cite{TW00} for simultaneous exact controllability and \cite[Section 6.2]{AKBGBdT14} for simultaneous null controllability.
With the Fattorini--Hautus test in place of the simultaneous approximate controllability, the remaining issue in the passage to several components is thus the verification of the simultaneous exact or null controllability within each group.


\subsection{Transport equation with nonlocal terms}\label{sec:nonloc t}

We now consider a transport equation controlled from an interior band, with \emph{nonlocal terms} in the equation and at the boundary.
Without the former, this is the Lotka--McKendrick system, controlled from a band in \cite{BIM01,HMT18,Mai19}.
This term, which arises in the linearization of the model of \cite{GM74}, couples at every time the uncontrolled regions to the whole state, so that the direct constructions of these works no longer apply.
Our main result in this section gives the \emph{minimal time} for null controllability, in every space $L^r$ with $1\le r<\infty$, \emph{including the natural space} $L^1$, as soon as the Fattorini--Hautus test holds, and, otherwise, it gives a \emph{spectral description} of the null controllable states after that time.

\subsubsection{Setting}\label{sec:nbsetting}

Inspired by the linearized model \cite[(3.10)]{GM74}, we investigate the control problem for the equation
\begin{equation}\label{eq:nbsys}
\left\{
\begin{aligned}
&\varrho_t+\varrho_x+\mu(x)\varrho=-k(x)\int_0^1\varrho(t,\xi)\,d\xi+v(t,x)\indic_{(a_1,a_2)}(x), &&(t,x)\in R_T,\\
&\varrho(t,0)=\int_0^1\beta(\xi)\varrho(t,\xi)\,d\xi, &&t\in(0,T),\\
&\varrho(0,x)=\varrho^0(x), &&x\in(0,1).
\end{aligned}
\right.
\end{equation}
Here $R_T:=(0,T)\times(0,1)$, $0\le a_1<a_2\le1$, and the coefficients are real valued.
In the terminology of \cite{GM74}, $x$ is the age, $\beta$ represents fertility, and $k$ the effect of the total population on mortality.
The parameter $\mu$ represents mortality and, as in \cite{BIM01,HMT18,Mai19}, we assume it satisfies
\begin{equation}\label{eq:nbmu}
\mu\ge0,\qquad \mu\in L^1(0,1-\delta)\, \text{ for every }\delta\in(0,1),\qquad \int_0^1\mu(x)\,dx=\infty,
\end{equation}
and we let $\pi(x):=\exp(-\int_0^x\mu(\xi)\,d\xi)$ for $x\in[0,1)$ and $\pi(1):=0$, which represents the probability to survive to the age $x$ \cite[(3.2)]{GM74}.
Thus $\pi$ is positive on $[0,1)$ and belongs to $W^{1,1}(0,1)$, with $\pi'=-\mu\pi$, so that $\int_0^1\mu\pi=1$, and $\pi$ is nonincreasing:
\begin{equation}\label{eq:nbmono}
\pi(x)\le\pi(\xi),\qquad \xi\le x.
\end{equation}
Our assumptions on the other coefficients are
\begin{equation}\label{eq:nbhyp}
\beta, \, \frac k\pi\in W^{1,1}(0,1),\qquad \beta\not\equiv 0 \text{ in } (a_1,a_2).
\end{equation}
In particular, $k\in L^\infty(0,1)$, since $\pi\le1$.

\paragraph{Removing the mortality term.}
From the properties of $\pi$, the change of unknowns
$$
y:=\frac{\varrho}{\pi},\qquad u:=\frac v\pi,\qquad h:=\frac k\pi,\qquad \gamma:=\beta\pi,
$$
transforms \eqref{eq:nbsys} into the following system, which has the same controllability properties:
\begin{equation}\label{eq:nbq}
\left\{
\begin{aligned}
&y_t+y_x=-h(x)\int_0^1\pi(\xi)y(t,\xi)\,d\xi+u(t,x)\indic_{(a_1,a_2)}(x), &&(t,x)\in R_T,\\
&y(t,0)=\int_0^1\gamma(\xi)y(t,\xi)\,d\xi, &&t\in(0,T),\\
&y(0,x)=y^0(x), &&x\in(0,1),
\end{aligned}
\right.
\end{equation}
with $y^0:=\varrho^0/\pi$.
Two quantities play an important role:
\begin{equation}\label{eq:def pn}
p(t):=\int_0^1\pi(\xi)y(t,\xi)\,d\xi,\qquad b(t):=\int_0^1\gamma(\xi)y(t,\xi)\,d\xi,
\end{equation}
which are, in the terminology of \cite{GM74}, the perturbations of the \emph{total population} and of the \emph{birth rate} around an equilibrium.

\paragraph{Functional setting.}
Let $r\in[1,\infty)$, and let $r'\in(1,\infty]$ be its conjugate exponent: $1/r+1/r'=1$.
We work in the setting
$$
Y=U:=L^r(0,1;\pi^r\,dx),\qquad \norm{y}_{Y}:=\norm{\pi y}_{L^r(0,1)},\qquad U_T:=L^r(0,T;U).
$$
The system \eqref{eq:nbq} reads $y'=Ay+Bu$, where
$$
(Ay)(x)=-y'(x) -h(x)\int_0^1\pi(\xi)y(\xi)\,d\xi,\qquad (Bu)(x)=u(x)\indic_{(a_1,a_2)}(x),
$$
with $D(A)=\ens*{y\in Y\st y'\in Y,\ y(0)=\int_0^1\gamma y}$.
Here, since $\pi\ge\pi(1-\delta)>0$ on $[0,1-\delta]$, every $y\in Y$ belongs to $L^r(0,1-\delta)$ for every $\delta\in(0,1)$, so that its derivative $y'$ in the sense of distributions on $(0,1)$ makes sense.
If $y,y'\in Y$, then $y\in W^{1,r}(0,1-\delta)\subset C^0([0,1-\delta])$ for every $\delta\in(0,1)$, so that the trace $y(0)$ makes sense as well, and finally $\int_0^1\gamma y$ is finite, since $\gamma y=\beta\, (y \pi)$, with $\beta$ bounded.

\begin{proposition}[the semigroup and its adjoint]\label{prop:nbsg}
Assume \eqref{eq:nbmu} and $\beta,k\in L^\infty(0,1)$.
\begin{enumerate}
\item\label{nbsg:gen}
The operator $A$ generates a $C_0$-semigroup $(S(t))_{t\ge0}$ on $Y$.

\item\label{nbsg:adj}
Identify $Y'$ with the space of the functions $z$ such that $z/\pi\in L^{r'}(0,1)$, through $\pair*{z}{y}_{Y',Y}:=\int_0^1z\,\conj{y}\,dx$.
Then
\begin{equation}\label{eq:nbadj}
(A^*z)(x)=z'(x)+z(0)\gamma(x)-\left(\int_0^1h(\xi)z(\xi) \,d\xi\right)\pi(x),
\end{equation}
with domain $D(A^*)=\ens{z\in Y'\st z'\in Y',\ z(1)=0}$.
\end{enumerate}
\end{proposition}

Note that the traces $z(0)$ and $z(1)$ make sense, since $\pi\le1$ implies that $z\in W^{1,r'}(0,1)\subset C^0([0,1])$ for every $z\in Y'$ with $z'\in Y'$.
We could not find a complete proof of this result in this generality in the literature, so we give one in Section~\ref{sec:nbother}.

Since $B$ is bounded and $U_T=L^r(0,T;U)$, the system \eqref{eq:nbq} fits into the framework of Section~\ref{sec:setting} (Remark~\ref{rem:bounded B} and Example~\ref{ex:UT}).
For $y^0\in Y$ and $u\in U_T$, we call the solution of \eqref{eq:nbq} in the sense of Definition~\ref{def:solution} its \emph{Duhamel solution}, to distinguish it, a priori, from the \emph{solution along the characteristics} introduced in Section~\ref{sec:nbproof} below.
By \eqref{duhamel} and Remark~\ref{rem:bounded B}, it is given by $y(t)=S(t)y^0+\int_0^tS(t-s)Bu(s)\,ds$, where the integral is a Bochner integral in $Y$, and it belongs to $C^0([0,T];Y)$, that is, $\pi y\in C^0([0,T];L^r(0,1))$.
In particular, $p,b\in C^0([0,T])$, by the H\"older inequality.

With the identification of Proposition~\ref{prop:nbsg}, $(B^*z)(x)=z(x)\indic_{(a_1,a_2)}(x)$, so that the set $\sigma_F$ and the spaces $E_\lambda$ of Theorem~\ref{thm:main} can be computed with \eqref{eq:nbadj} and this formula, and, for a subset $N\subset Y'$, the pre-annihilator $N_\perp$ is the set of all $y\in Y$ such that $\int_0^1z\,\conj{y}\,dx=0$ for every $z\in N$.

Let
\begin{equation}\label{equ:def T0}
T_0:=\max\ens{a_1,1-a_2},
\end{equation}
which is the time needed by every characteristic to enter the band or to reach $x=1$.

We first collect some properties of \eqref{eq:nbq}.
They show that its null controllability is a genuine problem, that it cannot be obtained through exact controllability and that no control time smaller than $T_0$ is possible.

\begin{proposition}[first properties]
\label{prop:nbfirst}
Assume \eqref{eq:nbmu} and \eqref{eq:nbhyp}.
\begin{enumerate}
\item\label{nbfirst:nonzero}
If $(1-\int_0^1\gamma)(1+\int_0^1\pi H)+\int_0^1\pi\int_0^1\gamma H\le0$, where $H(x):=\int_0^xh$, then $A$ has an eigenvalue $\lambda\ge0$, so that $S(t)$ does not tend to $0$.

\item\label{nbfirst:EC}
If $(a_1,a_2)\ne(0,1)$, then \eqref{eq:nbq} is never exactly controllable, whatever $T$.

\item\label{nbfirst:T0}
Equation \eqref{eq:nbq} is not null controllable in any time $T<T_0$.
\end{enumerate}
\end{proposition}

For the sake of readability, its proof is given in Section~\ref{sec:nbother}.

\subsubsection{Main result and comments}\label{sec:nbresult}

The main result of Section~\ref{sec:nonloc t} is the following.

\begin{theorem}[control of the nonlocal transport equation]\label{thm:nb}
Assume \eqref{eq:nbmu} and \eqref{eq:nbhyp}.
Let $r \in [1,\infty)$ and let $T_0$ be defined by \eqref{equ:def T0}.
Then, we have:
\begin{enumerate}
\item\label{nb:i1}
The set $\sigma_F$ of unobservable eigenvalues of the pair $(A^*,B^*)$ given by \eqref{eq:nbadj} is finite and, for every $\lambda\in\sigma_F$, the subspace $E_\lambda$ of unobservable generalized eigenvectors of $(A^*,B^*)$ associated with $\lambda$ is finite dimensional.

\item\label{nb:i2}
For every $T>T_0$ and $\varrho^0\in L^r(0,1)$, there exists $v\in L^r(R_T)$ such that the solution $\varrho$ of \eqref{eq:nbsys} satisfies $\varrho(T,\cdot)=0$ if, and only if, $\varrho^0$ satisfies
\[
\int_0^1\frac{z(x)}{\pi(x)}\,\varrho^0(x)\,dx=0\quad\text{for every }z\in E,\qquad\text{where }E:=\bigoplus_{\lambda\in\sigma_F}E_\lambda.
\]
In particular, if $(A^*,B^*)$ satisfies the Fattorini--Hautus test, then \eqref{eq:nbsys} is null controllable in every time $T>T_0$.
\end{enumerate}
\end{theorem}

The proof is given in Section~\ref{sec:nbproof}.

\begin{remark}[minimal control time]\label{rem:nbminimal}
If $(A^*,B^*)$ satisfies the Fattorini--Hautus test, then, by Theorem~\ref{thm:nb} and Proposition~\ref{prop:nbfirst}, the time $T_0$ is the \emph{minimal null control time} of \eqref{eq:nbsys}, that is, the infimum of the times $T$ for which \eqref{eq:nbsys} is null controllable in time $T$.
\end{remark}

\begin{remark}[failure of \eqref{eq:h2}]
\label{rem:nbh2}
As announced in Remark~\ref{rem:HO22}, the estimate in \eqref{eq:h2} fails for the equation \eqref{eq:nbq} in the Hilbert case $r=2$, as soon as $t_1<1\le t_2$.
It suffices to consider the sequence $z_j:=\sqrt j\,\pi\indic_{(t_1,t_1+1/j)}$.
This is clear for the pure transport equation, that is, for $k=\beta=0$, since the left-hand side of the estimate is bounded away from zero, whereas its right-hand side is $O(1/j)$.
The nonlocal terms do not destroy this, since they perturb the semigroup by a compact operator, whose values on $z_j$ tend to zero, as $z_j$ tends weakly to zero.
Consequently, \eqref{eq:h2} fails whenever $T_0<1$, and \cite[Theorem 7.7]{HO22-JDE} cannot reach the minimal null control time $T_0=\max\ens{a_1,1-a_2}<1$ of Remark~\ref{rem:nbminimal}.
\end{remark}

The set $\sigma_F$ can be computed explicitly.

\begin{proposition}[the Fattorini--Hautus test]\label{prop:nbFH}
Assume \eqref{eq:nbmu} and \eqref{eq:nbhyp}.
\begin{enumerate}
\item\label{nbFH:indep}
If $\beta$ is not constant on $(a_1,a_2)$, then $\sigma_F=\emptyset$.

\item\label{nbFH:dep}
If $\beta$ is equal to a constant $\beta_0$ on $(a_1,a_2)$, then $\lambda\in\sigma_F$ if, and only if,
\begin{equation}\label{eq:nbFHeq}
\psi(a_1)=\psi(1)=0,\qquad \int_0^1h(x)\psi(x) \, dx=\beta_0,
\end{equation}
where $\psi$ is the solution of $\psi'=\lambda\psi-(\beta-\beta_0)\pi$ on $(0,1)$ with $\psi(0)=1$.
In that case, $\ker(\lambda-A^*)\cap\ker B^*=\Span\ens*{\psi}$.
In particular, $\sigma_F=\emptyset$ if $a_1=0$, $k=0$ or $\beta$ constant on $(0,1)$.
\end{enumerate}
\end{proposition}

The proof, a direct computation on the eigenvalue equation, is postponed to Section~\ref{sec:nbother}.

The test may also fail.
Let $a_1>0$ and $\beta=\beta_0$ on $(a_1,1)$.
For $\lambda=0$, the function $\psi$ of Proposition~\ref{prop:nbFH} is $\psi(x)=1-\int_0^x(\beta-\beta_0)\pi$, which is constant on $[a_1,1]$, and thus vanishes there as soon as $\psi(a_1)=0$.
Hence, if we choose $\beta$ and then $h$ on $(0,a_1)$ so that $\int_0^{a_1}(\beta(x)-\beta_0)\pi(x)\,dx=1$ and $\int_0^{a_1}h(x)\psi(x)\,dx=\beta_0$, then \eqref{eq:nbFHeq} holds for $\lambda=0$, so that $0\in\sigma_F$.

\begin{remark}[independence of the exponent]\label{rem:nbr}
For $T>T_0$, the null controllability of \eqref{eq:nbsys} in time $T$ does not depend on $r$: it holds with initial states in $L^r(0,1)$ and controls in $L^r(R_T)$ for some $r$ if, and only if, it holds for every $r$.
Indeed, by Theorem~\ref{thm:nb}, it holds if, and only if, $\sigma_F=\emptyset$, and $\sigma_F$ does not depend on $r$, by Proposition~\ref{prop:nbFH}.
In fact, one can even show that the space $E$ does not depend on $r$, so that the null controllable states are characterized by the same conditions for every $r$.
\end{remark}

\begin{remark}[two regimes]\label{rem:nbregimes}
Proposition~\ref{prop:nbFH} distinguishes \emph{two regimes}.
If $\beta$ is not constant on the band, then $\sigma_F=\emptyset$, and the controls can essentially be constructed explicitly, without Theorem~\ref{thm:main}, see Remark~\ref{rem:nbexplicit} below.
However, if $\beta$ is equal to a constant on the band, then our proof requires Theorem~\ref{thm:main}, whether the test holds or fails.
\end{remark}

\paragraph{Comparison with the literature.}
Table~\ref{tab:nb} gathers the null controllability results for \eqref{eq:nbsys}.
None of the previous works treats a nonlocal term in the equation.
All of them assume either that the control band touches $x=0$ or that $\beta$ vanishes on a neighborhood of $[0,a_1]$, and all but one are set in the Hilbert space $L^2$.
Theorem~\ref{thm:nb} removes these restrictions at once: it allows the nonlocal term, an arbitrary band, and every space $L^r$, $r\in[1,\infty)$.
In return, our assumption $\beta\not\equiv0$ in $(a_1,a_2)$ excludes a function $\beta$ vanishing on the whole band, which is allowed in \cite{HMT18,Mai19}, and our coefficients are more regular.
The case $r=1$ deserves to be emphasized, since $L^1$ is the natural space for transport equations and for densities \cite{GM74}, and thus for the control.
Moreover, it is not reflexive, so that the classical duality with observability is not directly available.
The only result in $L^1$ we are aware of is \cite[Proposition~A.1]{HMT18}, for $k=0$, $a_1=0$ and $\beta=0$ near $0$.
Finally, when the Fattorini--Hautus test fails, there is no result in the literature, whereas Theorem~\ref{thm:nb} describes the null controllable states \emph{exactly}.

\begin{table}[ht]
\centering
\footnotesize
\renewcommand{\arraystretch}{1.5}
\begin{tabular}{@{}>{\raggedright}p{2.8cm}>{\raggedright}p{2.0cm}l>{\raggedright}p{2.9cm}l>{\raggedright\arraybackslash}p{2.1cm}@{}}
\toprule
Result & Final state & Kernel $k$ & Assumptions & Space for $\varrho$ & Control time\tabularnewline
\midrule
\cite[Thm.~1]{BIM01} & steady states\newline on $(\delta,1)$ & $0$ & $a_1=0$, $\beta\ge0$ & $L^2$ & $1\le T<1+\delta$\tabularnewline
\cite[Thm.~1.1]{HMT18} & steady states\newline $>0$ & $0$ & $a_1=0$, $\beta\ge0$, $\beta=0$ on $(0,a_b)$ & $L^2$ & $T>T_0$\tabularnewline
\cite[Prop.~A.1]{HMT18} & $0$ & $0$ & same as above & $L^1$ & $T>T_0$\tabularnewline
\cite[Thm.~1.1]{Mai19} & $0$ & $0$ & $\beta\ge0$, $\beta=0$ on $(0,a_b)$, $a_b>a_1$ & $L^2$ & $T>a_1+T_0$\tabularnewline
Thm.~\ref{thm:nb} & $0$, from $E_\perp$ & $k/\pi\in W^{1,1}$ & $\beta\in W^{1,1}$, $\beta\not\equiv0$ on $(a_1,a_2)$ & $L^r$, $r\in[1,\infty)$ & $T>T_0$, sharp\tabularnewline
\bottomrule
\end{tabular}
\caption{Null controllability results for \eqref{eq:nbsys}.}
\label{tab:nb}
\end{table}

\subsubsection{Proof of Theorem~\ref{thm:nb}}\label{sec:nbproof}

\begin{figure}[ht]
\centering
\begin{tikzpicture}[x=20cm,y=9cm,font=\small,>=latex]
\fill[gray!20] (0,0.3) rectangle (0.5,0.6);
\fill[pattern=north east lines,pattern color=gray!80] (0.15,0)--(0.5,0)--(0.5,0.35)--cycle;
\fill[pattern=north east lines,pattern color=gray!80] (0.05,0.55)--(0.5,0.55)--(0.5,1)--cycle;
\draw (0,0) rectangle (0.5,1);
\draw[->] (0.5,0)--(0.57,0) node[right] {$t$};
\draw[->] (0,1)--(0,1.07) node[above] {$x$};
\draw[densely dashed] (0,0.35)--(0.5,0.35);
\draw[densely dashed] (0,0.55)--(0.5,0.55);
\foreach \u/\l in {0/0,0.3/a_1,0.35/a_1+\varepsilon,0.55/a_2-\varepsilon,0.6/a_2,1/1} {\draw (0.004,\u)--(-0.004,\u) node[left] {$\l$};}
\draw (0.15,0.006)--(0.15,-0.006) node[below] {$T_1-(a_1+\varepsilon)$};
\draw (0.5,0.006)--(0.5,-0.006) node[below] {$T_1$};
\draw[line width=2.2pt] (0.15,0)--(0.5,0);
\draw[line width=2.2pt] (0.05,0.55)--(0.5,0.55);
\draw[very thick] (0.5,0)--(0.5,1);
\node[right] at (0.5,0.78) {$y(T_1,\cdot)=0$};
\node[right] at (0.5,0.17) {$y(T_1,\cdot)=0$};
\fill (0,0) circle (0.12pt);
\node[anchor=south west,font=\footnotesize] at (0.003,0.008) {$b(0)=\int_0^1\gamma y^1$};
\node[fill=white,inner sep=1.5pt,font=\footnotesize] at (0.385,0.045) {$b(t)=\int_t^{T_1}h(s-t)p(s)\,ds$};
\node[fill=white,inner sep=1.5pt,font=\footnotesize] at (0.32,0.505) {$g(t)=\int_t^{T_1}h(s-t+a_2-\varepsilon)p(s)\,ds$};
\draw[thick,->] (0.165,0)--(0.49,0.325);
\draw[thick,->] (0.15,0.55)--(0.49,0.89);
\node at (0.07,0.17) {$\yout$};
\node at (0.07,0.8) {$\yout$};
\node at (0.07,0.45) {$\yin$};
\node[font=\footnotesize] at (0.3,0.42) {control band, $\beta\not\equiv0$ on $(a_1,a_2)$};
\end{tikzpicture}
\caption{The construction for the proof of Theorem~\ref{thm:nb}.}\label{fig:nbconstruction}
\end{figure}
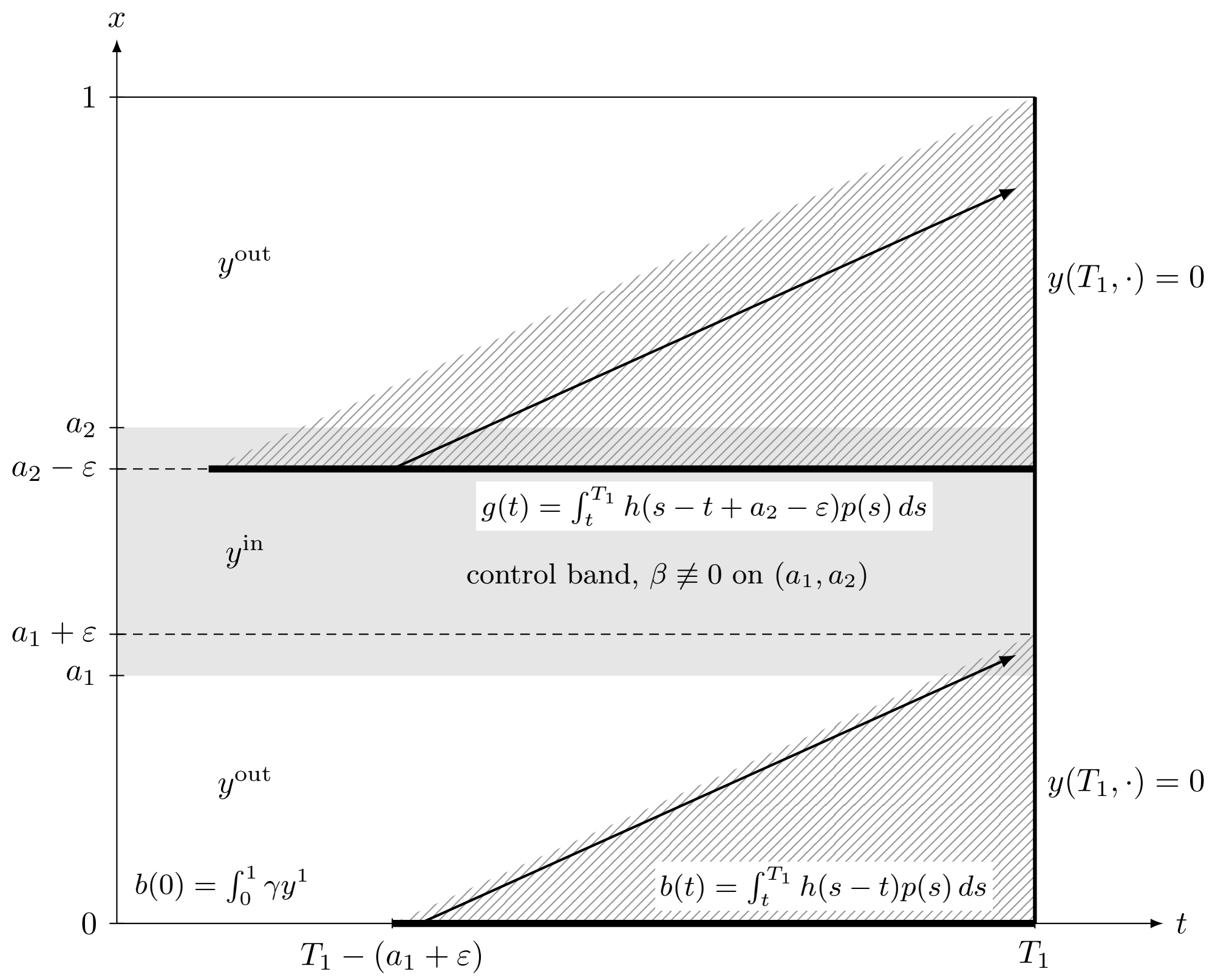

\paragraph{Strategy.}
We apply Theorem~\ref{thm:main}, whose assumption \eqref{eq:primalH} with $F=\ens*{0}$ is the null controllability of a subspace of finite codimension.
We establish it by constructing the controls directly, without going through the adjoint system, so that Theorem~\ref{thm:obs} is not used, and the proof is the same for every $r\in[1,\infty)$, including $r=1$, for which $U_T=L^1(0,T;U)$ is not reflexive, so that an observability estimate would not give \eqref{eq:primalH} through Theorem~\ref{thm:obs}.
The construction is based on the method of \cite[Proposition~3.2]{ABCO17} (see also \cite{HO24-COCV}), which we call \emph{control relocation}, and which goes back to \cite[Theorem~2.2]{AKBGBdT11} in the context of heat equations: we solve \emph{boundary control} problems on the uncontrolled intervals, we let the free solution evolve in the band, and we glue the two with a cut-off function, the control being the residual of the equation in the band.
Here, these boundary control problems are solved along the characteristics.
There are, however, two important differences.
First, the boundary value at $x=0$ is not a free datum but $\int_0^1\gamma y(t,\cdot)$, and we realize it by a \emph{correction in the band}.
Second, the nonlocal term $-hp$, with $p$ given by \eqref{eq:def pn}, couples the uncontrolled intervals to the \emph{whole state}, so that the construction, explicit along the characteristics, is first made for a given function $p$.
The loop is then closed by a Fredholm alternative, which only gives a subspace of finite codimension, and this is where Theorem~\ref{thm:main} is needed.
The construction is summarized in Figure~\ref{fig:nbconstruction}.

\paragraph{Solutions along the characteristics.}
Transport equations are naturally solved along the characteristics, and this yields explicit formulas that are convenient both to construct solutions and to estimate them.
In this paragraph, we define these solutions and collect the properties that we use in the proof.

Let $T>0$, and let $y^0$, $b$ and $f$ be functions on $(0,1)$, $(0,T)$ and $R_T$, respectively.
Then, modulo some regularity, any ``reasonable'' solution to the transport equation
\begin{equation}\label{equ:transp}
\left\{
\begin{aligned}
&y_t+y_x=f(t,x), &&(t,x)\in R_T,\\
&y(t,0)=b(t), &&t\in(0,T),\\
&y(0,x)=y^0(x), &&x\in(0,1),
\end{aligned}
\right.
\end{equation}
is given, after integration along the characteristic $s\mapsto s-t+x$ passing through $(t,x)$, by
\begin{equation}\label{eq:nbchar}
y(t,x)=
\left\{
\begin{aligned}
&y^0(x-t)+\int_0^tf(s,s-t+x)\,ds, &&x>t,\\
&b(t-x)+\int_{t-x}^tf(s,s-t+x)\,ds, &&x<t.
\end{aligned}
\right.
\end{equation}
Whenever the integrals make sense, \eqref{eq:nbchar} defines a function $y$, which we call the \emph{solution along the characteristics} of \eqref{equ:transp}.
Accordingly, we say that $y$ is a solution along the characteristics of \eqref{eq:nbq} if it is the solution along the characteristics of \eqref{equ:transp} with the initial state $y^0$, the source $f=-hp+u\indic_{(a_1,a_2)}$ and the boundary value $b$, where $p$ and $b$ are given by \eqref{eq:def pn} for this function $y$.

\begin{lemma}[characteristic form of the solutions]\label{lem:nbmild}
Let $T>0$, and let $y^0$ and $u$ be such that $\pi y^0\in L^r(0,1)$ and $\pi u\in L^r(R_T)$.
A function $y$ such that $\pi y\in L^r(R_T)$ is the Duhamel solution of \eqref{eq:nbq} if, and only if, it is a solution along the characteristics of \eqref{eq:nbq}.
\end{lemma}

We only sketch the proof, which is standard, except possibly for the treatment of the weight $\pi$.
For smooth data, the Duhamel solution is a classical solution \cite[Corollary~VI.7.6]{EN00}, and integrating it along the characteristics gives \eqref{eq:nbchar}, which then extends to general data by density, by \eqref{eq:nbmono}.
The converse follows from the direct implication and the uniqueness of solutions along the characteristics, which is obtained by multiplying \eqref{eq:nbchar} by $\pi$ and applying the Gronwall lemma to the difference of two of them.

In particular, for every $y^0$ and $u$ as in Lemma~\ref{lem:nbmild}, we also call the Duhamel solution of \eqref{eq:nbq} \emph{the} solution along the characteristics.
When $u=0$, it is $t\mapsto S(t)y^0$, and we keep this notation in what follows.

\begin{lemma}[transport with given data]\label{lem:nbtransport}
Let $T>0$, and let $y$ be given by \eqref{eq:nbchar}, with an initial state $y^0\in L^r(0,1)$, a source $f\in L^r(R_T)$ and a boundary value $b\in L^r(0,T)$.
\begin{enumerate}
\item\label{nbtr:est}
\textbf{Uniform estimate.}
We have
\begin{equation}\label{eq:nbfreebd}
\sup_{t\in[0,T]}\norm{y(t,\cdot)}_{L^r(0,1)}\le\norm{y^0}_{L^r(0,1)}+\norm{b}_{L^r(0,T)}+T^{1/r'}\norm{f}_{L^r(R_T)}.
\end{equation}

\item\label{nbtr:char}
\textbf{Characterization along the characteristics.}
For $(t,x)\in R_T$, let $J_{t,x}:=\ens*{s\in(0,T)\st 0<s-t+x<1}$ be the time interval of the characteristic $s\mapsto s-t+x$ passing through $(t,x)$.
A function $\tilde{y}\in L^r(R_T)$ coincides with $y$ if, and only if, for a.e.\ $(t,x)\in R_T$, the function $s\mapsto \tilde{y}(s,s-t+x)$ coincides, up to a null set, with a function $\psi\in W^{1,1}(J_{t,x})$, with $\psi'(s)=f(s,s-t+x)$, and with $\psi(0)=y^0(x-t)$ if $x>t$, and $\psi(t-x)=b(t-x)$ if $x<t$.

\item\label{nbtr:mean}
\textbf{Regularity of the weighted means.}
For every $w\in W^{1,1}(0,1)$ with $w(1)=0$, the function $m(t):=\int_0^1wy(t,\cdot)$ belongs to $W^{1,r}(0,T)$, with $m(0)=\int_0^1wy^0$ and
$$
\norm{m}_{W^{1,r}(0,T)}\le c\left(\norm{y^0}_{L^r(0,1)}+\norm{f}_{L^r(R_T)}+\norm{b}_{L^r(0,T)}\right),
$$
where $c$ only depends on $\norm*{w}_{W^{1,1}(0,1)}$, $T$ and $r$.
\end{enumerate}
\end{lemma}

We only sketch the proof, which is elementary.
Item~\ref{nbtr:est} follows from \eqref{eq:nbchar} by the triangle and H\"older inequalities and the Fubini theorem.
Item~\ref{nbtr:char} is \eqref{eq:nbchar} written along each characteristic, combined with the Fubini theorem.
For item~\ref{nbtr:mean}, an integration by parts using $w(1)=0$ gives, for smooth data, $m'=\int_0^1wf+\int_0^1w'y+w(0)b$, in which no derivative of $y$ appears.
The estimate then follows from item~\ref{nbtr:est} with the roles of $t$ and $x$ exchanged, and the general case follows by density, since $r<\infty$.

\paragraph{Reduction to initial states in $L^r(0,1)$.}
It is more convenient to work with initial states in the unweighted space $L^r(0,1)\subset Y$.
We reduce to this case by letting the control vanish during an arbitrarily short time, thanks to the following lemma.

\begin{lemma}[the free evolution regularizes]\label{lem:nbreg}
For every $\varepsilon>0$ and $y^0\in Y$, we have $S(\varepsilon)y^0\in L^r(0,1)$.
\end{lemma}

\begin{proof}
Let $y(t,\cdot):=S(t)y^0$, which is the solution along the characteristics with $u=0$, and let $p$ and $b$ be given by \eqref{eq:def pn}.
They are continuous on $[0,\varepsilon]$ (Section~\ref{sec:nbsetting}), so that $f=-hp$ is bounded on $R_\varepsilon$, since $h\in W^{1,1}(0,1)\subset L^\infty(0,1)$.
By \eqref{eq:nbchar} at the time $\varepsilon$, we have $y(\varepsilon,x)=y^0(x-\varepsilon)+\int_0^\varepsilon f(s,s-\varepsilon+x)\,ds$ for $x\in(\varepsilon,1)$, where the first term belongs to $L^r(\varepsilon,1)$, since $\pi(x-\varepsilon)\ge\pi(1-\varepsilon)>0$, and the second one is bounded.
For $x\in(0,\varepsilon)$, we have $y(\varepsilon,x)=b(\varepsilon-x)+\int_{\varepsilon-x}^\varepsilon f(s,s-\varepsilon+x)\,ds$, which is bounded.
\end{proof}

From now on, fix $0<\varepsilon<(a_2-a_1)/2$ and $T_1>T_0+\varepsilon$.
We construct a subspace $G\subset Y$ of finite codimension such that every $y^0\in G$ can be steered to $0$ in time $\varepsilon+T_1$ by a control of $U_{\varepsilon+T_1}$.
This is assumption \eqref{eq:primalH} with $\varepsilon+T_1$ in place of $T_0$ and $F=\ens*{0}$, so that Theorem~\ref{thm:main} gives the result for every $T>\varepsilon+T_1$, hence for every $T>T_0$, since $\varepsilon+T_1$ can be arbitrarily close to $T_0$.

Let $y^0\in Y$, and set
$$
y^1:=S(\varepsilon)y^0\in L^r(0,1).
$$
In Steps~\ref{step:nbnec} to~\ref{step:nbloop}, we construct a control steering $y^1$ to $0$ in time $T_1$, provided that $y^1$ belongs to a subspace $G_1$ of $L^r(0,1)$ of finite codimension, so that the total time is $\varepsilon+T_1$.
In these steps, the data are in the unweighted spaces $L^r$, so that all the solutions along the characteristics belong to $L^r(R_{T_1})$, by Lemma~\ref{lem:nbtransport}, item~\ref{nbtr:est}.

\refstepcounter{proofstep}\label{step:nbnec}
\paragraph{Step~\theproofstep: necessary conditions.}
We start by deriving necessary conditions for $y^1$ to be steered to $0$, which guide the construction.
Let $y$ be the solution along the characteristics of \eqref{eq:nbq} with $y(0,\cdot)=y^1$ and a control $u\in L^r(R_{T_1})$ such that $y(T_1,\cdot)=0$, and let $p$ and $b$ be given by \eqref{eq:def pn}.

\emph{A regularity condition.}
We first note that
\begin{equation}\label{eq:nbbirthsnec}
b\in W^{1,r}(0,T_1),\qquad b(0)=\int_0^1\gamma y^1.
\end{equation}
This follows from Lemma~\ref{lem:nbtransport}.
Indeed, $y$ is given by \eqref{eq:nbchar}, with the initial state $y^1\in L^r(0,1)$, the source $-hp+u\indic_{(a_1,a_2)}\in L^r(R_{T_1})$ and the boundary value $b\in L^r(0,T_1)$.
Moreover, $\gamma=\beta\pi$ belongs to $W^{1,1}(0,1)$ and $\gamma(1)=0$, since $\beta\in W^{1,1}(0,1)$, by \eqref{eq:nbhyp}, and $\pi\in W^{1,1}(0,1)$ with $\pi(1)=0$, by \eqref{eq:nbmu}.

\emph{Conditions outside the band.}
Since $p$ depends on the unknown solution, we regard it from now on as a parameter in $L^r(0,T_1)$, which is determined in Step~\ref{step:nbloop}.
Assume moreover, as a guideline, that the control vanishes outside $[a_1+\varepsilon,a_2-\varepsilon]$, that is, slightly inside the band, to leave room in Step~\ref{step:nbout} for a cut-off function.

\emph{The interval $(0,a_1+\varepsilon)$.}
Let $x\in(0,a_1+\varepsilon)$.
Since $a_1+\varepsilon \leq T_1$, the characteristic $s\mapsto s-T_1+x$ passing through $(T_1,x)$ leaves the boundary $\ens{x=0}$ at the time $T_1-x$ and stays in $(0,a_1+\varepsilon)$, so that
$$
y(T_1,x)=b(T_1-x)-\int_{T_1-x}^{T_1}h(s-T_1+x)p(s)\,ds.
$$
With $t:=T_1-x$, we deduce that $y(T_1,\cdot)=0$ on $(0,a_1+\varepsilon)$ if, and only if,
$$
b(t)=\int_t^{T_1}h(s-t)p(s)\,ds,\qquad T_1-a_1-\varepsilon<t<T_1.
$$
We define $b$ accordingly, and we connect it with a cut-off function to the value $b(0)$ prescribed by \eqref{eq:nbbirthsnec}.
Let $\theta\in C^\infty([0,T_1])$ with $\theta(0)=0$ and $\theta=1$ on $[T_1-a_1-\varepsilon,T_1]$, and set
\begin{equation}\label{equ:def b}
b(t):=\theta(t) \int_t^{T_1}h(s-t)p(s)\,ds +(1-\theta(t))\int_0^1\gamma y^1,\qquad 0<t<T_1.
\end{equation}
Here, if $T_1>1$, the function $h=k/\pi$, which belongs to $W^{1,1}(0,1)$ by assumption, has been extended to a function of $W^{1,1}(0,T_1)$, still denoted by $h$ (for instance by the constant $h(1)$ on $(1,T_1)$), so that $h(s-t)$ is defined for $0<t<s<T_1$.
Using that $h\in W^{1,1}(0,T_1)$, we check that $b\in W^{1,r}(0,T_1)$.

\emph{The interval $(a_2-\varepsilon,1)$.}
We argue similarly, with the boundary value $g(t)=y(t,a_2-\varepsilon)$ of $y$ at $x=a_2-\varepsilon$.
Let $x\in(a_2-\varepsilon,1)$.
Since $1-a_2+\varepsilon\le T_1$, the characteristic $s\mapsto s-T_1+x$ passing through $(T_1,x)$ crosses $\ens{x=a_2-\varepsilon}$ at the time $T_1-x+a_2-\varepsilon$ and stays in $(a_2-\varepsilon,1)$, so that
$$
y(T_1,x)=g(T_1-x+a_2-\varepsilon)-\int_{T_1-x+a_2-\varepsilon}^{T_1}h(s-T_1+x)p(s)\,ds.
$$
With $t:=T_1-x+a_2-\varepsilon$, we deduce that $y(T_1,\cdot)=0$ on $(a_2-\varepsilon,1)$ if, and only if,
$$
g(t)=\int_t^{T_1}h(s-t+a_2-\varepsilon)p(s)\,ds,\qquad T_1-1+a_2-\varepsilon<t<T_1.
$$
We define $g$ accordingly, and we extend it by zero, that is, we set
\begin{equation}\label{equ:def g}
g(t):=
\left\{
\begin{aligned}
&\int_t^{T_1}h(s-t+a_2-\varepsilon)p(s)\,ds, &&T_1-1+a_2-\varepsilon<t<T_1,\\
&0, && \text{ otherwise.}
\end{aligned}
\right.
\end{equation}
We have $g\in L^r(0,T_1)$, since $h\in L^\infty(0,1)$.

\emph{Linearity and continuity.}
It is important for what follows that all the functions we construct are linear and continuous in $(y^1,p)$.
For $b$ and $g$, this follows from their explicit formulas \eqref{equ:def b} and \eqref{equ:def g}, and, more precisely, we have
\begin{equation}\label{eq:nbnbd}
\norm{b}_{W^{1,r}(0,T_1)}+\norm{g}_{L^r(0,T_1)}\les\norm{y^1}_{L^r(0,1)}+\norm{p}_{L^r(0,T_1)}.
\end{equation}

\refstepcounter{proofstep}\label{step:nbout}
\paragraph{Step~\theproofstep: the state outside the band.}
With $b$ and $g$ defined by \eqref{equ:def b} and \eqref{equ:def g}, let $y^c$ be the solution along the characteristics of
$$
\left\{
\begin{aligned}
&y^c_t+y^c_x=-h(x)p(t), &&t\in(0,T_1),\ x\in(0,1)\setminus[a_1+\varepsilon,a_2-\varepsilon],\\
&y^c(t,0)=b(t),\qquad y^c(t,a_2-\varepsilon)=g(t), &&t\in(0,T_1),\\
&y^c(0,x)=y^1(x), &&x\in(0,1)\setminus[a_1+\varepsilon,a_2-\varepsilon].
\end{aligned}
\right.
$$
It satisfies $y^c(T_1,\cdot)=0$, by construction, and we extend it by zero on $(0,T_1)\times[a_1+\varepsilon,a_2-\varepsilon]$.
Let now $\zeta\in C^\infty([0,1])$ with $\zeta=1$ on $[0,a_1+\varepsilon/2]\cup[a_2-\varepsilon/2,1]$ and $\zeta=0$ on $[a_1+\varepsilon,a_2-\varepsilon]$.
Then, using Lemma~\ref{lem:nbtransport}, we can check that the function $\yout:=\zeta y^c$ is the solution along the characteristics of
$$
\left\{
\begin{aligned}
&\yout_t+\yout_x=-\zeta(x)h(x)p(t)+\zeta'(x)y^c(t,x), &&(t,x)\in R_{T_1},\\
&\yout(t,0)=b(t), &&t\in(0,T_1),\\
&\yout(0,x)=\zeta(x)y^1(x), &&x\in(0,1),
\end{aligned}
\right.
$$
and we have $\yout(T_1,\cdot)=0$.

\refstepcounter{proofstep}\label{step:nbin}
\paragraph{Step~\theproofstep: the state inside the band.}
Let $y^f(t,\cdot):=S(t)y^1$ and $p^f(t):=\int_0^1\pi y^f(t,\cdot)$.
Let $\eta\in C^\infty([0,T_1])$ with $\eta(0)=1$ and $\eta(T_1)=0$.
Then, using Lemma~\ref{lem:nbtransport}, we can check that the function $\yin:=(1-\zeta)\eta y^f$ is the solution along the characteristics of
$$
\left\{
\begin{aligned}
&\yin_t+\yin_x=(1-\zeta)\left(\eta'y^f-\eta hp^f\right)-\zeta'\eta y^f, &&(t,x)\in R_{T_1},\\
&\yin(t,0)=0, &&t\in(0,T_1),\\
&\yin(0,x)=(1-\zeta(x))y^1(x), &&x\in(0,1),
\end{aligned}
\right.
$$
and we have $\yin(T_1,\cdot)=0$.

\refstepcounter{proofstep}\label{step:nbglue}
\paragraph{Step~\theproofstep: gluing.}
Define $\hat y:=\yin+\yout$.
By Steps~\ref{step:nbout} and~\ref{step:nbin}, it is the solution along the characteristics of
$$
\left\{
\begin{aligned}
&\hat y_t+\hat y_x=-h(x)p(t)+\hat u(t,x), &&(t,x)\in R_{T_1},\\
&\hat y(t,0)=b(t), &&t\in(0,T_1),\\
&\hat y(0,x)=y^1(x), &&x\in(0,1),
\end{aligned}
\right.
$$
with
$$
\hat u:=(1-\zeta)\left(\eta'y^f-\eta hp^f+hp\right)+\zeta'\left(y^c-\eta y^f\right).
$$
Now observe that $\hat u$ vanishes outside $(a_1+\varepsilon/2,a_2-\varepsilon/2)$, since $\zeta=1$ and $\zeta'=0$ there, so that, in particular, $\hat u=\hat u\indic_{(a_1,a_2)}$.
Moreover, $\hat y(T_1,\cdot)=0$, and $\hat u\in L^r(R_{T_1})$, since $p,p^f\in L^r(0,T_1)$ and $y^f,y^c \in L^r(R_{T_1})$.
However, $\hat y$ is not yet a solution along the characteristics of \eqref{eq:nbq}: the quantities $\int_0^1\gamma\hat y(t,\cdot)$ and $\int_0^1\pi\hat y(t,\cdot)$ have no reason to be equal to $b(t)$ and $p(t)$.

\refstepcounter{proofstep}\label{step:nbcorr}
\paragraph{Step~\theproofstep: correction in the band.}
We now use the assumption that $\gamma=\beta\pi \not\equiv 0$ in $(a_1,a_2)$.
This gives the existence of $\chi\in C_c^\infty(a_1,a_2)$ such that $\int_0^1\gamma\chi=1$, and we then define
$$
\alpha(t):=b(t)-\int_0^1\gamma\hat y(t,\cdot),\qquad y:=\hat y+\alpha\chi .
$$
For the moment, $y$ is the solution along the characteristics of
$$
\left\{
\begin{aligned}
&y_t+y_x=-h(x)p(t)+ u(t,x), &&(t,x)\in R_{T_1},\\
&y(t,0)=\int_0^1\gamma y(t,\cdot), &&t\in(0,T_1),\\
&y(0,x)=y^1(x), &&x\in(0,1),
\end{aligned}
\right.
$$
with $u:=\hat u+\alpha'\chi+\alpha\chi'$.
Indeed, $y(t,0)=\hat y(t,0)=b(t)$, since $\chi$ vanishes near $0$, and $\int_0^1\gamma y(t,\cdot)=b(t)$ by construction.
The equality $y(0,\cdot)=y^1$ comes from $\alpha(0)=b(0)-\int_0^1\gamma y^1=0$ by \eqref{eq:nbbirthsnec}.
Finally, the correction $\alpha\chi$ is supported in the band, and so is its source $(\alpha\chi)_t+(\alpha\chi)_x=\alpha'\chi+\alpha\chi'$, which can thus be produced by the control, since $\alpha\in W^{1,r}(0,T_1)$, by \eqref{eq:nbnbd} and Lemma~\ref{lem:nbtransport}, item~\ref{nbtr:mean}.

We have $y(T_1,\cdot)=0$ since $\alpha(T_1)=b(T_1)-\int_0^1\gamma\hat y(T_1,\cdot)=0$.
Hence, as soon as $p=\int_0^1\pi y$, the function $y$ is the solution along the characteristics of \eqref{eq:nbq}, with the control $u$, and it vanishes at the time $T_1$.

\refstepcounter{proofstep}\label{step:nbloop}
\paragraph{Step~\theproofstep: closing the loop.}
We write the condition $p=\int_0^1\pi y(\cdot\,,\cdot)$ as a fixed point of the map $\Psi:L^r(0,1)\times L^r(0,T_1)\to L^r(0,T_1)$ defined by
$$
\Psi(y^1,p)(t):=\int_0^1\pi y(t,\cdot).
$$
All the functions constructed so far are linear in $(y^1,p)$, so that $\Psi$ is linear, and
$$
\Psi(y^1,p)=Ly^1+Kp,\qquad Ly^1:=\Psi(y^1,0),\qquad Kp:=\Psi(0,p).
$$
Thus, the fixed point condition $p=\Psi(y^1,p)$ reads
$$
(\Id-K)p=Ly^1\quad\text{in }L^r(0,T_1).
$$
Consequently, such a $p$ exists if, and only if, $y^1$ belongs to the space
$$
G_1:=\ens{y^1\in L^r(0,1)\st L y^1\in\Ima(\Id-K)}.
$$
It remains to see that this is a subspace of finite codimension of $L^r(0,1)$.

Writing
$$
\Psi(y^1,p)(t)=\int_0^1\pi\hat y(t,\cdot)+b(t)\int_0^1\pi\chi -\int_0^1\gamma\hat y(t,\cdot)
\int_0^1\pi\chi,
$$
and using Lemma~\ref{lem:nbtransport} together with the estimates \eqref{eq:nbfreebd} and \eqref{eq:nbnbd}, we obtain
$$
\norm{Kp}_{W^{1,r}(0,T_1)}\les \norm{p}_{L^r(0,T_1)}.
$$
Thus, $K$ is bounded from $L^r(0,T_1)$ into $W^{1,r}(0,T_1)$, and, since the embedding of $W^{1,r}(0,T_1)$ into $L^r(0,T_1)$ is compact \cite[Theorem~8.8]{Bre11}, $K$ is a compact operator on $L^r(0,T_1)$.
By the Fredholm alternative \cite[Theorem~6.6]{Bre11}, $\Ima(\Id-K)$ has finite codimension in $L^r(0,T_1)$, since $\Ima(\Id-K)=\left(\ker(\Id-K^*)\right)_\perp$ and $\ker(\Id-K^*)$ is finite dimensional.
Hence $G_1=L^{-1}\left(\Ima(\Id-K)\right)$ has finite codimension in $L^r(0,1)$.

\refstepcounter{proofstep}\label{step:nbconc}
\paragraph{Step~\theproofstep: conclusion.}
We can now put the previous steps together.
Let $G$ be the subspace of all $y^0\in Y$ such that $S(\varepsilon)y^0\in G_1$.
It has finite codimension in $Y$, since $G=S(\varepsilon)^{-1}(G_1)$, $G_1$ is a subspace of finite codimension of $L^r(0,1)$, and $S(\varepsilon)$ is a linear map from $Y$ into $L^r(0,1)$ by Lemma~\ref{lem:nbreg}.
For $y^0\in G$, introduce $y^1:=S(\varepsilon)y^0\in G_1$, let $p \in L^r(0,T_1)$ solve $(\Id-K)p=Ly^1$, and let $y$ and $u$ be the corresponding functions of Step~\ref{step:nbcorr}.
We shift them by $\varepsilon$ and extend them to $(0,\varepsilon)$ by $y(t,\cdot):=S(t)y^0$ and $u(t,\cdot):=0$, still denoted by $y$ and $u$.
We have $u\in L^r(R_{\varepsilon+T_1})\subset U_{\varepsilon+T_1}$.
Since $y(\varepsilon,\cdot)=y^1$, splitting the integrals in \eqref{eq:nbchar} at $s=\varepsilon$ shows that $y$ is the solution along the characteristics of \eqref{eq:nbq} on $(0,\varepsilon+T_1)$, with the initial state $y^0$ and the control $u$, and it satisfies $y(\varepsilon+T_1,\cdot)=0$.
Finally, since $\pi y\in L^r(R_{\varepsilon+T_1})$, $y$ is also the Duhamel solution of \eqref{eq:nbq}, by Lemma~\ref{lem:nbmild}.
This completes the proof of Theorem~\ref{thm:nb}. \qed

\begin{remark}[the two regimes in the proof]\label{rem:nbexplicit}
A correction $\chi\in C_c^\infty(a_1,a_2)$ in the band modifies the pair $(b,p)$ of \eqref{eq:def pn} by $(\int_0^1\gamma\chi,\int_0^1\pi\chi)$.
Thus, whether $b$ and $p$ can be modified independently amounts to whether $\gamma$ and $\pi$ are linearly independent on $(a_1,a_2)$, that is, since $\gamma=\beta\pi$, whether $\beta$ is not constant there.

\emph{If $\beta$ is not constant on the band.}
Then the map $\chi\mapsto(\int_0^1\gamma\chi,\int_0^1\pi\chi)$ is onto $\C^2$.
Otherwise, $c_1\gamma+c_2\pi=0$ on $(a_1,a_2)$ for some $(c_1,c_2)\neq0$, so that $c_1\neq0$, since $\pi>0$, and $\beta=-c_2/c_1$ there.
We can thus take $\chi_2\in C_c^\infty(a_1,a_2)$ with $\int_0^1\gamma\chi_2=0$ and $\int_0^1\pi\chi_2=1$, which prescribes $p$ without modifying $b$.
Precisely, let $p\in C^\infty([0,T_1])$ with $p(0)=\int_0^1\pi y^1$ and $p(T_1)=0$, and let $y$ be as in Step~\ref{step:nbcorr}.
Let $\alpha_2:=p-\Psi(y^1,p)$, with $\Psi$ as in Step~\ref{step:nbloop}, which is the gap between the prescribed $p$ and the actual $\int_0^1\pi y$.
By Lemma~\ref{lem:nbtransport}, the function $\alpha_2$ belongs to $W^{1,r}(0,T_1)$, and it vanishes at $0$ and $T_1$.
Then, as in Step~\ref{step:nbcorr}, $y+\alpha_2\chi_2$ is a solution along the characteristics of \eqref{eq:nbq} with the control $u+\alpha_2'\chi_2+\alpha_2\chi_2'$, and it vanishes at the time $T_1$.
Hence $G_1=L^r(0,1)$, and neither Step~\ref{step:nbloop} nor Theorem~\ref{thm:main} is needed.

\emph{If $\beta$ is equal to a constant $\beta_0$ on the band.}
Then $\int_0^1\gamma\chi=\beta_0\int_0^1\pi\chi$ for every $\chi\in C_c^\infty(a_1,a_2)$, so that a correction in the band cannot modify $p$ without modifying $b$, and the loop has to be closed through Step~\ref{step:nbloop}, hence Theorem~\ref{thm:main}.
\end{remark}

\subsubsection{Proofs of the other results}\label{sec:nbother}

This section gathers the proofs postponed.

\paragraph{Proof of Proposition~\ref{prop:nbsg}.}
Let us write $A=A_0+A_1$ with $A_0y:=-y'$, $D(A_0):=D(A)$ and $A_1y:=-h\int_0^1\pi y$.
The operator $A_1$ is bounded on $Y$: $\norm*{A_1y}_{Y}=\norm*{k}_{L^r(0,1)}\abs*{\int_0^1\pi y}\le\norm*{k}_{L^r(0,1)}\norm*{y}_{Y}$.
Thus, by the perturbation theory \cite[Theorem~III.1.3]{EN00}, it suffices to show that $A_0$ generates a $C_0$-semigroup on $Y$.
For item~\ref{nbsg:gen}, we apply the Lumer--Phillips theorem \cite[Theorem~II.3.15]{EN00} to $A_0-\omega$, for some $\omega>0$ to be determined below.

\emph{Density of the domain.}
The space $C_c^\infty(0,1)$ is dense in $Y$, since $\pi$ is bounded from above and from below on compact subsets of $(0,1)$.
Let $\varphi\in C_c^\infty(0,1)$ and, for $\varepsilon\in(0,1)$ with $\varepsilon\norm*{\beta}_{L^\infty(0,1)}<1$, let $\chi_\varepsilon\in C^\infty([0,1])$ with $\chi_\varepsilon(0)=1$, $\abs*{\chi_\varepsilon}\le1$ and $\chi_\varepsilon=0$ on $[\varepsilon,1]$.
Then $\abs*{\int_0^1\gamma\chi_\varepsilon}\le\varepsilon\norm*{\beta}_{L^\infty(0,1)}<1$, and $y:=\varphi+c\chi_\varepsilon$, with $c:=\int_0^1\gamma\varphi/(1-\int_0^1\gamma\chi_\varepsilon)$, satisfies $y(0)=c=\int_0^1\gamma y$, so that $y\in D(A_0)$, and $\norm*{y-\varphi}_{Y}\le\abs*{c}\varepsilon^{1/r}\to0$ as $\varepsilon\to0$.

\emph{Quasi-dissipativity.}
Let $\lambda>0$, $y\in D(A_0)\setminus\ens*{0}$, $f:=(\lambda-A_0)y$ and $\varrho:=\pi y$, which belongs to $W^{1,1}(0,1-\delta)$ for every $\delta\in(0,1)$ and satisfies $\varrho'+(\lambda+\mu)\varrho=\pi f$ on $(0,1)$ and $\varrho(0)=\int_0^1\beta\varrho$.
Let $\varrho^*:=|\varrho|^{r-2}\conj{\varrho}$, with $\varrho^*:=0$ on $\ens*{\varrho=0}$, so that $\varrho^*\varrho=|\varrho|^r$ and $\norm*{\varrho^*}_{L^{r'}(0,1)}\le\norm*{\varrho}_{L^r(0,1)}^{r-1}$.
By the H\"older inequality,
$$
\norm{\varrho}_{L^r(0,1)}^{r-1}\norm{\pi f}_{L^r(0,1)}\ge\Re\int_0^1\varrho^*\pi f\,dx=\Re\int_0^1\varrho^*\left(\varrho'+\mu\varrho\right)dx+\lambda\norm{\varrho}_{L^r(0,1)}^r ,
$$
where $\varrho^*(\varrho'+\mu\varrho)$ is integrable, since $\varrho'+\mu\varrho=\pi f-\lambda\varrho\in L^r(0,1)$.
Using that $|\varrho|^r \in W^{1,1}(0,1-\delta)$ with $(|\varrho|^r)'=r\Re(\varrho^*\varrho')$ and $\mu\ge0$, we get
$$
\Re\int_0^{1-\delta}\varrho^*\left(\varrho'+\mu\varrho\right)dx\ge-\frac1r\abs{\varrho(0)}^r
=-\frac1r \abs{\int_0^1\beta\varrho\,dx}^r
\ge-\omega\norm{\varrho}_{L^r(0,1)}^r ,
$$
with $\omega:=\norm*{\beta}_{L^\infty(0,1)}^r/r$.
Letting $\delta\to0$, combining the two previous estimates and dividing by $\norm*{\varrho}_{L^r(0,1)}^{r-1}\ne0$, we get $\norm*{\pi f}_{L^r(0,1)}\ge(\lambda-\omega)\norm*{\varrho}_{L^r(0,1)}$, that is, $\norm*{(\lambda-A_0)y}_{Y}\ge(\lambda-\omega)\norm*{y}_{Y}$, so that $A_0-\omega$ is dissipative.

\emph{The range condition.}
It remains to show that $\lambda-A_0$ is onto for some $\lambda>\omega$.
Let $\lambda>\omega$ and $f\in Y$.
The equation $(\lambda-A_0)y=f$ means that $y'+\lambda y=f$ on $(0,1)$ and $y(0)=\int_0^1\gamma y$, that is,
$$
y(x)=y(0)e^{-\lambda x}+\int_0^xe^{-\lambda(x-\xi)}f(\xi)\,d\xi,\qquad \left(1-\hat\gamma(\lambda)\right)y(0)=\int_0^1\gamma(x)\int_0^xe^{-\lambda(x-\xi)}f(\xi)\,d\xi\,dx,
$$
where $\hat\gamma(\lambda):=\int_0^1\gamma(x)e^{-\lambda x}\,dx$.
Assume moreover that $\lambda>\norm*{\beta}_{L^\infty(0,1)}$.
Then, $|\hat\gamma(\lambda)|\le\norm*{\beta}_{L^\infty(0,1)}/\lambda<1$, so that this determines $y(0)$, hence $y$, which we denote by $R(\lambda)f$.
It remains to check that $y=R(\lambda)f$ belongs to $D(A_0)$.
By \eqref{eq:nbmono} and since $\lambda \geq 0$, we have $|\pi(x)\int_0^xe^{-\lambda(x-\xi)}f(\xi)\,d\xi|\le\int_0^x |\pi f|(\xi)\,d\xi$, so that $y\in Y$ and $y'=f-\lambda y\in Y$, while $y(0)=\int_0^1\gamma y$ by construction.
Hence $\lambda-A_0$ is bijective, with inverse $R(\lambda)$, which is bounded, by the quasi-dissipativity of $A_0$.
In particular, $A_0$ is closed, and the Lumer--Phillips theorem applies and shows that $A_0$ generates a $C_0$-semigroup.

\emph{The adjoint.}
We recall that $Y'$ is identified with the space of the functions $z$ such that $z/\pi\in L^{r'}(0,1)$, through $\pair*{z}{y}_{Y',Y}=\int_0^1z\,\conj{y}\,dx$, and that the adjoints are computed in this duality.
Since $A_1$ is bounded, with $A_1^*z=-(\int_0^1hz)\pi$, we have $A^*=A_0^*+A_1^*$, so that it suffices to show that $D(A_0^*)=D$, where $D:=\ens{z\in Y'\st z'\in Y',\ z(1)=0}$, and that $A_0^*z=z'+z(0)\gamma$ for $z\in D$.
Let $\lambda$ as in the previous paragraph, so that $(\lambda-A_0^*)^{-1}=R(\lambda)^*$ and $D(A_0^*)$ is thus the range of $R(\lambda)^*$, which we now compute.
Let $w\in Y'$.
Inserting the formula for $R(\lambda)f$ into $\pair*{w}{R(\lambda)f}_{Y',Y}$ and exchanging the order of integration, which is allowed by the Fubini theorem, since $|w(x)f(\xi)|\le|w/\pi|(x)|\pi f|(\xi)$ for $\xi\le x$, by \eqref{eq:nbmono}, and similarly with $\gamma$ in place of $w$, we get $\pair*{w}{R(\lambda)f}_{Y',Y}=\pair*{z_w}{f}_{Y',Y}$ for every $f\in Y$, that is, $R(\lambda)^*w=z_w$, where
$$
z_w(x):=\int_x^1e^{-\lambda(\xi-x)}\left(w(\xi)+c_w\gamma(\xi)\right)d\xi,\qquad c_w:=\frac{1}{1-\hat\gamma(\lambda)}\int_0^1e^{-\lambda\xi}w(\xi)\,d\xi .
$$
By \eqref{eq:nbmono} and since $\lambda \geq 0$, we have $|z_w(x)/\pi(x)|\le\norm*{w/\pi}_{L^1(0,1)}+|c_w|\norm*{\beta}_{L^1(0,1)}$, so that $z_w/\pi\in L^\infty(0,1)$.
It follows that $z_w'/\pi=\lambda z_w/\pi-w/\pi-c_w\beta \in L^{r'}(0,1)$.
Thus $z_w\in D$ and $A_0^*z_w=\lambda z_w-w=z_w'+z_w(0)\gamma$.
It remains to show that every $z\in D$ is of the form $z_w$, and the only candidate is $w:=\lambda z-z'-z(0)\gamma$, which clearly belongs to $Y'$.
Let $\zeta:=z-z_w\in D$.
It solves $\zeta'=\lambda\zeta-\zeta(0)\gamma$, $\zeta(1)=0$, so that $\zeta(x)=\zeta(0)\int_x^1e^{-\lambda(\xi-x)}\gamma(\xi)\,d\xi$, and thus $\zeta(0)=\zeta(0)\hat\gamma(\lambda)$, that is, $\zeta(0)=0$, since $\hat\gamma(\lambda)\ne1$, hence $\zeta=0$, that is, $z=z_w$. \qed

\paragraph{Proof of Proposition~\ref{prop:nbfirst}, item~\ref{nbfirst:nonzero}.}
If $Ay=\lambda y$ with $\lambda\ge0$ and $y\in D(A)\setminus\ens*{0}$, then $S(t)y=e^{\lambda t}y$ for every $t\ge0$, since the derivative of $t\mapsto e^{-\lambda t}S(t)y$ is $e^{-\lambda t}S(t)(A-\lambda)y=0$, so that $S(t)$ does not tend to $0$.
Let $\lambda\ge0$, and let $y\in D(A)$ with $Ay=\lambda y$.
Setting $c:=y(0)$ and $p:=\int_0^1\pi y$, this means that $y'+\lambda y=-hp$ on $(0,1)$, with $c=\int_0^1\gamma y$.
Since $y\in W^{1,r}(0,1-\delta)$ for every $\delta\in(0,1)$, solving this differential equation with $y(0)=c$ gives
$$
y(x)=ce^{-\lambda x}-pH_\lambda(x),\qquad H_\lambda(x):=\int_0^xe^{-\lambda(x-\xi)}h(\xi)\,d\xi,
$$
for $x\in[0,1)$.
Conversely, for every $(c,p)\in\C^2$, this function satisfies $y'=-\lambda y-hp$, and $\pi y$ and $\pi y'$ are bounded, since $\pi(x)|h(\xi)|\le|k(\xi)|$ for $\xi\le x$, by \eqref{eq:nbmono}.
Hence $y\in D(A)$ and $Ay=\lambda y$ as soon as $c=\int_0^1\gamma y$ and $p=\int_0^1\pi y$.
These two conditions form a homogeneous linear system in $(c,p)$, whose determinant is
$$
D(\lambda):=\det
\begin{pmatrix}
1-\int_0^1\gamma(x)e^{-\lambda x}\,dx & \int_0^1\gamma H_\lambda\\[1ex]
-\int_0^1\pi(x)e^{-\lambda x}\,dx & 1+\int_0^1\pi H_\lambda
\end{pmatrix}.
$$
If $D(\lambda)=0$, this system has a solution $(c,p)\ne(0,0)$, and then the corresponding function $y$ is nonzero.
The function $D$ is real valued, since the coefficients are real.
Since $H_\lambda(x)\to0$ as $\lambda\to+\infty$ for every $x\in(0,1)$ and $\pi|H_\lambda|\le\norm*{k}_{L^1(0,1)}$, the dominated convergence theorem shows that $D$ is continuous on $[0,\infty)$ and tends to $1$ at infinity.
The assumption of item~\ref{nbfirst:nonzero} says exactly that $D(0)\le0$, so that, by the intermediate value theorem, $D(\lambda)=0$ for some $\lambda\ge0$, and then $A$ has the eigenvalue $\lambda$, by the above. \qed

\paragraph{Proof of Proposition~\ref{prop:nbfirst}, item~\ref{nbfirst:EC}.}
Let $T>0$ and $u\in U_T$, let $y$ be the solution of \eqref{eq:nbq} with $y^0=0$, and let $f:=-hp+u\indic_{(a_1,a_2)}$, where $h$ is bounded, by \eqref{eq:nbhyp}.
If $a_1>0$, then, for $x<\min\ens{a_1,T}$, \eqref{eq:nbchar} gives $y(T,x)=b(T-x)-\int_{T-x}^Th(s-T+x)p(s)\,ds$, which is bounded, so that $x\mapsto x^{-1/(2r)}$, which belongs to $Y$, is not reachable.
If $a_2<1$, then $f\in L^r(R_T)$, since $\pi\ge\pi(a_2)>0$ on $(a_1,a_2)$, and \eqref{eq:nbchar} gives $y(T,x)=b(T-x)\indic_{\{x<T\}}+\int_{(T-x)_+}^Tf(s,s-T+x)\,ds$, so that $y(T,\cdot)\in L^r(0,1)$, by the H\"older inequality.
But $\mu^{1/r}\in Y\setminus L^r(0,1)$, since $\int_0^1\pi^r\mu\le\int_0^1\pi\mu=1$ and $\int_0^1\mu=\infty$, so that $\mu^{1/r}$ is not reachable. \qed

\paragraph{Proof of Proposition~\ref{prop:nbfirst}, item~\ref{nbfirst:T0}.}
Let $T<a_1$, the case $T<1-a_2$ being similar, let $y^0\in Y$ and $u\in U_T$, let $y$ be the solution of \eqref{eq:nbq}, and let $p$ be given by \eqref{eq:def pn}.
For $x\in(T,a_1)$, the characteristic through $(T,x)$ starts at the time $0$ from $x-T\in(0,a_1)$, and it stays in $(0,a_1)$, so that it never meets the band: $y(T,x)=y^0(x-T)-\int_0^Th(s-T+x)p(s)\,ds$ by \eqref{eq:nbchar}.
The integral is bounded, as in the proof of item~\ref{nbfirst:EC}.
Thus, if $y(T,\cdot)=0$, then $y^0$ is bounded near $\xi=0$, so that $y^0(\xi):=\xi^{-1/(2r)}$, which belongs to $Y$, is not a null controllable state in time $T$. \qed

\paragraph{Proof of Proposition~\ref{prop:nbFH}.}
Let $\lambda\in\C$ and $z\in D(A^*)$ with $A^*z=\lambda z$ and $B^*z=0$, that is, $z=0$ on $(a_1,a_2)$, and let $m:=\int_0^1hz$.
By \eqref{eq:nbadj}, the eigenvalue equation reads
\begin{equation}\label{eq:nbeig}
z'=\lambda z-z(0)\gamma+m\pi\quad\text{on }(0,1).
\end{equation}
On $(a_1,a_2)$, we have $z'=0$, so that $z(0)\gamma=m\pi$, that is, $z(0)\beta=m$.

\emph{Item~\ref{nbFH:indep}.}
If $\beta$ is not constant on $(a_1,a_2)$, this gives $z(0)=m=0$, so that \eqref{eq:nbeig} gives $z=0$.

\emph{Item~\ref{nbFH:dep}.}
We get $m=\beta_0z(0)$, and \eqref{eq:nbeig} becomes $z'=\lambda z-z(0)(\beta-\beta_0)\pi$ on $(0,1)$.
If $z(0)=0$, then $z=0$ as above.
Thus, if $z\ne0$, we may normalize $z(0)=1$, so that $z=\psi$ and $m=\beta_0$.
Since $z$ vanishes on $(a_1,a_2)$ and $z(1)=0$, this proves that $\lambda\in\sigma_F$ implies \eqref{eq:nbFHeq}, and that $\ker(\lambda-A^*)\cap\ker B^*=\Span\ens*{\psi}$.

Conversely, assume that $\psi$ satisfies \eqref{eq:nbFHeq}.
Since $\psi'=\lambda\psi-(\beta-\beta_0)\pi$, $\psi(a_1)=0$ and $\beta=\beta_0$ on $(a_1,a_2)$, the function $\psi$ vanishes on $[a_1,a_2]$, so that $B^*\psi=0$.
Solving the differential equation on $(x,1)$ with $\psi(1)=0$, we get $\psi(x)=e^{\lambda x}\int_x^1e^{-\lambda\xi}(\beta(\xi)-\beta_0)\pi(\xi)\,d\xi$, and we deduce $|\psi(x)/\pi(x)|\le e^{|\lambda|}\norm*{\beta-\beta_0}_{L^\infty(0,1)}$, by \eqref{eq:nbmono}.
Hence $\psi/\pi$ is bounded, and so is $\psi'/\pi=\lambda\psi/\pi-(\beta-\beta_0)$, so that $\psi\in D(A^*)$.
Since $\psi(0)=1$ and $\int_0^1h\psi=\beta_0$, we get $A^*\psi=\psi'+\gamma-\beta_0\pi=\psi'+(\beta-\beta_0)\pi=\lambda\psi$, and $\psi\ne0$, that is, $\lambda\in\sigma_F$. \qed


\section*{Acknowledgments}

This project was supported by National Science Centre, Poland UMO-2023/50/E/ST1/00081.


\renewcommand{\bibfont}{\small}
\bibliographystyle{amsalpha}
\bibliography{biblio}

\end{document}